%% file: article_DEN.tex
\documentclass[12pt]{article}
\usepackage{amsfonts}
\usepackage{amsmath}
\usepackage{verbatim}
\usepackage{amssymb}
\usepackage[latin1]{inputenc}
\usepackage{a4wide}
\usepackage{amsthm}
\usepackage{graphicx}

\DeclareMathAlphabet{\mathpzc}{OT1}{pzc}{m}{it}

\newcommand{\N}{\ensuremath{\mathbb N}}

\newcommand{\R}{\ensuremath{\mathbb R}}

\newcommand{\T}{\ensuremath{\mathbb T}}
\newcommand{\Y}{\ensuremath{\mathbb Y}}
\newcommand{\Pp}{\ensuremath{\mathbb P}}
\newcommand{\1}{\mathbf{1}}

\newcounter{comptage}[part]
\newtheorem{lem}[comptage]{Lemma}
\newtheorem{theo}[comptage]{Theorem}
\newtheorem{defin}[comptage]{Definition}

\newtheoremstyle{remarque}
  {3pt}
  {3pt}
  {}
  {}
  {\bf}
  {.}
  {.5em}
  {}
\theoremstyle{remarque}

\author{\vspace{1cm} Antoine Lemenant \\
Université Paris XI \\
antoine.lemenant@math.u-psud.fr}

\title{Energy improvement for energy minimizing functions in the complement of generalized Reifenberg-flat sets.}

\begin{document}
\maketitle

  {\bf Abstract.} Let $P$ be an hyperplane in $\R^N$, and denote by $d_H$ the Hausdorff distance.
    We show that for all positive radius $r<1$ there is an $\varepsilon>0$, such
  that if $K$ is a Reifenberg-flat set in $B(0,1)\subset
  \R^N$ that contains the origin, with
  $d_H(K,P)\leq \varepsilon$, and if $u$ is an energy minimizing
  function in $B(0,1)\backslash K$ with restricted values on $\partial B(0,1)\backslash
  K$, then the  normalized energy of $u$ in $B(0,r)\backslash K$ is bounded by the
  normalized energy of  $u$ in $B(0,1)\backslash K$. We also prove the same result in  $\R^3$ when $K$
  is a $\varepsilon$-minimal set, that is a generalization of Reifenberg-flat sets
  with minimal cones of type $\Y$ and $\T$. Moreover, the result is still true for a further generalization of
  sets called  $(\varepsilon,\varepsilon_0)$-minimal. This article is a preliminary study
  for a forthcoming paper where a regularity result for the
  singular set of the Mumford-Shah functional close to minimal cones in $\R^3$ is proved by the same author.




\section*{Introduction}

Let $B$ be the unit ball in $\R^N$ and let $K$ be a closed set
with locally finite $H^{N-1}$ measure. We say that $u$ locally
minimizes the energy in $B \backslash K$ if for every ball
$B(x,r)$ included in $\bar B$ and every function $v \in W^{1,2}(B
\backslash K)$ such that $v=u$ in $B  \backslash (B(x,r)\cup K)$
we have
$$\int_{B(x,r)\backslash K}|\nabla u(x)|^2dx\leq \int_{B(x,r)\backslash K}|\nabla v(x)|^2dx.$$
Then we know that $u$ is harmonic in $B\backslash K$, and morally
has its normal derivative equal to $0$ on $K$, at least if $K$ is
smooth enough.

If $u$ is an energy minimizer in $B\backslash K$, we denote by
$$\omega_2(x,r):=\frac{1}{r^{N-1}}\int_{B(x,r)\backslash K}|\nabla u(x)|^2dx$$
the normalized energy of $u$ in $B(x,r)$. It is well known (see
for instance exercice 7.6 of  \cite{afp}) that if $K=\emptyset$,
then for a locally energy minimizing function in $B$ (that means
that $u$ is harmonic in $B$) we have for all $r>0$,
\begin{eqnarray}
\omega_2(0,r)\leq r^{\gamma}\omega_2(0,1) \label{dect}
\end{eqnarray}
with $\gamma=1$. It is a consequence of the proof of the mean
value inequality for subharmonic functions, applied to $|\nabla
u|^2$. By a reflection argument, this is also true if $K$ is a
hyperplane in $\R^N$ containing $0$.  In this this paper, we want
to find some conditions on $K$ that imply \eqref{dect} for all
energy minimizers in $B\backslash K$ and for some positive
exponent $\gamma$. For instance, we want to prove that
\eqref{dect} is true if $K$ is flat enough, or close enough to a
minimal cone. In fact we would not obtain exactly \eqref{dect},
but we will give some conditions on $K$ for which we know that any
energy minimiser in $B(0,1)\backslash K$ has its normalized energy
smaller in a smaller ball centered at the origin. This will be
enough in many cases because if the set is almost flat, then we
could apply the same result in $B(0,r)$ and do an iteration.

Notice that we cannot expect to have \eqref{dect} when $K$ is any
set at distance less than $\varepsilon$ to a hyperplane. Indeed,
consider a little tube of size $\varepsilon$ in the unit disc of
$\R^2$,
$$K:=\{ (x,y) \in \R^2; y=\pm\varepsilon \}\cap B(0,1)$$
and take $u_0$ a function on the unit circle, that is equal to
zero everywhere, except at one side of the tube where $u_0$ is
equal to a constant $M$.

\begin{center}
\input{cex.tex}
\end{center}

If we minimize the Dirichlet integral over all the functions that
are equal to $u_0$ at the boundary, we probably get a linear
function that goes from $M$ to $0$ in the tube, and which is equal
to $0$ anywhere else. For this minimizer $u$ and for all radius
$r<1$ we have that the normalized energy is almost constant
$$\frac{1}{r}\int_{B(0,r)}|\nabla u|^2\simeq {\frac{\varepsilon}{M^2}}$$
thus we could find a constant $M$, depending on $\varepsilon$ and
$\gamma$, in such a way that $\omega_2(0,r)$ will  obviously not
be bounded by a power of any radius between $0<r<1$.

Therefore, in order to get some decreasing of energy in smaller
balls, we have to make sure that the set $K$ does not contain some
little tubes that could carry the energy from outside to inside.

The first class of sets for which the decreasing of energy will be
true, is the class of Reifenberg-flat sets with small constant.
Let us give some definitions.We denote by $D_{x,r}$ the normalized
Hausdorff distance between two closed sets $E$ and $F$ in $B(x,r)$
defined by
\begin{eqnarray}
D_{x,r}(E,F):=\frac{1}{r}\Big\{\max\{\sup_{y\in E\cap
B(x,r)}d(y,F),\sup_{y\in F\cap B(x,r)}d(y,E)\}\Big\}.
\label{hausdorffdist}
\end{eqnarray}

\begin{defin}\label{epsmin} Let $B$ be a ball in  $\R^3$. A closed set $E\subset
B$, containing the origin is said to be
$\varepsilon_0$-Reifenberg-flat in $B$ if for all $x\in E$ and for
all $r$ such that $B(x,r)\subset B$ we have that
\begin{eqnarray}
\inf_{P \ni x}D_{x,r}(E,P)\leq \varepsilon_0 \label{infh}
\end{eqnarray}
where the infimum is taken over all hyperplanes $P$ that contains
$x$.
\end{defin}

Reifenberg-flat sets are introduced in \cite{re}, where a
regularity theorem is stated: it is proved that every
$\varepsilon_0$-Reifenberg-flat set in $B(0,1)$ with
$\varepsilon_0$ small enough, is the bi-Hölderian image of the
unit disc.

A first result of this paper is the following.

\begin{theo}\label{theorem1} For all  $\gamma <1$ and
$0<r<\frac{1}{2}$, there is an  $\varepsilon_1 >0$ such that for
every $10^{-6}$-Reifenberg-flat set $K$ in the unit ball of $\R^N$
with $P$ a hyperplane through the origin satisfying
$$\sup\{d(y,P); y \in K\cap B(0,1)\}\leq \varepsilon_1,$$
we have that
\begin{eqnarray}
\omega_2(0,r)\leq r^{\gamma}\omega_2(0,1) \label{inl1}
\end{eqnarray}
for all locally energy minimizing function in $B(0,1)\backslash
K$.
\end{theo}

Note that we don't require to have \eqref{inl1} for all $r$ but
just to have  \eqref{inl1}  for a given $r$ if $\varepsilon_2$ is
small enough, depending on $r$.

In fact, we will concentrate on dimension 3 because we are
principally interested in a same statement of Theorem
\ref{theorem1} but for $\varepsilon_0$-minimal sets that will be
defined just after. Thus the proof will be done only in dimension
3 and with more general sets than Reifenberg-flat sets. However,
one can easily see that for the case of Reifenberg-flat sets the
same proof holds in any dimension so Theorem \ref{theorem1} is
also true.

Before defining $\varepsilon_0$-minimal sets, we have to define
the 3 minimal cones in $\R^3$. Cones of type $1$ are planes in
$\R^3$, also called $\Pp$. Cones of types $2$ and $3$ and their
spines are defined as in  \cite{ddpt} by the following way.

\begin{defin}\label{prop} Define $Prop\subset \R^2$
by
$$Prop=\{(x_1,x_2);x_1 \geq 0, x_2=0\} $$
$$\hspace{4cm} \cup\{(x_1,x_2);x_1 \leq 0, x_2=-\sqrt{3}x_1\}$$
$$\hspace{6.5cm}\cup\{(x_1,x_2);x_1 \leq 0, x_2=\sqrt{3}x_1\}.$$
Then let $Y_0=Prop\times \R \subset \R^3.$ The spine of $Y_0$ is
the line $L_0=\{x_1=x_2=0\}$. A cone of type $2$ (or of type $\Y$)
is a set $Y=R(Y_0)$ where  $R$  is the composition of a
 translation and a rotation. The spine of $Y$ is then the line $R(L_0)$.
 We denote by $\Y$ the set of all the cones of type 2. Sometimes we also may use
 the expression ``of type $\Y$".
\end{defin}

\begin{defin}\label{defT} Let $A_1=(1,0,0)$, $A_2=(-\frac{1}{3},\frac{2\sqrt{2}}{3},0)$,
$A_3=(-\frac{1}{3},-\frac{\sqrt{2}}{3},\frac{\sqrt{6}}{3})$, and
$A_4=(-\frac{1}{3},-\frac{\sqrt{2}}{3}, -\frac{\sqrt{6}}{3})$ the
four vertices of a regular tetrahedron centered at $0$. Let $T_0$
be the cone over the union of the $6$ edges $[A_i,A_j]$ $i\not
=j$. The spine of $T_0$  is the union of the four half lines
$[0,A_j[$. A cone of type $3$ (or of type $\T$) is a set
$T=R(T_0)$ where $R$ is the composition of a translation and a
rotation. The spine of $T$ is the image by $R$  of the spine of
$T_0$. We denote by $\T$ the set of all the cones of type 3.
\end{defin}

\begin{center}
\includegraphics[width=6cm]{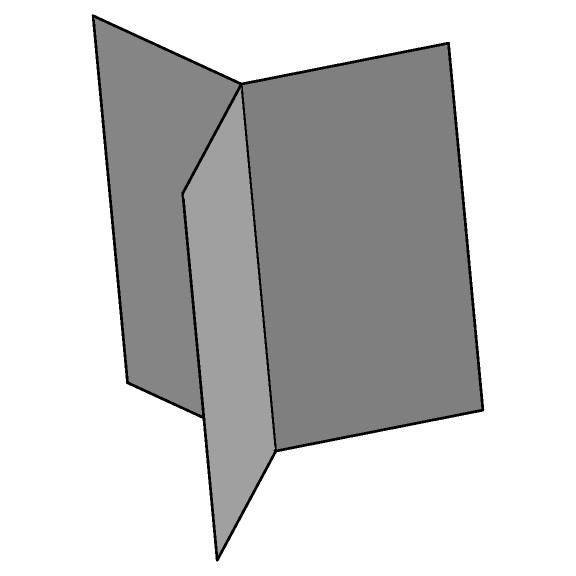} \hspace{2cm}
\includegraphics[width=6cm]{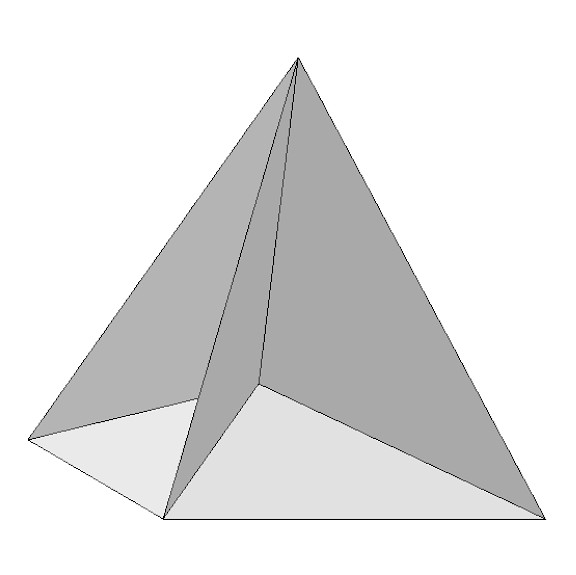}
\nopagebreak[4]

\nopagebreak[4] Cones\footnote{Thanks to Ken Brakke for those
pictures.} of type $\Y$ and $\T$.
\end{center}

Cones of type $\Pp$, $\Y$ and $\T$ are the only sets (except the
empty set) in $\R^3$ that locally minimizes the Hausdorff measure
of dimension 2 under topological conditions (i.e. every competitor
keep the same connected components outside the competitor ball).
This fact is proved in \cite{d5}. That is why in the following we
will say ``minimal cones'' to design cones of type $\Pp$, $\Y$ and
$\T$.

So here is now the definition of $\varepsilon_0$-minimal sets.

\begin{defin}\label{epsmin} Let $B$ be a ball in  $\R^3$. A closed set $E\subset B$
is said to be $\varepsilon_0$-minimal in $B$ if for all $x\in E$
and for all $r$ such that $B(x,r)\subset B$ we have that
\begin{eqnarray}
\inf_{Z \ni x}\left\{\frac{1}{r}\sup\left\{d(y,Z);y \in E\cap
B(x,r)\right\} \right\}\leq \varepsilon_0  \label{infh}
\end{eqnarray}
where the infimum is taken over all the minimal cones of type
$\Pp$,$\Y$, and $\T$ that contain $x$ (but are not necessarily
centered at $x$).
\end{defin}

Note that $\varepsilon_0$-minimal sets have nothing in common with
minimal sets or almost minimal sets.  The name ``minimal'' only
comes from  the minimal cones. Moreover, in this paper we don't
use the fact that cones $\Y$ and $\T$ are minimal. One could prove
a similar result with other cones that have good topological,
flatness and hierarchy properties.

We will also use this definition of ``separating condition".

\begin{defin}[Separating]\label{separex} Let $Z$ be a minimal cone in
$\R^3$ and  $B$ a ball of radius $r$ such that  $B\cap Z\not =
\emptyset$. For all $a>0$ we define $Z_{a}$ by
$$Z_{a}:= \{ y \in B ; d(y,Z)\leq a\}.$$
 Let $E$ be a closed set in $B$ such that $E$ is contained in
$Z_{r\varepsilon_0}$ for some $\varepsilon_0<10^{-5}$. We say that
``$E$ is separating in $B$'' if the connected components  of
$B\backslash Z_{r\varepsilon_0}$  are contained in different
connected components of $B\backslash E$. We denote by
$\mathpzc{k}^B$ the number of connected component of $B\backslash
Z_{2\varepsilon_0}$ (thus $\mathpzc{k}^B$ is equal to $type(Z)+1$
if $Z$ is not centered too close to $\partial B$).
\end{defin}

  Definition
\ref{epsmin} is introduced in \cite{ddpt} but in a slightly
different way. In \cite{ddpt} the inequality \eqref{infh} is
replaced by
$$\inf_{Z \ni x}D_{x,r}(E,Z)\leq \varepsilon_0$$
 where $D_{x,r}$ is the normalized Hausdorff distance defined in
 \eqref{hausdorffdist}. With this modification, $\varepsilon_0$-minimal sets will be
called ``strong $\varepsilon_0$-minimal sets''. In our case
(``weak $\varepsilon_0$-minimal sets'') we consider the first half
of the Hausdorff distance so we allow $E$ to contain some holes.
However, we will always suppose that our set is also separating in
$B$. In general, and for technical reasons, we will always
consider separately the topological separating condition and the
closeness to minimal cones.

The main result on  strong $\varepsilon_0$-minimal sets is Theorem
2.1 of \cite{ddpt} which says that for  $\varepsilon_0$ small
enough, every strong $\varepsilon_0$-minimal set is locally the
bi-hölderian image of a minimal cone. This is a generalization of
the Reifenberg's topological disc theorem that we mention below.
So $\varepsilon_0$-minimal sets can be seen as a generalization of
Reifenberg-flat sets.

In particular, a strong $\varepsilon_0$-minimal set is a
separating set. This is a consequence of \cite{ddpt} but one might
prove it without using the whole result in \cite{ddpt}. A
consequence is that a strong $\varepsilon_0$-minimal set in $B$ is
a weak $\varepsilon_0$-minimal set that separates in $B$. That is
why our result, that will be stated later for weak
$\varepsilon_0$-minimal with the separating condition, applies
directly to strong $\varepsilon_0$-minimal set in $B$. As a
result, it applies for $\varepsilon_0$-Reifenberg-flat sets hence
we deduce Theorem \ref{theorem1}. And we also have

\begin{theo}\label{theorem2} For all  $\gamma <1$ and
$0<r<\frac{1}{2}$, there is an  $\varepsilon_1 >0$ such that for
every strong $10^{-6}$-minimal set $K$ in the unit ball of $\R^N$
with $Z^0$ a minimal cone satisfying
$$\sup\{d(x,Z^0); x \in K \cap B(0,1)\}\leq \varepsilon, $$
 we have that
\begin{eqnarray}
\omega_2(0,r)\leq r^{\gamma}\omega_2(0,1) \label{inl}
\end{eqnarray}
for all locally energy minimizing function in $B(0,1)\backslash
K$.
\end{theo}

We will use the notation
\begin{eqnarray}
\beta(x,r):=\inf_{Z \ni x}\left\{\frac{1}{r}\sup\left\{d(y,Z);y
\in E\cap B(x,r)\right\} \right\}  \label{defbeta}
\end{eqnarray}

For readers who are familiar with $\beta$-numbers, note that here
we don't take a bilateral definition. We are now ready to define
$(\varepsilon_0,\varepsilon)$-minimal sets, and state the main
theorem that will in particular imply Theorem \ref{theorem1} and
Theorem \ref{theorem2}.

\begin{defin} \label{presqueepsplat}Let $E$ be a closed set with locally finite $H^2$ measure in $\R^3$.
Let $\varepsilon$ and $\varepsilon_0$ be two positive  constants
such that $0<\varepsilon<\varepsilon_0<10^{-5}$. We say that $E$
is $(\varepsilon_0,\varepsilon)$-minimal if there is a constant
$C_0$ and a family of balls
  $\{B_i\}_{i \in I}$:=$\{B(x_i,r_i)\}_{i \in I}$ such that $\{2B_i\}$ is of bounded cover with constant $C_0$,
   centered on  $E$, and such that :\\
$i)$ $\forall i \in I;$ $ r_i \leq \varepsilon $.\\
$ii)$  $E \backslash \bigcup_{i \in I}B(x_i,r_i)$  is $\varepsilon_0$-minimal in $B$. \\
$iii)$ There is a minimal cone $Z$ centered at the origin such
that
$$E\subset \{y \in B(0,1);d(y,Z)\leq \varepsilon\}.$$
$iv)$ For all $i\in I$ and for all $r>r_i$ with $B(x_i,r)\subset B$, we have $\beta(x_i,r)\leq \varepsilon_0$.\\
$v)$ $E$ is separating in $B$.
\end{defin}

The separating condition in $v)$ uses the cone of $iii)$.

So basically, a $(\varepsilon_0,\varepsilon)$-minimal set is a
$\varepsilon_0$-minimal set except in a collection of tiny bad
balls $B_i$ of radius less than $\varepsilon$. In general,
$(\varepsilon_0,\varepsilon)$-minimal sets will be obtained by a
stopping time argument. If we take a closed set $E$ and we do a
stopping time argument on $E$ with the stopping condition  of not
being close to a cone, and if we manage to control the radii of
all the stopping balls by $\varepsilon$, then we would have a
$(\varepsilon_0,\varepsilon)$-minimal set. This is what we will do
in a second paper to prove regularity for the singular set of the
Mumford-Shah functional.

So we want to prove a decay of normalized energy in the complement
of $(\varepsilon_0,\varepsilon)$-minimal sets. In fact this type
of sets is too general to hope to obtain such a result since
$(\varepsilon_0,\varepsilon)$-minimal sets allows the existence of
little tubes hidden in the bad balls $B_i$ and that could carry
some energy from exterior to interior. This is why we consider a
``cutting sphere'' of a certain radius $\rho$ that is not fixed
but just belongs to $[\frac{1}{2},\frac{3}{4}]$.

For all $i \in I$ we call
\begin{eqnarray} \label{defbpr}
B'_i:=C_1B_i=B(x_i,C_1r_i)
\end{eqnarray}
 where $C_1>1$ is a constant that will be chosen later, depending on $C_0$ and other
 geometric constants.

For all $\rho \in [\frac{1}{2},\frac{3}{4}]$ and for all
$(\varepsilon_0,\varepsilon)$-minimal set $E$ we define
$$I_\rho:=\{i \in I ; B_i \cap \partial B(0,\rho) \not = \emptyset \}.$$
 We define also
$$E^\rho := (E\backslash  \bigcup_{i \in I_\rho} B'_i) \cup \bigcup_{i \in I_\rho}\partial B'_i$$
and

$$U(E^\rho):=\Big\{u \in W^{1,2}(B\backslash E^\rho) ; u=
\text{argmin} \big\{\int_{B\backslash E^\rho}|\nabla v|^2;
v|_{\partial B\backslash E^\rho}=u \big\},
 \int_{B\backslash E^\rho}|\nabla u|^2=1\Big\}.$$
Since $\varepsilon$ is smaller than $(1-\rho)/C_1$, all the
functions in $U(E^\rho)$ are constant in each $B'_i$ for $i \in
I_\rho$.

Let  $\varepsilon_0<10^{-5}$ be fixed. For all $\varepsilon
<\frac{1}{4} $ we may introduce
$$\Lambda(\varepsilon):=\{(u,E,\rho) \text{ that verify  } (*) \}  $$
$$
\begin{array}{cl}
(*) & \left\{
\begin{array}{rl}
 E &\text{ is } (\varepsilon_0,\varepsilon )\text{-minimal in }  \bar B(0,1) \text{ and } 0 \in E\\
 \rho& \in [\frac{1}{2}, \frac{3}{4}] \\
 u & \in U(E^\rho) \\
\end{array}
\right.
\end{array}
$$
 For $i\in \{1,2,3\}$ we denote by
$\Lambda_i(\varepsilon)$ the elements of $\Lambda(\varepsilon)$
such that the cone $Z$ of condition  $iii)$ of Definition
\ref{presqueepsplat} is of type $i$.

For all $(u,E,\rho)\in \Lambda(\varepsilon)$  and all $0<r<1$,
recall that
$$\omega_2(0,r):= \frac{1}{r^2}\int_{B(0,r)\backslash E^\rho}|\nabla u(x)|^2dx . $$

We now come to the main result.

\begin{theo} \label{decreifenberg} For all $i\in \{1,2,3\}$, $\varepsilon_0\leq 10^{-5}$, $\gamma <0,8$ and
$0<r<\frac{1}{2}$, there is an  $\varepsilon_2 >0$ such that for
all $(u,E,\rho)\in \Lambda_i(\varepsilon_2)$  we have
\begin{eqnarray}
\omega_2(0,r)\leq r^{\gamma}.
\end{eqnarray}
\end{theo}

The proof of Theorem \ref{decreifenberg} is by compactness. So the
first step is to prove a decay estimate on the energy when $K$ is
a minimal cone. This is done in Section \ref{section3cone}. In
fact we will prove that in this case the normalized energy
increases like a power of radius. To show this, we use an argument
that A. Bonnet used in dimension 2, that link the decay of
normalized energy and the spectrum of the spherical Laplacian on
$\partial B(0,r)\backslash K$. In particular we give a lower bound
for the first eigenvalue on those domains.

Section \ref{section1} and Section \ref{section2whit} are devoted
to a general  method to obtain extensions of functions $u$ near
$K$ by a Whitney type construction. This work gives some useful
competitors for energy minimizing functions in the complement of
$(\varepsilon_0,\varepsilon)$-minimal sets. This approach might
give a powerful tool that could by applied in other problems. In
particular it will be used twice in a next paper that proves a
regularity result for the singular set of the Mumford-Shah
functional \cite{l2}.

Finally in Section \ref{section4proff}, we give the proof of
Theorem \ref{decreifenberg}. The argument is by contradiction and
compactness. It is inspired by what L. Ambrosio, N. Fusco, and D.
Pallara did in their theorem about the decay of energy in
\cite{afp1} (see also Theorem 8.19 in \cite{afp}). The difficulty
here, and also the key ingredient of the proof, is to estimate the
energy close to the set $E$ (inequality \eqref{lemmefonda}). In
their case,  L. Ambrosio, N. Fusco, and D. Pallara obtained this
estimate by approaching the set $E$ with a lipschitz surface and
by controlling  the difference with the ``Tilt estimate". Here, we
shall use the Whitney extension from the preceding sections.


\section{Monotoniciy in the complement of a minimal cone.}
\label{section3cone}

We want to prove first  that if $K=Z$ with $Z$ a cone of type $\Y$
or $\T$ centered at the origin and $u$ locally minimizes the
energy in $B(0,1)\backslash Z$ then
\begin{eqnarray}
\omega_2(0,r)\leq r^{0,8}\omega_2(0,1)\quad \forall r<1.
\label{dect1}
\end{eqnarray}
 To prove \eqref{dect1}, we will adapt an argument that A.
Bonnet \cite{b} used in dimension $2$, and to do this, we will
need the following lemma.

\begin{lem} \label{eigen} Let $Z$ be a minimal cone in $\R^3$ centered at $0$
and let $\Omega_r$ be a connected component of $\partial
B(0,r)\backslash Z$.  Then for all function $f\in
W^{1,2}(\Omega_r)$ we have
\begin{eqnarray}
\int_{\Omega_r}|f-m_f|^2dw\leq \frac{1}{2}r^2\int_{\Omega_r}
|\nabla_\tau f|^2dw. \label{coco}
\end{eqnarray}
\end{lem}

{\bf Proof :} Let $\lambda_1$ be the first positive eigenvalue for
$-\Delta_S$ (spherical Laplacian) in $\Omega_r$ with Neumann
condition on the boundary of $\Omega_r$. Then we have
$$\int_{\Omega_r}|f-m_f|^2dw\leq \frac{1}{\lambda_1}r^2\int_{\Omega_r} |\nabla_\tau f|^2dw. $$
If $Z$ if of type $\Y$, then by Lemma 4.1. of \cite{da} applied
with $\omega=\frac{2\pi}{3}$ and with Neumann boundary conditions,
we get $\lambda_1=2$ and \eqref{coco} follows.

So we have to consider the case when $Z$ is of type $\T$. Let $f$
be a eigenvector for the first positive Neumann eigenvalue. For
$1\leq i\leq 3$ we denote by $\delta_i$ the three symmetry axis of
$\Omega_1$ and we denote by $s_i$ the corresponding symmetries. If
$f$ is symmetric by all the three axis, then by reflection we
could extend $f$ to the entire sphere for which it is well known
that the first eigenvalue is equal to $2$.

If $f$ is not symmetric by the three axis, then there is one, for
instance $\delta_1$, such that $f$ is not symmetric. Then we
consider the anti-symmetric function
$$g=f-f\circ s_1.$$
We can suppose that $g$ is positive in $\Omega_1$, and $g$ is
still an eigenvector associated to $\lambda_1$. Moreover, we have
that $g$ vanishes on $\delta_1$. The axis $\delta_1$ cut
$\Omega_1$ in two isometric triangles. We call  $\Omega'$ one of
them. Now we consider $\Sigma$ a connected component of  $S^2
\backslash Y$ where $Y$ is of type $\Y$ and we suppose that
$\Sigma$ contains
 $\Omega$. $\Sigma$ is also cut by  $\delta_1$ and we denote by $\Sigma'$ the connected component
 that contains
  $\Omega'$. Now we apply Proposition 4.3. of \cite{da} with $G=\Omega'$, $\partial_D
G=\delta_1\cap \partial \Omega'$, $\partial_N G=\partial \Omega'
\backslash \delta_1$, $G'=\Sigma'$, $\partial_NG'=\partial
\Sigma'\backslash \delta_1$ and $\partial_D G'=\delta_1 \cap
\partial \Sigma'$. We obtain that
$$\lambda_1\geq \mu(\Sigma')$$
where $\mu(\Sigma')$ is the first positive eigenvalue in $\Sigma'$
with Neumann condition on  $\partial \Sigma'\backslash \delta_1$
and Dirichlet condition on $\delta_1 \cap
\partial \Sigma'.$
Now applying Lemma 4.1. of \cite{da} again with
$\omega=\frac{2\pi}{3}$, but now with mixed boundary conditions,
we get $\mu(\Omega')=2$, thus
$$\lambda_1\geq 2$$
and the lemma follows.\qed

\begin{lem}\label{bonnet} Let $Z$ be a minimal cone in $\R^3$ centered at
$0$. Then for all local energy minimizer  $u$ in $B(0,1)\backslash
Z$ and for all $a,r<1$ we have
$$\omega_2(0, ar)\leq a^{0,8}\omega_2(0,r).$$
Moreover, $r\mapsto \omega_2(0,r)$ is increasing.
\end{lem}

{\bf Proof :} Set
$$E(r) = \int_{B(0,r)\backslash Z}|\nabla u|^2.$$
For almost every $r$ the derivative of  $r\mapsto E(r)$ exists,
$$E'(r)= \int_{\partial B(0,r) \backslash Z}|\nabla u|^2\quad \text{
and } \quad E(r)= \int_{0}^r E'(t)dt$$ (see Lemma 47.4 page 316 of
\cite{d}).

We want to prove an inequality of type
$$E(r)\leq CrE'(r)$$
where $C$ is a contant that will be explicited later.

Firstly, since $u$ is a harmonic function and since
$\frac{\partial }{\partial n}u =0$ on $Z$ (where $Z$ is regular),
an integration by parts gives
\begin{eqnarray}
E(r)=\sum_{j=1}^{J}\int_{S_j} u \frac{\partial u}{\partial n}
\label{ipp1}
\end{eqnarray}
with $S_j$ the connected components of $\partial B(0,r) \backslash
Z$. To justify the integration by parts in this "regular
polyhedral domain", one could find a lemma in \cite{l1}. In
addition, one might prove \eqref{ipp1} without integrating by
parts but just by using that $u$ is energy minimizing as in
\cite{d} page 320.

Denote by $A_j$ the connected components of $B(0,r)\backslash Z$
which boundary contains $S_j$. An other integration by parts in
$S_j$ gives
$$\int_{S_j}\frac{\partial u}{\partial
n}=\int_{A_j} \Delta u =0.$$ Thus we can subtract by a constant
and we find
\begin{eqnarray}
\int_{S_j} u \frac{\partial u}{\partial
n}&=&\int_{S_j}[u-c_j(u)]\frac{\partial u}{\partial r} \notag \\
&\leq&\left[\int_{S_j}[u-c_j(u)]^2\right]^{\frac{1}{2}}\left[\int_{S_j}\left(\frac{\partial
u}{\partial r}\right)^2\right]^{\frac{1}{2}}. \notag
\end{eqnarray}
Then by use of $ab\leq \frac{1}{2}[\lambda^{-1}a^2+\lambda b^2] $
with $\lambda $ a positive constant to be chosen later,
\begin{eqnarray}
\int_{S_j} u \frac{\partial u}{\partial n}&\leq
&\frac{1}{2\lambda}\left[\int_{S_j}[u-c_j(u)]^2\right]+\frac{\lambda
}{2}\int_{S_j}\left(\frac{\partial u}{\partial r}\right)^2 \notag
\\
&\leq&\frac{1}{4\lambda}r^2\int_{S_j}|\nabla_{\tau
}u|^2+\frac{\lambda }{2}\int_{S_j}\left(\frac{\partial u}{\partial
r}\right)^2 \notag
\end{eqnarray}
then by setting $\lambda =\frac{r}{\sqrt{2}}$,
$$\int_{S_j}u\frac{\partial u}{\partial n}\leq \frac{1}{2\sqrt{2}} r\int_{S_j}|\nabla u|^2.$$
Finally, summing over $j$,
$$E(r)\leq \frac{1}{2\sqrt{2}} r E'(r).$$
This estimate shows that the derivative of $\omega_2(0,r)$ is
positive, in other words $r\mapsto\omega_2(0,r)$ is increasing. To
have now the estimate about the speed of increasing, consider
$g(r):=\ln(E(r))$. By absolute continuity we have
$$g(r)-g(ar)=\int_{ar}^{r}\frac{E'(t)}{E(t)}dt \geq \int_{ar}^{r}\frac{2\sqrt{2}}{t}dt \geq 2\sqrt{2}\ln(\frac{1}{a}).$$
Hence
$$E(ar)\leq a^{2\sqrt{2}}E(r).$$
Then we divide by $(ar)^2$
$$\frac{1}{(ar)^2}E(ar)\leq a^{2(\sqrt{2}-1)}\frac{1}{r^2}E(r)$$
and this implies the lemma because $2(\sqrt{2}-1)\geq 0,8$.\qed


\section{Some geometric lemmas}
\label{section1}

Before doing the Whitney extension, we have to discuss about some
geometric facts. We begin this section with a lemma that will
allow us to work with ``almost centered cones''.

\subsection{The Recentering Lemma}

\begin{lem}[Recentering]\label{recentrage} Let $Z$
be a minimal cone in $\R^3$ that contains $0$ (but is not
necessarily centered at $0$). Then for all $r_0>0$ and for all
constant $V\geq 1$ there is a $r_1$ such that
$$r_1\in \{r_0,Vr_0, V^2 r_0\}$$
and such that we can find a cone $Z'$, containing $0$ and centered
in $B(0,\frac{1}{V}r_1)$ with $Z\cap B(0,r_1)=Z'\cap B(0,r_1)$.
\end{lem}

{\bf Proof :} Let consider the ball $B(0,r_0)$. If  $type(Z)=1$
then $Z$ is a plane that contains $x$, thus its center is $x$ and
we can take  $r_1=r_0$.

Suppose now that  $Z$ is of type $\Y$ and assume that all the
points on the spine of $Z$ are in $\R^3\backslash
B(0,\frac{1}{V}r_0)$ (otherwise we could take $r_1=r_0$ that is
case 1 of Figure 1). If there is no point on the spine of $Z$ in
$B(0,r_0)$, we can take $Z'$ a plane that is equal to $Z$ in
$B(0,r_0)$. Otherwise we are in case 2 of Figure 2 and we can take
$r_1=Vr_0$.

\begin{center}
Figure 1

 \nopagebreak[4]
\input{dYcas1.txt}\hspace{3cm}
\input{dYcas2.txt}

 \nopagebreak[4]
\vspace{1cm} Case 1 \hspace{4cm} Case 2
\end{center}

Now it remains to consider the case  when $Z$ is of type $\T$. We
discuss it in the same way. If the center of $Z$ is not in
$B(0,r_0)$, and if there is no point on the spine of a cone of
type $\Y$ in $B(0,r_0)$, then we take $r_1=r_0$ and for $Z'$ we
take a plane.

Now if there is some point of type $\Y$ in $B(0,r_0)$ but no point
of type $\T$, and if in addition there is a $\Y$ spine passing
trough $B(0,\frac{1}{V}r_0)$ we can take for $Z'$ a cone of type
$\Y$ and fix $r_1=r_0$. Now if there is some point of type $\Y$ in
$B(0,r_0)\backslash B(0,\frac{1}{V}r_0)$, we try $r_1=Vr_0$. If
there is no point of type $\T$ in $B(0,Vr_0)$, we can take for
$Z'$ a cone of type $\Y$ and $r_1=Vr_0$ (case 1 of Figure 2).
Otherwise we take for $Z'$ a cone of type $\T$ and we fix
$r_1=V^2r_0$.

\begin{center}
Figure 2

 \nopagebreak[4]
\input{dTcas1e.txt} \hspace{1cm}
\input{dTcas2.txt}

\nopagebreak[4]  Case 1 \hspace{6cm} Case 2
\end{center}

Finally if  $Z$ is of type $\T$ and its center lies in $B(0,r_0)$,
then if the center lies in $B(0,\frac{1}{V}r_0)$ we take $r_0=r_1$
and $Z=Z'$, otherwise we take $r_1=Vr_0$ and $Z'=Z$.\qed

\begin{defin}[Almost Centered] Let $Z$ be a minimal cone and  $B$ a ball that meets  $Z$. We say that $Z$ is
almost centered with constant $V$ if the center of $Z$ lies in
$\frac{1}{V}B$. If $V=2$ we just say that $Z$ is almost centered
in $B$.
\end{defin}

\subsection{The geometric function and general assumptions}

Here we describe the general situation that will appear in next
sections. Let $K$ be a closed set in $\bar B(x_0,r_0)$ such that
$H^{2}(K\cap \bar B(x_0,r_0))< + \infty$. Suppose that there is a
positive constant $\varepsilon_0 <10^{-5}$ and a minimal cone $Z$,
centered at $x_0$, such that
\begin{eqnarray}
\sup \{ d(x,Z) ; x \in K\cap B(x_0,r_0)\} \leq r_0\varepsilon_0
\label{sect1}
\end{eqnarray}
 and that $K$ is
separating in $B(x_0,r_0)$. For all $x \in K\cap B(x_0,r_0)$ and
$r>0$ such that $B(x,r)\subset B(0,r_0)$ recall that
$$\beta(x,r)=\inf_{Z \ni x}\frac{1}{r}\sup \{d(x,Z); x \in K\cap B(x,r)\}$$
Let $\rho \in [\frac{1}{2}r_0,\frac{3}{4}r_0]$ and assume that we
have an application
\begin{eqnarray}
\delta: B(x_0,\rho) \to [0,\frac{1}{4}r_0] \label{sect12}
\end{eqnarray}
 with the property that
\begin{eqnarray}
\beta(x,r)\leq \varepsilon_0, \text{ for all } x \in K\cap
B(x_0,\rho) \text{ and } r \text{ such that  } \delta(x)\leq r
\leq \frac{1}{4}r_0 \label{sect13}.
\end{eqnarray}
In addition we suppose that
\begin{eqnarray}
\delta \text{ is } C_0-\text{Lipschitz}\label{sect14}.
\end{eqnarray}
 The application $\delta$ will be called the ``geometric
function''.

\begin{defin}[Hypothesis $\mathcal{H}$] \label{defhyph} We will say that a closed set $K\subset B(x_0,r_0)$ with finite
$H^2$ measure is satisfying hypothesis
$\mathcal{H}$ if \\
i) There is a minimal cone $Z$ that verify \eqref{sect1} for a
``geometric constant"
 $\varepsilon_0<10^{-5}$ and a ``Lipschitz constant" $C_0$. \\
ii)  $K$ is separating in $B(x_0,r_0)$.\\
iii) There is a geometric function $\delta$ satisfying
\eqref{sect12}, \eqref{sect13} and \eqref{sect14} for a radius
$\rho \in [\frac{1}{2}r_0, \frac{3}{4}r_0]$.
\end{defin}

For example, if we have $\beta(x_0,r_0)\leq \varepsilon$ with
$\varepsilon<\frac{1}{4}\varepsilon_0$, and if in addition $K$ is
separating, then we have Hypothesis $\mathcal{H}$ with
$\delta(x)=\frac{\varepsilon}{\varepsilon_0}r_0$ everywhere.

An other example is given by a Reifenberg-flat set $K$ included in
$B(x_0,r_0)$, containing $x_0$ and with constant $\varepsilon_0$
less than $10^{-6}$. Then we have Hypothesis $\mathcal{H}$ on $K$
with $\delta=0$ everywhere.

\begin{sloppypar}
 Under Hypothesis $\mathcal{H}$ we will always denote by
$A_k(x_0,r_0)$ for $k \in \N \cap [1,\mathpzc{k}^{B(x_0,r_0)}]$
the connected components of
 $$B(x_0,r_0)\backslash \{ y \in B(x,r); d(y,Z)\leq \varepsilon_0r_0\}$$
 and we will call $\Omega_k(x_0,r_0)$ the connected component of
 $B(x_0,r_0)\backslash K$ that contains $A_k(x_0,r_0)$.
\end{sloppypar}

\subsection{The orientation lemma}

Now we have to discuss orientation and separation.

\begin{lem}[Orientation]\label{bienordef} Let $K$ be a closed set in $B(x_0,r_0)$ satisfying Hypothesis
$\mathcal{H}$ with a geometric function $\delta$, a minimal cone
$Z$ and a constant $\varepsilon_0<10^{-5}$.
 Let $B(x,r)$ be a ball
included in $B(x_0,r_0)$ such that $x\in K$ and
$$\delta(x)\leq r \leq \frac{r_0}{32}.$$
Let $r_1$ be the radius (equal to $r$, $2r$ or $4r$) such that
$Z(x,r_1)$ is almost centered in $B(x,r_1)$ where $Z(x,r_1)$
denotes the minimal cone $\varepsilon_0$-close to $K$ in
$B(x,r_1)$. We also call
   $A_k(x,r_1)$, $k \in \N \cap
[1,\mathpzc{k}^{B(x,r_1)}]$, the connected components of
 $$\{ y \in B(x,r_1); d(y,Z(x,r_1))\geq \varepsilon_0r_1\}.$$
Then $B(x,r_1)$ is well oriented in $B(x_0,r_0)$. That means that
$B(x,r_1)$ verifies the two following points:\\
i) $K$ is separating in $B(x,r_1)$.\\
ii) There is an injective application
$$l:\N \cap
[1,\mathpzc{k}^{B(x,r_1)}]\to\N \cap
[1,\mathpzc{k}^{B(x_0,r_0)}]$$
 such that $A_k(x,r_1)\subset \Omega_{l(k)}(x_0,r_0)$.
\end{lem}

{\bf Proof :} We consider the balls $B^p:=B(x,2^pr_1)$ for $p\in
\N \cap [0,P]$ where $P$ is such that
$$\frac{1}{16}r_0\leq 2^Pr_1 \leq \frac{1}{8}r_0.$$
It is always possible because $r_1 \leq \frac{1}{8}r_0$.
Therefore, every $B^p$ is included in $B(x_0,r_0)$. We define also
$B^{P+1}:=B(x_0,2r_0)$ and $Z^p$ is the minimal cone that is
$\varepsilon_0$-close to $K$ in $B^p$  (we know that there is one
for all $p$ because the radius of $B^p$ is larger than
$\delta(x)$). Now we use Lemma \ref{recentrage} to extract among
the $B^{p}$ a subsequence $B^{\sigma(p)}$ such that for all $p$,
$B^{\sigma(p)}$ is almost centered and such that the radius of
each ball is not larger than eight times the radius of the
preceding one.
 We still denote by $B^p$ this subsequence (instead of $B^{\sigma(p)}$).
 $B^0=B(x,r_1)$  is the beginning ball   and   $B^{P+1}$ is $B(x_0,r_0)$.
Nevertheless, the radius of $B^p$ is not exactly $2^pr_0$ but
equivalent with a factor $4$.

Now we are going to prove by induction that $B^p$ is well oriented
for all $p$.

Hypothesis $\mathcal{H}$ in $B(x_0,r_0)$ clearly shows that
$B^{P+1}:=B(x_0,r_0)$ is well oriented. Next we consider  $p \in
\N \cap [1,P+1]$. We have to show that if $B^{p}$ is well
oriented, then  $B^{p-1}$ is well oriented. We denote by
$\mathpzc{k}^p$ instead of $\mathpzc{k}^{B^p}$ the number of
connected components and we denote by $r_p$ the radius of $B^p$.
We know that
\begin{eqnarray}
K \cap B^p \subset {Z^p_{r_{p}\varepsilon_0}}:=\{y  ; d(y,Z^p)\leq
r_p\varepsilon_0\} \label{rp}
\end{eqnarray}
and in addition $K$ is separating in $B^p$. Let  $A^p_k$ be for $k
\in \N \cap [1,\mathpzc{k}^p]$  the connected components of
$B^p\backslash Z^p_{r_p\varepsilon_0}$.

Now consider $B^{p-1}$. Firstly note that  $K$ is separating in
$B^{p-1}$. To see this we have to show that connected components
of $B^{p-1}\backslash Z^{p-1}_{r_{p-1}\varepsilon_0}$ are in
different connected components of $B^{p-1}\backslash K$. Let
$A^{p-1}_k$ be the connected components of $B^{p-1}\backslash
Z^{p-1}_{r_{p-1}\varepsilon_0}$ for $k \in \N\cap
[1,\mathpzc{k}^{p-1}]$. Since
$$K \cap B^{p-1} \subset Z^{p-1}_{r_{p-1}\varepsilon_0}$$
and $\varepsilon_0 <  10^{-5}$, we can choose for all $k$ a point
$a^{p-1}_k$ such that $a^{p-1}_k \in A^{p-1}_k$, and $d(a_k^{p-1},
K)\geq \frac{r_{p-1}}{10} $ (because the cones $Z^p$ are almost
centered).

 Since  $d(a_k^{p-1}, K)\geq \frac{r_{p-1}}{10} $, then
$d(a_k^{p-1}, Z^{p-1}_{r_{p-1}\varepsilon_0})\geq
r_{p-1}(\frac{1}{10}-\varepsilon_0)$ and by use of \eqref{rp} we
can deduce that $d(a_k^{p-1},Z^{p}_{r_{p} \varepsilon_0})\geq r_p
(\frac{1}{10}-10\varepsilon_0)$. It follows that for $k \in \N\cap
[1,\mathpzc{k}^{p-1}]$ there is an $l \in \N
\cap[1,\mathpzc{k}^{p}] $ such that $a_k^{p-1}\in A^{p}_l$.

In addition,
\begin{eqnarray}
 \text{the } a_{k}^{p-1} \text{ are in different connected
components of } B^{p}\backslash Z^p_{\varepsilon_0r_p}.
\label{finiouquoi}
\end{eqnarray} Indeed, suppose that there is $k_1$ and $k_2$ such
that $a_{k_1}^{p-1}$ and $a^{p-1}_{k_2}$ are both in the same
connected component of $B^{p}\backslash Z^p_{\varepsilon_0r_p}$.
Therefore, since $d(a_k^{p-1},Z^{p}_{r_p\varepsilon_0})\geq
r_p\frac{1}{11} $ and since the $Z^p$ are almost centered, we can
deduce that there is a continuous path  $\Gamma$ from
$a_{k_1}^{p-1}$ to $a^{p-1}_{k_2}$ such that all the points of
$\Gamma$ are situated at a distance larger than $\frac{1}{11}r_p$
from $Z^p$, consequently larger than $\frac{1}{100}r_{p-1}$ from
$Z^{p-1}$ which is a contradiction with the definition of
$a_k^{p-1}$. Consequently
$$\mathpzc{k}^{p-1}\leq \mathpzc{k}^p \leq \mathpzc{k}^{B(x_0,r_0)}.$$
Now if $K$ is not separating in $B^{p-1}$, then there is $k_1$ and
$k_2$, and there is a continuous path from $a^{p-1}_{k_1}$ to
$a^{p-1}_{k_2}$ without meeting $K$. However we know that $K$ is
separating in $B^{p}$ hence we get the contradiction.

Since all the $A_k^{p}$ are included in a certain $\Omega_l$ and
since $K$ is separating in $B^{p-1}$, we can deduce that every
$A_k^{p-1}$ is also in an $\Omega_{l(k)}$. Moreover $l(k)$ is
injective by \eqref{finiouquoi} and by induction, so the
conclusion follows.\qed


\section{Whitney extension from a geometric function}
\label{section2whit}

We still assume that $K$ is a closed set in $B(x_0,r_0)$
satisfying Hypothesis $\mathcal{H}$ with a geometric function
$\delta$, a minimal cone $Z$, a constant $\varepsilon_0<10^{-5}$
and a radius $\rho\in [\frac{1}{2}r_0, \frac{3}{4}r_0 ]$. Let
$U>1$ be a constant that will be fixed later, depending on $C_0$
and a dimensional constant. In addition we assume that
$\varepsilon_0$ is very small compared to $U^{-1}$. For all $t>0$
we define
\begin{eqnarray}
\mathpzc{V}(t):= \bigcup_{x \in K\cap B(0,\rho)}B(x,
\frac{t}{U}\delta(x)). \label{defmore}
\end{eqnarray}
 We also set
\begin{eqnarray}
\mathpzc{V}:=\mathpzc{V}(10) \label{defQ1}
 \end{eqnarray} and
\begin{eqnarray}
\mathpzc{Z}:=\mathpzc{V}(30)\backslash \mathpzc{V}(\frac{1}{10}).
\label{defQ2}
\end{eqnarray}
Finally we define
\begin{eqnarray}
\mathpzc{V}_{\rho}:= \bigcup_{x ; B(x, \frac{10}{U}\delta(x))\cap
\partial B(x_0,\rho)\not = \emptyset}B(x, \frac{10}{U}\delta(x)). \label{defQ3}
\end{eqnarray}

Recall that  by hypothesis, $K$ is separating in $B(x_0,r_0)$ and
that for all $k \in [1,\mathpzc{k}^{B(x_0,r_0)}]$ we have denoted
by $A_k(x_0,r_0)$ the connected components of
$B(x_0,r_0)\backslash Z_{\varepsilon_0r_0}$ and by
$\Omega_k(x_0,r_0)$ the connected component of
$B(x_0,r_0)\backslash K$ that contains $A_k(x_0,r_0)$. We also set
\begin{eqnarray}
\Delta_k:=B(x_0,\rho)\cap (\Omega_k(x_0,r_0)\cup \mathpzc{V}).
\label{defQ4}
\end{eqnarray}

The purpose of this section is to prove the following lemma.

\begin{lem}\label{extentionW}{\rm(Whitney Extension)} Let $K$ be a closed
set in $B(x_0,r_0)$ satisfying
Hypothesis $\mathcal{H}$ with a geometric function $\delta$, a
minimal cone $Z$, a constant $\varepsilon_0<10^{-5}$ and a radius
$\rho\in [\frac{1}{2}r_0, \frac{3}{4}r_0]$. Then for all function
$u \in W^{1,2}(B(0,r_0)\backslash K)$, and for all $k \in
[1,\mathpzc{k}^{B(x_0,r_0)}]$, there is a function
$$v_k \in
W^{1,2}(\Delta_k \backslash \mathpzc{V}_{\rho})$$ such that
$$v_k = u \text{ in }B(x_0,\rho)\backslash \mathpzc{V}$$
and
\begin{eqnarray}
\int_{\Delta_k\backslash \mathpzc{V}_{\rho}}|\nabla v_k|^2dx \leq
\int_{\Delta_k \backslash \mathpzc{V(\frac{1}{3})}}|\nabla u|^2dx+
C\int_{\mathpzc{Z}}|\nabla u|^2dx \label{inwh1}
\end{eqnarray}
where $C$ is a constant depending only on dimension and where
$\mathpzc{V}$, $\mathpzc{V}_{\rho}$, $\mathpzc{Z}$ and $\Delta_k$
are defined in \eqref{defQ1}, \eqref{defQ3}, \eqref{defQ2} and
$\eqref{defQ4}$ with constant $U>30C_0$ depending also on
dimension.
 \end{lem}

{\bf Proof :} We will use a Whitney type extension. For all $x \in
K \cap B(x_0,\rho)$ set
$$B_x:= B(x,\frac{1}{U}\delta(x))$$ where $U \geq
30C_0$ is a constant that will be fixed later. Then we choose a
subfamily $\{W_j\}_{j \in J}$ of balls from $\{B_x\}$, maximal for
the property that
$$\forall j \not = j',\quad\frac{1}{100}W_j \cap
\frac{1}{100}W_{j'} = \emptyset.$$ We denote by $r_j$ the radius
of the ball $W_j$. The $\{W_j\}_{j\in J}$ is our Whitney family of
balls, and we denote by $x_j$ and $r_j$ the center and radius of
$W_j$. We have the following proprieties about the Whitney balls.

\begin{lem}\label{lemwhit}
\begin{eqnarray}
\text{Whitney Property} &i)&  10W_j \cap 10W_{j'} \not = \emptyset \Rightarrow \frac{1}{20}r_{j'} \leq r_j \leq 20 r_{j'}  \notag \\
\text{Covering }  \mathpzc{V}&ii)&   \text{ for all } t>0,\;
\bigcup_{j \in J} tW_j \subset  \mathpzc{V}(t) \subset \bigcup_{j
\in J}
 (20t+\frac{3}{10})W_j \notag \\
\text{The cover is bounded} &iii)& \exists C_2; \forall x \in
B(x_0,\rho),  \sharp\{j \in J; x \in 10W_j\}\leq C_2 \notag \\
\text{Geometry is under control} &iv)& \forall j \in J, \forall r
\in [r_j, \frac{1}{4}r_0], \beta(x_j,r)\leq U\varepsilon_0\notag
\end{eqnarray}
\end{lem}

{\bf Proof :} We begin with $i)$. Let $x \in B(x_0,\rho)$ and let
$j$ and $k$ be two indices such that $10W_{j}$ and $10W_k$ contain
$x$. Then if $x_j$ and $x_k$ are the centers of $W_j$ and $W_k$,
since $\delta$ is $C_0$-lipschitz we have
\begin{eqnarray}
\delta(x_j)&\leq& \delta(x_k)+C_0|x_k-x_j| \notag \\
&\leq& Ur_k + C_0|x_k-x|+C_0|x-x_j| \notag \\
&\leq & Ur_k+10C_0r_k+10C_0r_j .\notag
\end{eqnarray}
Therefore,
$$(U-10C_0)r_j\leq (U+10C_0)r_k.$$
With the same argument exchanging  $r_j$ and $r_k$ we can finally
deduce that
$$\frac{U-10C_0}{U+10C_0}r_j \leq r_k \leq \frac{U+10C_0}{U-10C_0}r_j$$
thus $i)$ follows if we consider that $U\geq 30C_0$.

The first inclusion of $ii)$ is trivial, by definition of
$\mathpzc{V}(t)$. So we have to prove the second inclusion. Let $y
\in \mathpzc{V}(t) $ and let $x$ be a point of $K$ such that $y
\in tB_x:= tB(x,\frac{1}{U}\delta(x))$. We denote by $r_x$ the
radius of $B_x$. Since $\{W_j\}$ is a maximal family, there is a
$j_0$ such that $\frac{1}{100}B_x\cap \frac{1}{100}W_{j_0}\not =
\emptyset$ (otherwise we could add $B_x$ in the family $\{W_j\}$).
Let $z$ be a common point. By the same argument as for $i)$, we
can deduce that $20^{-1}r_{j_0}\leq r_x \leq 20r_{j_0}$. Now if
$x_{j_0}$ is the center of $W_{j_0}$ we have
\begin{eqnarray}
|y-x_{j_0}|&\leq& |y-x|+|x-z|+|z-x_{j_0}| \notag \\
&\leq& tr_x +\frac{1}{100}r_x+\frac{1}{100}r_{j_0} \notag\\
&\leq & 20tr_{j_0}+\frac{20}{100}r_{j_0}+\frac{1}{100}r_{j_0} \notag \\
&\leq & (20t+\frac{3}{10}) r_{j_0}\notag
\end{eqnarray}
so $y \in (20t+\frac{3}{10}) W_{j_0}$ and this proves the second
inclusion of $ii)$.

For $iii)$ it is just a simple consequence of a geometric fact in
$\R^N$. Consider a family of balls in $\R^N$ that are all
containing a same point, with radius equivalent to $1$ and
centered at distance more than $\frac{1}{100}$ to each other, then
the number of these balls is finite. The proof of  $iii)$ follows.

Finally we have to prove $iv)$. Let $j \in J$ and  $r>0$ be such
that $r_j\leq r \leq r_0$. By definition of $\delta$, we know that
if $r\geq \delta(x_j)$, then  $\beta(x_j,r)\leq \varepsilon_0$ and
this is that we want. Now if $r_j\leq r \leq \delta(x_j)$ we have
$$\beta(x_j,r)\leq  \frac{\delta(x_j)}{r}\beta(x_j, \delta(x_j))\leq U\varepsilon_0$$
because $r\geq r_j=\frac{\delta(x_j)}{U}$, and $iv)$ is
proved.\qed

With help of Lemma \ref{recentrage}, for every $j \in J$ we can if
necessary change $W_j$ to $2W_j$ or $4W_j$ in order to have that
all the cones $Z^j$ associated to the  $W_j$ in $iv)$ are almost
centered in $W_j$. To prove this, using a translation we can
suppose that $W_j$ is centered at $0$. Then we apply Lemma
\ref{recentrage} to the cone $Z$ associated to $W_j$, which is
$U\varepsilon_0r_j$-close to $K$ in $W_j$. If
$U\varepsilon_0<10^{-5}$, Lemma \ref{recentrage} says that we can
choose among $W_j$, $2W_j$ or $4W_j$, a ball such that $Z$ is
almost centered. Since the cone $\tilde Z$ associated to this
choice of ball is close to $Z$, it is also almost centered itself.
This new family of balls still verify proprieties of Lemma
\ref{lemwhit} with constant that may be slightly different  (by
multiplying by $4$).

Moreover, applying Lemma \ref{bienordef} to $W_j$ we can assume
that $W_j$ is well oriented.

We are now ready to make our Whitney extension $v$ from $u$ in
$B(x_0,\rho)$. For every ball $W_j$ consider a function $\varphi_j
\in C^\infty$, with compact support in $10W_j$, equal to $1$ on
$8W_j$ and to $0$ out of $10W_j$.

\begin{lem}\label{defred} There is a function $\varphi_0 \in C^\infty$ such that
\begin{eqnarray}
\varphi_0 = 1 & \text{ in } & B(x_0,\rho) \backslash \bigcup_{j \in J} 10W_j \label{propvoulue1}\\
\varphi_0 = 0 & \text{ in }  & \bigcup_{j \in J} 8W_j \label{propvoule2} \\
\varphi_0+\sum_{j \in J}\varphi_j\geq 1 & \text{ in }& B(x_0,\rho)
\label{propvoulue3}
\end{eqnarray}
and in addition there is a constant $C$ such that for all $j \in J
$ and for all $x \in 10W_j\backslash 8 W_j$,
\begin{eqnarray}
|\nabla \varphi_0(x)|\leq C\frac{1}{r_j}. \label{propvoulue4}
\end{eqnarray}
\end{lem}

{\bf Proof  :} Let $l \in C^{\infty}(\R^+)$ such that $l$ is equal
to $0$ on $[0,8]$, equal to $1$ on $[10,+\infty[$ and
\begin{eqnarray}
l'(x)\leq 10. \label{estime}
\end{eqnarray}
Then the function
$$\varphi_0(x):=\prod_{j \in J}l\left(\frac{|x-x_j|}{r_i}\right)$$
satisfies the hypothesis of Lemma \ref{defred}. Indeed,
\eqref{propvoulue1} and \eqref{propvoule2} is obvious from the
definition of $\varphi_0$, \eqref{propvoulue3} is easy to prove
depending on a good construction of  the $\varphi_j$, and
\eqref{propvoulue4} follows from \eqref{estime}, Property $i)$ of
Lemma \ref{lemwhit} and also from the fact that the $\{10W_j\}$
are in bounded cover with equivalent radius when they meet each
other.\qed

Set
$$\theta_j=\frac{\varphi_j}{\varphi_0+\sum_{j\not = 0} \varphi_j} $$
so that we have a partition of unity on $B(x_0,r_0)$. Since the
$10W_j$ are in bounded cover, the sum is locally finite.

In each $10W_j$ there is a cone $Z^j$ such that
$$K \cap 10W_j \subset Z^j_{10r_jU\varepsilon_0}$$
with
$$Z^j_{10r_jU\varepsilon_0}:=\{y; d(y,Z^j)< \varepsilon_010Ur_j\}.$$
Denote by  $A_k^j$, $1\leq k \leq \mathpzc{k}^j$ the connected
components of $10W_j\backslash Z^j_{r_j\varepsilon_0}$. Since the
$10W_j$ are well oriented, we know that each $A_k^j$ is contained
in one of the big connected components $\Omega_{l(k)}$ of
$B(x_0,r_0)\backslash K$. And we know that $l(k)$ is injective.
Then we rename $A_k^j$ to $A_l^j$ such that $A_l^j$ is now
contained in $\Omega_l^j$. By convention, if $l$ is such that
$\Omega_l$ does not meet $W_j$, we set $A_l^j=\emptyset$.
 Hence for all $j$ and for $1\leq k \leq \mathpzc{k}^{B(x_0,r_0)}$, we have defined $A_k^j$.
 In each $A_k^j$ we choose a point $a_k^j \in 10W_j \backslash 8 W_j$ at greater  distance than $7r_j$
  from $K$ and we also consider $D_k^j$ a ball centered at $a_k^j$ and of radius  $\frac{1}{100}r_j$.
  We denote by  $m_k^j$ the mean value of $u$ on $D_k^j$. It is
  always possible because the cones $Z^j$ are almost centered.

Now for all $k \in \N \cap [1,\mathpzc{k}^{B(x_0,r_0)}]$ set
$$v_k(x):= \varphi_0(x) u(x)+\sum_{j>0} m^j_k \theta_j (x).$$
The $v_k$ are well defined on $B(x_0,r_0)\backslash K$ and since
$\bigcup_{j\in J}10 W_j \subset \mathpzc{V}$, the $v_k$ are equal
to $u$ in the exterior of $\mathpzc{V}$.

We want to check that the functions $v_k$ belong to
$W^{1,2}(\Delta_k\backslash \mathpzc{V}_{\rho})$ and we want to
estimate their energy with the energy of $u$. Let $x$ be a fixed
point in $\bigcup_{j \in J}10W_j$ and call $J_x$ the set of all
indices $j$ such that $10W_j$ contains $x$. We know by Lemma
\ref{lemwhit} $i)$ that all these balls have equivalent diameters.

On the other hand, all the $10W_j$ for $j \in J_x$ are included in
$30W_{j_0}$ where $j_0$ is any index fixed in $J_x$. By Property
$v)$ of Lemma \ref{lemwhit}, we know that there is a  cone $Z_x$
containing $x_{j_0}$ (the center of $B_{j_0}$) and such that every
point of  $K\cap 30B_{j_0}$ is at a distance less than
$\varepsilon_0 120Ur_{j_0}$ from $Z_x$ where $r_{j_0}$ is the
radius of $B_{j_0}$. We also know that $Z_x$ is almost centered in
$B_{j_0}$. Thus $K$ is at a distance less than $400U\varepsilon_0
r_{j}$ from $Z_x$  in all the $W_j$ which contains $x$. Therefore,
if we consider the connected components of $\R^3 \backslash Z_x$,
each one contains one and only one $D^j_k$ for all $j\in J_x$.
Thus we can define for all $k$, a polyhedral domain $D^x_k$ that
contains each $D^j_k$ for all $j \in J_x$. For instance we could
define $D^x_k$ as being the smaller convex domain containing all
the $D^j_k$ for $j \in J_x$. In addition the diameter of $D^x_k$
is equivalent to   $\frac{C}{U}\delta(x)$ because all the balls
$W_j$ for $j$ in $J_x$  have a radius equivalent to
$\frac{C}{U}\delta(x)$. Here the constant $C$ is just a geometric
constant.
  Finally, $D^x_k$ does not meet $K$, and, as $\varepsilon_0$ is small enough,  we also
   have
  $$D^x_k\subset \bigcup_{j \in J}10W_j\backslash 6 W_j \subset \mathpzc{Z}.$$

\begin{center}
\input{dx.txt}

 \nopagebreak[4]
Picture of the situation when $Z_x$ is a  hyperplane.
\end{center}

\begin{center}
\input{dxy.tex}

 \nopagebreak[4]
Picture of the situation when $Z_x$ is a $\Y$.
\end{center}

  Let $m^x_k$ be the average of $u$ on $D^x_k$.   Since the $\theta_j$ are a partition of
  unity, we know that $\sum \nabla \theta_j =0$ hence we can subtract   $m^x_k$
  and then
\begin{eqnarray}
\nabla v_k(x)= \theta_0\nabla u(x)+\nabla
\theta_0(u-m^x_k)+\sum_{j>0}(m^j_k-m^x_k)\nabla \theta_j (x).
\label{somme}
\end{eqnarray}
On the other hand using Poincaré inequality (with  $C$ a
dimensional constant),
$$|m^j_k-m^x_k|\leq C\frac{1}{r_j^3}\int_{D^j_k}|u(y)-m^x_k|dy\leq \frac{1}{r_j^3}\int_{D^x_k}|u(y)-m^x_k|dy\leq C
\frac{1}{r_j^{2}}\int_{D^x_k}|\nabla u(y)|dy.$$ In addition all
the $r_j$ for $j\in J_x$ are equivalent, and since $|\nabla
\theta_j|\leq C r_j^{-1}$,  every term in the sum in \eqref{somme}
is bounded by $C\frac{1}{r_j^3}\int_{D^x_k}|\nabla u(y)|dy$.
Furthermore, $r_j$ is bounded from below by $C\delta(x)$ then
finally since the sum has only $C_2$ terms,
\begin{eqnarray}
|\nabla v_k(x)|\leq |\nabla u (x) \theta_0(x)|+ |\nabla
\theta_0(x)
||u(x)-m^x_k|+C\frac{1}{\delta(x)^3}\int_{D^x_k}|\nabla u(y)|dy.
\label{estim1}
\end{eqnarray}
We will show later that when $x \in \mathpzc{V}\backslash
\mathpzc{V}_\rho$,
\begin{eqnarray}
|u(x)-m^x_k|\leq C\frac{1}{\delta(x)^{2}}\int_{D^x_k}|\nabla u
(y)|dy. \label{rfinir}
\end{eqnarray}
It follows that
\begin{eqnarray}
|\nabla v_k(x)|&\leq& \underbrace{|\nabla u(x)
|\theta_0(x)}_{f_1}+
\underbrace{C\frac{1}{\delta(x)^3}\int_{D^x_k}|\nabla
u(y)|dy}_{f_2} \notag
\end{eqnarray}
Therefore
$$|\nabla v_k(x)|^2\leq f_1^2+ f_2^2+2f_1f_2 \leq 2(f_1^2+f_2^2). $$
Set
$$\Theta:=(\bigcup_{j \in J} 10W_j)\backslash (\bigcup_{j \in J} 8W_j)\subset \mathpzc{Z}$$
so that $\nabla \theta_0$ has its support in $\Theta$. It will be
convenient to also define
$$V:=\bigcup_{j\in J}10 W_j \subset \mathpzc{V}$$
and
$$Q_k:=\bigcup_{x \in V}D^x_k \subset \mathpzc{Z}.$$
Taking the integral on $x \in V$ immediately give for $f_1$
$$\int_{V}f_1^2(x)dx\leq \int_{\Theta} |\nabla u(x)|^2 dx$$
because $\theta_0=0$ on $V\backslash \Theta$.  For $f_2$, applying
Hölder's inequality we obtain
$$f_2^2\leq C \frac{1}{\delta(x)^{3}}\int_{D^x_k}|\nabla u(y)|^2dy . $$
Then integrating and applying Fubini leads to
\begin{eqnarray}
\int_{V}f_2^2(x)dx &\leq&  C\int_{V} \left(\frac{1}{\delta(x)^{3}}
\int_{D^x_k}\textbf{1}_{Q_k}|\nabla u(y)|^2 dy\right) dx \notag \\
&\leq&  C \int_{Q_k}|\nabla u(y)|^2 \left( \int_{\{x; y\in
D^x_k\}} \delta(x)^{-3} dx\right)dy. \notag
\end{eqnarray}
The point is now to show that
\begin{eqnarray}
 \int_{\{x; y\in D^x_k\}} \delta(x)^{-3} dx \leq C
 \label{rpresquefini}.
\end{eqnarray}
Here $C$ is a dimensional constant. As $x$ is fixed, $D^x_k$ is a
polyedral domain of diameter $\frac{C}{U}\delta(x)$ and  at a
distance less than $\frac{C}{U}\delta(x)$ from $x$. Hence $D^x_k$
is included in a ball $B(x, \frac{C}{U}\delta(x))$. Therefore, the
set $\{x; y \in D^x_k \}$ is included in $A:=\{x; d(x,y)\leq
\frac{C}{U}\delta(x)\}$. We want to show that for all $x$ in $A$,
$\delta(x)$ is equivalent to $\delta(y)$. This is where the choice
of constant $U$ is important. Indeed, if $x$ is in $A$ then since
$\delta$ is $C_0$-lipschitz,
\begin{eqnarray}
 \delta(x)\leq  \delta(y)+C_0|x-y| \notag
\end{eqnarray}
hence
\begin{eqnarray}
 \delta(y)\geq \delta(x)-\frac{CC_0}{U}\delta(x).  \label{rdy}
\end{eqnarray}
Recall that   $U\geq 30 C_0$ is big as we want compared to $C_0$.
Thus if $U$ is chosen large enough, $\frac{CC_0}{U}$ is less than
$1$ and then \eqref{rdy} gives
$$\delta(x)\leq C\delta(y).$$
On the other hand,
\begin{eqnarray}
\delta(y) \leq \delta(x)+C_0|x-y|\leq C \delta(x) .\notag
\end{eqnarray}
Thus $\delta(x)$ and $\delta(y)$ are equivalent on $A$ and
$$ \int_{\{x; y\in D^x_k\}} \delta(x)^{-3} dx \leq C \delta(y)^{-3}\int_{\{x; d(x,y)\leq C\delta(y)\}} dx\leq C$$
hence \eqref{rpresquefini} is true and finally
\begin{eqnarray}
\int_{V\backslash \mathpzc{V}_\rho}|\nabla v_k(x)|^2dx &\leq&
\int_{\Theta}|\nabla u(x)|^2dx + C \int_{Q_k}|\nabla u(y)|^2
dy\notag
\end{eqnarray}
thus, since $\Theta \cup Q_k  \subset \mathpzc{Z}$,
\begin{eqnarray}
 \int_{V\backslash \mathpzc{V}_\rho}|\nabla v_k(x)|^2dx &\leq &C\int_{\mathpzc{Z}}|\nabla
 u(x)|^2dx\notag
\end{eqnarray}
and this implies
\begin{eqnarray}
\int_{\Delta_k\backslash \mathpzc{V}_\rho}|\nabla v_k|^2 &\leq&
\int_{\Delta_k\backslash (V\cup \mathpzc{V}_\rho)}|\nabla u|^2
+\int_{V\backslash \mathpzc{V}_\rho}|\nabla v_k|^2
\notag \\
&\leq &\int_{\Delta_k\backslash \mathpzc{V(\frac{1}{3})}}|\nabla
u|^2 +C\int_{\mathpzc{Z}}|\nabla u|^2. \notag
\end{eqnarray}
So to finish the proof we have now to show \eqref{rfinir}. In fact
we want to estimate
$$|u(x)-m^x_k|\nabla \theta_0(x).$$
We may assume that  $x$ is in the support of $\nabla \theta_0$,
hence that $x \in \Theta$. We will use the mean value theorem on
$u$. For this, it is fundamental that for all  $x\in \Theta$ and
for all $y \in D^x_k$, the segment $[x,y]$ does not meet the
singular set $K$.  Consequently, the next estimate is not true for
all $x\in \Theta$, but it is true for all $x \in \Theta_k
\backslash \mathpzc{V}_{\rho}$ where $\Theta_k$ is the connected
component of $\Theta \backslash K$ that is included in
$\Omega_k(x_0,r_0)$. Indeed, if $10W_{j_0}$ does not meet
$\partial B(x_0,\rho)$, and if $J_{j_0}$ is the set of indices $j$
such that $W_j$ meets $10W_{j_0}$, we claim that
\begin{eqnarray}
K \cap 10W_{j_0} \subset K \cap \bigcup_{j \in
J_{j_0}}\frac{1}{10}W_j.
\end{eqnarray}
To see this, we denote $W_{j_0}:=B(x_{j_0},r_{j_0})$ and we use
the fact that $\delta$ is $C_0$-lipschitz. Thus for all $x \in K
\cap 10W_{j_0} $ we have
$$\delta(x)\geq \delta(x_{j_0})-10C_0r_{j_0}\geq Ur_{j_0}-10C_0r_{j_0}\geq 20C_0r_{j_0}.$$
because remember that $U\geq 30 C_0$. Then by Lemma \ref{lemwhit}
$ii)$ we know that
\begin{eqnarray}
\bigcup_{j \in J_{j_0}}\frac{6}{10}W_j \supseteq  \bigcup_{x \in
K\cap 10 W_{j_0}}B(x,\frac{3}{200 U}\delta(x)) \supseteq
\bigcup_{x \in K\cap 10 W_{j_0}}B(x,\frac{3C_0}{10 U}r_{j_0})
\label{roro}
\end{eqnarray}
Note that \eqref{roro} is not true if $W_{j_0}\cap \partial
B(x_0,\rho)\not = \emptyset$ because if $x \in 8 W_{j_0}\backslash
B(x_0,\rho)$, $x$ does not belong to any of the $\frac{6}{10}W_j$.
Denote by $J_x$ the set of indices such that $10W_j$ contains
$x\in \Theta_k \backslash \mathpzc{V}_{\rho}$. By doing the same
for all the $W_j$ with $j \in J_x$ it follows from $\eqref{roro}$
that we have an universal and small constant $c$ such that
$$K \cap \bigcup_{j \in J_x} 10W_{j}\subset \bigcup_{y\in K \cap \bigcup_{j \in J_x}10W_{j}}B(y,c\delta(x))
\subset \bigcup_{j \in J_x}(\frac{6}{10})W_j.$$
 Moreover, if
$\varepsilon_0$ is small enough, the smaller convex domain that
contains $D_x^k$ and $10W_j\backslash \bigcup_{j \in
J_x}\frac{6}{10}W_j$
 does not meet $K$ as
well. This is because all the $W_j$ for $j \in J_x$ are almost
``aligned'' on a same minimal cone $Z$. By passing to the
complement, and by the fact that $\Theta_k \cap \bigcup_{j \in
J}(\frac{6}{10})W_j= \emptyset$ we obtain that if $x \in \Theta_k
\backslash \mathpzc{V}_{\rho}$ and $y \in D^x_k$, the segment
$[x,y]$ does not meet $K$.
\begin{center}
\includegraphics[width=13cm]{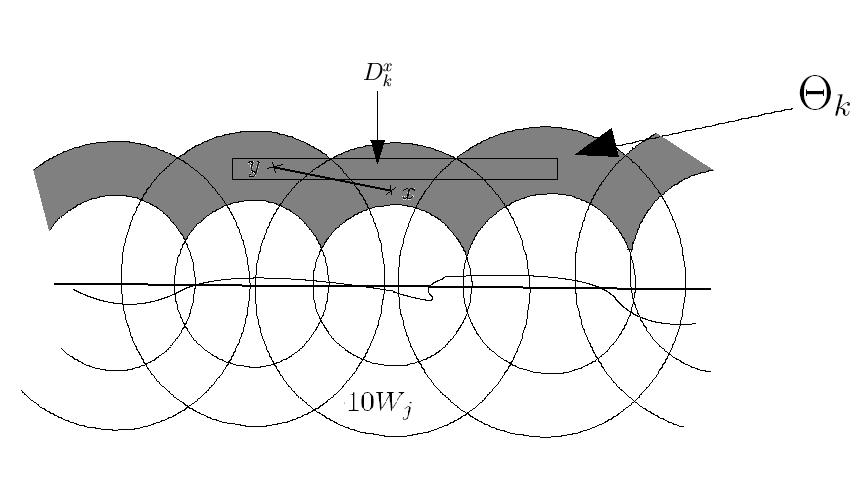}

 \nopagebreak[4]
Figure 3
\end{center}

 Thus the following inequality is true
(firstly we suppose that $u$ is $C^\infty$ in $D^x_k$ and then we
use a density argument)

\begin{eqnarray}
|u(x)-m^x_k|&\leq& \frac{1}{|D^x_k|}\int_{D^x_k}|u(x)-u(y)|dy \notag \\
&\leq & \frac{1}{|D^x_k|}\int_{D^x_k}\int_{[0,1]}|x-y||\nabla u(x+t(y-x))|dtdy \notag \\
&\leq & \frac{1}{|D^x_k|}\int_{D^x_k}\int_{D^x_k} |\nabla u_n(z)|\frac{dz}{|x-y|^{2}}dy \notag \\
&\leq& C C_1\frac{1}{\delta(x)^{2}}\int_{D^x_k}|\nabla u(z)|dz
\notag
\end{eqnarray}

because $\frac{1}{|x-y|^2} \in L^{1}_{loc}(\R^3)$ and
$\int_{B(x,R)}\frac{1}{|x-y|^2}\simeq R$. This ends the
construction of $v^k$, and the proof of Lemma \ref{extentionW} is
achieved. \qed


\section{Decay of energy in the complement of
$(\varepsilon_0,\varepsilon)$-minimal sets}
\label{section4proff}
We can now prove the main result.

 {\bf Proof of Theorem \ref{decreifenberg}:} The argument is by contradiction and compactness. The beginning
 is very close to Theorem 8.19 of \cite{afp}. The main changes
 will come when we will have to prove Inequality
 \eqref{lemmefonda}.

We suppose that Theorem \ref{decreifenberg} is false. Then there
is an $i_0 \in \{1,2,3\}$, an $\varepsilon_0 <10^{-5}$, a
$\gamma<0,8$ and there is a positive radius $r<\frac{1}{2}$ such
that for all $\varepsilon > 0$ there is a triplet
$(u_\varepsilon,E_{\varepsilon},\rho_\varepsilon)\in
\Lambda_{i_0}(\varepsilon)$  such that
\begin{eqnarray}
\int_{B(0,r)\backslash E^\rho_{\varepsilon}}|\nabla u_{\varepsilon}|^2
> r^{2+\gamma}. \label{nondec}
\end{eqnarray}

In the following, $i_0$ is a fixed index. We call $Z^0$ a minimal
cone of type $i_0$ centered at $0$. By rotation invariance we
suppose that all the sets $E_\varepsilon$ are included in $\{y \in
B;d(y,Z^0)\leq \varepsilon\}$.

For every positive number $a$ such that $0<a<10^{-5}$ we call
$Z^0_a$ the region
$$Z^0_a:=\{y \in B; d(x,Z^0)\leq a\}.$$
Then for $j \in [1,i_0+1]\cap \N$ we denote by $A_j$ the different
connected components of $B\backslash Z^0_a$ and
$m_{a,\varepsilon}^j$ the mean value of  $u_\varepsilon$ on
$A_a^j$. Since for all $\varepsilon$ we have
$$\int_{B\backslash E^\rho_\varepsilon}|\nabla u_{\varepsilon}|^2=1,$$
the Poincaré inequality says that for all $0<a<10^{-5}$ the
sequences $u_\varepsilon- m_{10^{-5},\varepsilon}^j$ are uniformly
bounded for the norm $W^{1,2}(A_a^{j})$. We also know that the
sequence $\nabla u_\varepsilon$ is uniformly bounded in $L^2(B)$.
Therefore, possibly extracting a further subsequence obtained by a
diagonal argument, we may conclude that there is a subsequence
$u_n$ that weakly converges in $W^{1,2}$ on all the $A_a^{j}$ to a
function $u\in W^{1,2}(B\backslash Z^0)$, and such that $\nabla
u_n$ also weakly converges in $L^2(B)$ to a certain function,
which by uniqueness of limits in the distributional space is equal
to $\nabla u$. We denote by $E_n$ the set associated to $u_n$, and
also let $\varepsilon_n$ and $\rho_n$ be such that
 $(u_n,E_n,\rho_n)\in \Lambda(\varepsilon_n)$. In addition, passing if
 necessary to a subsequence, we can assume that
$\rho_n$ converges to some real $\rho_\infty \in
[\frac{1}{2},\frac{3}{4}]$.

 We claim that $u$ is an energy minimizing function in $B\backslash
 Z^0$, or in other words, $u$ is harmonic in $B\backslash Z^0$ and
 have zero normal derivative on both sides of $Z^0$. Let $\varphi$ be
 a $C^\infty$ function with compact support in
$B\backslash Z^0$. Since for all $n$ we know that $u_n$ is an
energy minimizer we have that
$$\int_{B}\langle \nabla \varphi, \nabla u_n\rangle=0.$$
Then passing to the limit  (because $\nabla u_n$ weakly converges
in $L^2(B)$) we obtain
$$\int_{B\backslash Z^0}\langle \nabla \varphi , \nabla u\rangle =0$$
and this proves that $u$ is harmonic in $B\backslash Z^0$.

Now we want to show that $u$ has zero normal derivative on $Z^0$.
So we have to show that for all $\varphi \in C^1_0(B)$,
$$\int_{\Omega^k}\langle \nabla \varphi , \nabla u\rangle =0$$
where $\Omega^k$ are connected components of $B\backslash Z^0$.
For all $n$ we know that
$$E_n
\subset Z^0_{\varepsilon_n}:=\{y \in B ; d(y, Z^0)\leq
\varepsilon_n\}$$ and that $E_n$ is separating in $B$. Recall that
$A^k_{\varepsilon_n}$ denotes the connected components of
$B\backslash Z^0_{\varepsilon_n}$ and we call $\Omega^k_n$ the
connected components of $B\backslash E_n$ that contains
$A^k_{\varepsilon_n}$. Let $k_0$ be fixed and consider a sequence
$c_n$ that converges to $+\infty$. Denote by $v_n$ the sequence in
$W^{1,2}(B \backslash K)$ defined by
$$v_n=\left\{
\begin{array}{cl}
u_n+c_n & \text{ in } \Omega^{k_0}_n\\
u_n-c_n & \text{ in } \Omega^k_n \text{ for  } k \not = k_0 \\
u_n & \text{ in } B \backslash \bigcup_{k}\Omega^k_n
\end{array}\right. .
$$
Now let $\eta>0$ be fixed, and let $\chi(t)$ be a function on $\R$
such that for all $t$, $\chi'(t)\leq \eta$ and such that $\lim_{t
\to - \infty}\chi(t)=0$ and $\lim_{t \to + \infty}\chi(t)=1$.  If
$\varphi$ is a function in $C^1_0(B)$, the function $\varphi
\chi(v_n) \in W^{1,2}(B \backslash E_n)$. Comparing energy of
$u_n$ with energy of $u_n+\varepsilon_n\varphi \chi(v_n)$ we
obtain that
$$\liminf_{n\to \infty}\left[\int_{B}\varphi\langle  \nabla u_n, \nabla \chi(v_n) \rangle
+\int_{B} \chi(v_n)\langle  \nabla u_n, \nabla \varphi
\rangle\right]\geq 0.$$

On the other hand, $\nabla \chi(v_n)=\chi'(v_n)\nabla u_n$,
$|\chi'|\leq \eta$ and $\int_{B}|\nabla u_n|^2=1$. Moreover,
$\chi(v_n)$ strongly converges in $L^{2}(B)$ to $\1_\Omega^{k_0}$.
We deduce
$$\int_{\Omega^{k_0}}\langle\nabla u, \nabla \varphi \rangle \geq -\eta \|\varphi\|_\infty$$
and we conclude by the fact that $\eta$ is an arbitrary constant.
We do the same for all $\Omega^k$ and that proves that $u$ is an
energy minimiser in $B\backslash Z^0$. Hence, we know by the
preceding section that for all $r\in ]0,1]$,
\begin{eqnarray}
\frac{1}{r^2}\int_{B(0,r)\backslash Z^0}|\nabla u|^2\leq
r^{0,8}\int_{B(0,1)\backslash Z^0}|\nabla u|^2=r^{0,8}.
\label{dec}
\end{eqnarray}

 Therefore, if we show that  $\nabla u_n$ strongly converges in
 $L^2(B(0,r))$, the contradiction will follow by passing to the
 limit in \eqref{nondec}. The purpose of all the following of the proof is to justify this convergence.

 So we consider the measures $\mu_n:=|\nabla u_n|^2dx$ on $B$. Since for all  $n$
 we have $\mu_n(B)=1$, passing if necessary to a subsequence we may assume that $\mu_n$
weakly converge to some measure $\mu$. Since $\nabla u_n$ weakly
converge in $L^2$ to $\nabla u$, all we can say is that (see for
instance proposition 1.62(b) of \cite{afp})
\begin{eqnarray}
|\nabla u|^2dx\leq \mu . \label{mu}
\end{eqnarray}
So we have to show that \eqref{mu} is an equality in $B(0,r)$.
This would be enough because the Radon-Riesz Theorem says that if
a sequence $f_n$ of $L^p$ functions weakly converge to some
function $f$, and if in addition the sequence of $L^p$ norms
$\|f_n\|_p$ converge to $\|f\|_p$ then $f_n$ strongly converge to
$f$ in $L^p$.

First of all, let us show that the regular part of $\mu$ with
respect to Lebesgue measure is $|\nabla u|^2dx$ and that the
singular part is concentred on $Z^0$. To show that, let us recall
that  $u_n$ is a sequence of harmonic functions with $L^2$ norms
uniformly bounded and then are uniformly bounded for the
$L^\infty$ norm on all compact sets of $B\backslash Z^0$. Thus by
covering $B\backslash Z^0$ with a countably union of compact sets
and by using a diagonal argument we can say that $u_n$ converges
to $u$ for the $L^p$ norm on all the compact sets and for all $p$.
Moreover since the $u_n$ are harmonic we know that their gradients
are also uniformly bounded. Then extracting a further subsequence
we may assume that $\nabla u_n$ converges to $\nabla u$ in $L^2$
strongly on all compact sets of $B\backslash Z^0$. It follows that
for every compact set $U$ in $B\backslash Z^0$,
$$\mu(U)=\lim_{n \to +\infty}\mu_n(U)=\lim_{n \to +\infty}\int_{U}|\nabla u_n|^2dx= \int_{U}|\nabla u|^2dx.$$

Consequently, to prove that $\mu=|\nabla u|^2dx$ everywhere, we
have to show that $\mu( Z^0)=0$. So we consider a domain
$\mathcal{C}_h$ containing $K$. We define the base
$$G:=B(0,\rho_\infty)\cap Z^0$$
and then set
$$\mathcal{C}_h:= \{y \in B ; d(y,G) \leq h \}.$$
Define also
$$\mathcal{C}^n_h:=\{y \in B ; d(y,G_n) \leq h \}$$
based on
$$G_n:=B(0,\rho_n)\cap Z^0$$
such that $\mathcal{C}_h^n$ converges to $\mathcal{C}_h$ for the
Hausdorff distance. The strategy is to estimate
$\mu(\mathcal{C}_h)$ and then let $h$ tend to $0$.

Let us sketch the next ideas that will be used to finish the
proof. The main point is to prove an inequality of the following
type
\begin{eqnarray}
\int_{\mathcal{C}_h}|\nabla u_n|^2dx \leq C \int_{S_n}|\nabla
u_n|^2 dx\label{lemmefonda}
\end{eqnarray}
where $S_n$ is a closed set such that $E^{\rho_n}_n \cap
B(0,\rho_n)\subset B(0,\rho_n)\backslash S_n$ and $S_n$ converges
to
$$S_\infty:=\{x\in \mathcal{C}_h ; d(x,Z^0)\geq C_5d(x, \partial \mathcal{C}_h)  \}.$$
If  \eqref{lemmefonda} is true, then passing to the limit (that
will be justify later),
$$\mu(\mathcal{C}_h)\leq C \mu(S_\infty)=\int_{S_\infty}|\nabla u|^2 dx \leq
Ch\|\nabla u\|^2_\infty$$ and then letting $h$ tend to $0$ we
conclude that  $\mu(Z^0\cap B(0, \rho_\infty))=0$.

Now the purpose of the following is to justify the above
arguments. To prove \eqref{lemmefonda} we will use the Whitney
extension of preceding section. So let $n$ be fixed. We have to
construct a function $\delta$ in order to apply Lemma
\ref{extentionW} to $u_n$ and $E_n$.

Let $\{B_i\}$ be the family of balls in the definition of
$(\varepsilon_0, \varepsilon_n)$-minimality of $E_n$. For all $i
\in I$ denote by $\psi_i$ a positive function, 1-lipschitz, such
that
$$\psi_i(x):=
\left\{
\begin{array}{c}
r_i \text{ in } B_i  \\
0 \text{ out of  } 2B_i
\end{array}
\right. .
$$
Then for all $x \in B$, let $d(x)$ be defined by
$$d(x):=\sum_{i \in I}\psi_i(x).$$
The sum is locally finite so  $d(x)$ is $C_0$-Lipschitz where
$C_0$ is the constant from the bounded cover propriety of
$\{2B_i\}_{i\in I}$.

Then for all $x\in B(0,1)$ we define
$$\delta(x):=\max(d(x,\partial \mathcal{C}_h^n),d(x)).$$

By construction and by definition of
$(\varepsilon_0,\varepsilon_n)$-minimality of $E_n$, we can deduce
that $E_n$ is satisfying Hypothesis $\mathcal{H}$ in $B(0,1)$ with
geometric function $\delta$, constant $\varepsilon_0$ and radius
$\rho_n$. Therefore applying Lemma \ref{extentionW}, we obtain for
all $k \in \N \cap [1, \mathpzc{k}^{B(0, 1)}]$, a function $v_k
\in W^{1,2}(B(0,\rho_n)\cap(\Omega_k\cup \mathpzc{V}^n) \backslash
\mathpzc{V}^n_{\rho_n})$ with
$$\mathpzc{V}^n:= \bigcup_{x\in K \cap B(0,\rho_n)}B(x,\frac{10}{U}\delta(x))$$
and
$$\mathpzc{V}^n_{\rho_n}:=\bigcup_{x; B(x,\frac{10}{U}\delta(x))\cap \partial B(0,\rho_n)\not = \emptyset}B(x,\frac{10}{U}\delta(x))$$
and where $U$ is a constant satisfying $U\geq 30 C_0$.

By construction, the $v^k$ are equal to $u_n$ in $B(0,1)\backslash
 \mathpzc{V}^n$. Moreover, if $\varepsilon$ is small enough,
$(\mathpzc{V}\backslash \mathpzc{V}_{\rho_n})\cap
B(0,\rho_n)\subset \mathcal{C}_h^n$ and note that we can fix $C_1$
in such a way that
\begin{eqnarray}
\mathpzc{V}_{\rho_n}^n\subset \bigcup_{i \in I_{\rho_n}}C_1B_i.
\label{boundary}
\end{eqnarray}
see \eqref{defbpr} for the definition of $C_1$. Indeed, let $x$ be
a point in $\mathpzc{V}_{\rho_n}^n$ and $y$ be such that $B(y ,
\frac{10}{U}\delta(y)) \cap
\partial B(0,\rho_n)\not = \emptyset $ and such that $x \in B(y,
\frac{10}{U}\delta(y))$. Set $r_y:=\frac{10}{U}\delta(y)$. Since
$\delta(y)=\max(d(y),d(y,\partial \mathcal{C}_h^n))$ and since
$U\geq 30C_0$, we can deduce that $\delta(y)=d(y):=\sum_{i \in
I}\psi_i(y)$. Therefore there is an index $i_0$ such that
$\psi_{i_{0}}(y)$ is not equal to zero. This implies that $y\in 2
B_{i_0}$. We can also suppose that $i_0$ is the index for which
$\psi_{i_0}(y)$ is the maximum of all the $\psi_i(y)$. On the
other hand, $r_y =\frac{1}{U}\sum_{i ; y \in 2B_i}\psi_{i}(y)\leq
\frac{C_0}{U}\psi_{i_0}(y)\leq \frac{C_0}{U}r_{i_0}$ (because the
$\{2B_i\}$ are in bounded cover with constant $C_0$). So we have
proved
$$B(y,r_y) \subset (2+10\frac{C_0}{U})B_{i_0}$$
and  \eqref{boundary}  follows if we choose
$C_1=(2+10\frac{C_0}{U})$. Lemma \ref{extentionW} gives
$$\int_{\mathpzc{V}^n\backslash \mathpzc{V}^n_{\rho_n}}|\nabla v^k|^2dx\leq
C\int_{\mathpzc{Z^n}}|\nabla u|^2dx.$$ with $\mathpzc{Z}^n$
defined by \eqref{defQ2}. Therefore we can compare the energy of
$u_n$ with the energy of $v$
 defined by
$$
v= v^k \text{  in } B(0,\rho_n)\cap \Omega_k.
$$
On other components we can assume that $v=0$. Recall that $u_n$ is
by definition the energy minimizing function in $B\backslash
E_n^{\rho_n}$ and note that  $v$ is a competitor in $B(0,\rho_n) $
that is equal to $u_n$ out of $\mathpzc{V}^n$ thus
\begin{eqnarray}
\int_{\mathpzc{V}^n\backslash E_n^{\rho_n}}|\nabla u_n|^2&\leq &
 \sum_{k=1}^{\mathpzc{k}} \int_{\mathpzc{V}^n\backslash \mathpzc{V}^n_{\rho_n}}|\nabla v^k|^2   \leq C\int_{ \mathpzc{Z}^n}|\nabla u_n|^2 .\notag
\end{eqnarray}

Hence

\begin{eqnarray}
\int_{\mathcal{C}_h}|\nabla u_n|^2 &\leq& \int_{\mathcal{C}_h
\backslash \mathpzc{V}^n}|\nabla u
_n|^2+\int_{\mathpzc{V}^n}|\nabla u_n|^2 \notag \\
&\leq& \int_{\mathcal{C}_h \backslash \mathpzc{V}^n}|\nabla u
_n|^2+C\int_{\mathpzc{Z}^n}|\nabla u_n|^2. \label{lemmefonda2}
\end{eqnarray}

Now we define the closed set
$$S_n:=\mathcal{C}_h\backslash \overset{\circ}{\mathpzc{V}^n }\cup \mathpzc{Z}^n.$$

Then we have proved
\begin{eqnarray}
\int_{\mathcal{C}_h}|\nabla u_n|^2 \leq C\int_{S_n}|\nabla u _n|^2
.\label{ouff}
\end{eqnarray}

The point is now to show that
\begin{eqnarray} S_n \text{ converges
to } S_\infty :=\{x \in \mathcal{C}_h ; d(x,Z^0)\geq C_5
d(x,\partial \mathcal{C}_h)\} \label{convtn}
\end{eqnarray}
and
\begin{eqnarray} \limsup_{n \to +\infty}\bigg(\int_{S_n}|\nabla u_n|^2dx\bigg)\leq
\mu(S_\infty).\label{convitn}
\end{eqnarray}
Here $C_5$ is a constant that will be chosen later.

For all  $\eta>0$ we denote by $S_\infty^{\eta} $ a
$\eta$-neighborhood of $S_\infty$. In other words
$$S_\infty^{\eta}:= \{x\in B; d(x,S_\infty)< \eta \}.$$ We want to
show that for $n$ big enough and if the constant $C_5$ is chosen
properly, then all the $S_n$ are included in $S_\infty^{\eta}$. So
let $n_0$ be fixed in  such a way that for all $n\geq n_0$ we have
$\varepsilon_n <\eta $ and $|\rho-\rho_n|\leq \eta $. Hence, for
all $n\geq n_0$ and for all $x\in E_n$ we have $d(x,Z^0)\leq
\eta$. Let now  $x$  be a point in $S_n$ for $n \geq n_0$. We
claim that $x$ is in $S_\infty^{C_6\eta}$ for a certain constant
$C_6$. Indeed, if $x$ is in $S_n$, then by definition of
$\mathpzc{Z}^n$ and $\mathpzc{V}^n$ (see \eqref{defQ1} and
\eqref{defQ2}) $x$ does not belong to any of the $B(y,
\frac{1}{10U}\delta(y))$. In other words,
$$x \in \mathcal{C}_h\backslash \mathpzc{V}(1/10)_n$$
with
$$\mathpzc{V}(t)_n:=\bigcup_{y \in E_n\cap B(0,\rho_n)}B(y,\frac{t}{U}\delta(y)).$$
Denote by
 $$\mathpzc{V}(1/10)_n^\eta:=\{y \in \mathcal{C}_h; d(y,\mathpzc{V}(1/10)_n)\leq
C_6\eta\}$$ with $C_6$ a constant to be fixed later. We want to
show that there is a constant $c$ (depending on $C_0$ and $U$)
such that
\begin{eqnarray}
\mathpzc{V}(1/10)_n^\eta \supseteq \{y \in \mathcal {C}_h;
d(y,Z^0)\leq cd(y,\partial \mathcal {C}_h)+C_6\eta\}.
\label{amontrer2}
\end{eqnarray}
Let $y\in \mathcal {C}_h$ be such that $d(y,Z^0)\leq cd(y,\partial
\mathcal {C}_h)$. We call $\bar y \in Z^0$ a point such that
$$d(y,\bar y)=d(y,Z^0).$$
Since $E_n$ is separating and since $E_n \subset
Z^0_{\varepsilon_n}$, there is a point $z \in E_n$ such that
$$d(\bar y, z)\leq \eta$$
otherwise the ball $B(\bar y,\varepsilon_n)$ would not contain any
point  of $E_n$.  Since $\delta$ is $C_0$-lipschitz,
\begin{eqnarray} \delta(y)&\leq&
\delta(z)+C_0|y-z| \notag \\
&\leq& \delta(z)+C_0(|y-\bar y|+|\bar y -z|)\notag \\
&\leq & \delta(z)+C_0(cd(y,\partial \mathcal {C}_h)+\eta).
\end{eqnarray}
Recall that
$$\delta(y)=\max(d(y,\partial
\mathcal{C}_h^n),d(x))\geq d(y,\partial \mathcal{C}_h^n) \geq
d(y,\partial \mathcal{C}_h)-\eta.$$ Hence,
$$(1-cC_0)d(y,\partial \mathcal{C}_h)-(1+C_0)\eta\leq \delta(z).$$
Now,
\begin{eqnarray}
|y-z|&\leq&|y-\bar y| +|\bar y-z| \notag \\
&\leq &c d(y ,\partial \mathcal{C}_h)+\eta \notag \\
&\leq& c(1-cC_0)^{-1}[\delta(z)+(1+C_0)\eta]+\eta. \notag
\end{eqnarray}
Then we can choose  $c<1$ small enough in order to have
$$|y-z|\leq \frac{1}{10U}\delta(z)+C_6\eta$$
where $C_6$ is depending on $c$. This implies that $y\in
B(z,\frac{1}{10U}\delta(z)+C_6 \eta)$, thus  \eqref{amontrer2} is
proved. Then by passing \eqref{amontrer2} to the complementary we
deduce
$$\mathcal{C}_h\backslash \mathpzc{V}(1/10)_n \subset
\{y \in \mathcal{C}_h; d(y,Z_0) \geq C_5d(y,\partial
\mathcal{C}_h)+C_6\eta\}$$ with $C_5$ and  $C_6$ depending on $U$
and $C_0$. Thus $x \in S_\infty^{C_6\eta}$ and this proves
\eqref{convtn}.

It is time now to prove \eqref{convitn} and to finish the proof of
the theorem. We keep the same notations
 as before. That is for all $\eta$, we know that there is an
 integer  $n_0$ such that for all $n \geq n_0$, $S_n\subset
 S_\infty^{C_6\eta}$. By \eqref{ouff} we know  that
$$\int_{\mathcal{C}_h}|\nabla u_n|^2dx\leq C\int_{S_n}|\nabla u_n|^2dx. $$
In other words if $\mu_n$ is the measure $|\nabla u_n|^2dx$ we
have
$$\mu_n(\mathcal{C}_h)\leq C \mu_n(S_n).$$
On the other hand, since $\mu_n$ weakly converges to $\mu$ and
since $\mathcal{C}_h$ is a fixed closed set, we can deduce that
$\mu_n(\mathcal{C}_h)$ converges to $\mu(\mathcal{C}_h)$. Now
since $S_n$ is included in $S_\infty^{\eta}$ for $n$ big enough,
$$\mu(\mathcal{C}_h)\leq C \liminf_{n \to \infty}\mu_n(S_\infty^{\eta})$$
thus since now $S_\infty^{\eta}$ does not depend on $n$, we can
apply the weak convergence and get
$$\mu(\mathcal{C}_h)\leq C \mu(S_\infty^{\eta}).$$
Finally, letting  $\eta$ tend to $0$ we obtain
\begin{eqnarray}
\mu(\mathcal{C}_h) \leq C\mu(S_\infty). \label{egalfonda}
\end{eqnarray}

This last limit is justified because $S_\infty$ is the decreasing
intersection of $S_\infty^\eta$.

To finish the proof of the theorem we use the fact that $\mu$ is
regular with respect to the Lebesgue measure on $S_\infty$ so
\eqref{egalfonda} gives
$$\mu(\mathcal{C}_{h})\leq C \int_{S_\infty}|\nabla u|^2 dx\leq Ch\|\nabla u\|_\infty^2.$$
We know that $\nabla u$ is bounded on $S_\infty$ because
 $u$ is an harmonic function with Neumann condition on $Z^0$.

Now by letting $h$ tend to $0$ we obtain that $\mu(Z^0\cap
B(0,\rho_\infty))=0$ and this gives that $\nabla u_n$ strongly
converges to $\nabla u$ in $L^2(B(0,\rho_\infty))$. So passing to
the limit in \eqref{nondec} we get a contradiction with
\eqref{dec} and this proves the theorem.\qed





\bibliographystyle{alpha}
\bibliography{biblio}

ADDRESS :

Antoine LEMENANT \\
e-mail : antoine.lemenant@math.u-psud.fr

Université Paris XI\\
Bureau 15 Bâtiment 430 \\
ORSAY 91400 FRANCE

Tél: 00 33 169157951

\end{document}

%% file: cex.tex
\unitlength 1mm
\begin{picture}(81.20,45.30)(0,0)

\linethickness{0.15mm}
\put(61.86,24.58){\line(0,1){0.94}}
\multiput(61.81,26.45)(0.05,-0.93){1}{\line(0,-1){0.93}}
\multiput(61.72,27.38)(0.09,-0.93){1}{\line(0,-1){0.93}}
\multiput(61.59,28.31)(0.14,-0.93){1}{\line(0,-1){0.93}}
\multiput(61.40,29.23)(0.09,-0.46){2}{\line(0,-1){0.46}}
\multiput(61.18,30.14)(0.11,-0.46){2}{\line(0,-1){0.46}}
\multiput(60.90,31.04)(0.14,-0.45){2}{\line(0,-1){0.45}}
\multiput(60.59,31.93)(0.11,-0.30){3}{\line(0,-1){0.30}}
\multiput(60.23,32.80)(0.12,-0.29){3}{\line(0,-1){0.29}}
\multiput(59.83,33.66)(0.13,-0.29){3}{\line(0,-1){0.29}}
\multiput(59.39,34.49)(0.11,-0.21){4}{\line(0,-1){0.21}}
\multiput(58.90,35.31)(0.12,-0.20){4}{\line(0,-1){0.20}}
\multiput(58.38,36.11)(0.13,-0.20){4}{\line(0,-1){0.20}}
\multiput(57.82,36.88)(0.11,-0.15){5}{\line(0,-1){0.15}}
\multiput(57.22,37.62)(0.12,-0.15){5}{\line(0,-1){0.15}}
\multiput(56.59,38.34)(0.13,-0.14){5}{\line(0,-1){0.14}}
\multiput(55.92,39.03)(0.11,-0.12){6}{\line(0,-1){0.12}}
\multiput(55.22,39.70)(0.12,-0.11){6}{\line(1,0){0.12}}
\multiput(54.48,40.33)(0.15,-0.13){5}{\line(1,0){0.15}}
\multiput(53.72,40.92)(0.15,-0.12){5}{\line(1,0){0.15}}
\multiput(52.93,41.49)(0.16,-0.11){5}{\line(1,0){0.16}}
\multiput(52.11,42.02)(0.20,-0.13){4}{\line(1,0){0.20}}
\multiput(51.26,42.51)(0.21,-0.12){4}{\line(1,0){0.21}}
\multiput(50.39,42.96)(0.22,-0.11){4}{\line(1,0){0.22}}
\multiput(49.51,43.38)(0.30,-0.14){3}{\line(1,0){0.30}}
\multiput(48.60,43.76)(0.30,-0.13){3}{\line(1,0){0.30}}
\multiput(47.67,44.10)(0.31,-0.11){3}{\line(1,0){0.31}}
\multiput(46.73,44.39)(0.47,-0.15){2}{\line(1,0){0.47}}
\multiput(45.77,44.65)(0.48,-0.13){2}{\line(1,0){0.48}}
\multiput(44.81,44.86)(0.48,-0.11){2}{\line(1,0){0.48}}
\multiput(43.83,45.04)(0.98,-0.17){1}{\line(1,0){0.98}}
\multiput(42.85,45.17)(0.98,-0.13){1}{\line(1,0){0.98}}
\multiput(41.86,45.25)(0.99,-0.09){1}{\line(1,0){0.99}}
\multiput(40.87,45.29)(0.99,-0.04){1}{\line(1,0){0.99}}
\put(39.88,45.29){\line(1,0){0.99}}
\multiput(38.88,45.25)(0.99,0.04){1}{\line(1,0){0.99}}
\multiput(37.89,45.17)(0.99,0.09){1}{\line(1,0){0.99}}
\multiput(36.91,45.04)(0.98,0.13){1}{\line(1,0){0.98}}
\multiput(35.94,44.86)(0.98,0.17){1}{\line(1,0){0.98}}
\multiput(34.97,44.65)(0.48,0.11){2}{\line(1,0){0.48}}
\multiput(34.01,44.39)(0.48,0.13){2}{\line(1,0){0.48}}
\multiput(33.07,44.10)(0.47,0.15){2}{\line(1,0){0.47}}
\multiput(32.15,43.76)(0.31,0.11){3}{\line(1,0){0.31}}
\multiput(31.24,43.38)(0.30,0.13){3}{\line(1,0){0.30}}
\multiput(30.35,42.96)(0.30,0.14){3}{\line(1,0){0.30}}
\multiput(29.48,42.51)(0.22,0.11){4}{\line(1,0){0.22}}
\multiput(28.64,42.02)(0.21,0.12){4}{\line(1,0){0.21}}
\multiput(27.82,41.49)(0.20,0.13){4}{\line(1,0){0.20}}
\multiput(27.03,40.92)(0.16,0.11){5}{\line(1,0){0.16}}
\multiput(26.26,40.33)(0.15,0.12){5}{\line(1,0){0.15}}
\multiput(25.53,39.70)(0.15,0.13){5}{\line(1,0){0.15}}
\multiput(24.83,39.03)(0.12,0.11){6}{\line(1,0){0.12}}
\multiput(24.16,38.34)(0.11,0.12){6}{\line(0,1){0.12}}
\multiput(23.52,37.62)(0.13,0.14){5}{\line(0,1){0.14}}
\multiput(22.92,36.88)(0.12,0.15){5}{\line(0,1){0.15}}
\multiput(22.36,36.11)(0.11,0.15){5}{\line(0,1){0.15}}
\multiput(21.84,35.31)(0.13,0.20){4}{\line(0,1){0.20}}
\multiput(21.36,34.49)(0.12,0.20){4}{\line(0,1){0.20}}
\multiput(20.91,33.66)(0.11,0.21){4}{\line(0,1){0.21}}
\multiput(20.51,32.80)(0.13,0.29){3}{\line(0,1){0.29}}
\multiput(20.15,31.93)(0.12,0.29){3}{\line(0,1){0.29}}
\multiput(19.84,31.04)(0.11,0.30){3}{\line(0,1){0.30}}
\multiput(19.57,30.14)(0.14,0.45){2}{\line(0,1){0.45}}
\multiput(19.34,29.23)(0.11,0.46){2}{\line(0,1){0.46}}
\multiput(19.16,28.31)(0.09,0.46){2}{\line(0,1){0.46}}
\multiput(19.02,27.38)(0.14,0.93){1}{\line(0,1){0.93}}
\multiput(18.93,26.45)(0.09,0.93){1}{\line(0,1){0.93}}
\multiput(18.88,25.52)(0.05,0.93){1}{\line(0,1){0.93}}
\put(18.88,24.58){\line(0,1){0.94}}
\multiput(18.88,24.58)(0.05,-0.93){1}{\line(0,-1){0.93}}
\multiput(18.93,23.65)(0.09,-0.93){1}{\line(0,-1){0.93}}
\multiput(19.02,22.72)(0.14,-0.93){1}{\line(0,-1){0.93}}
\multiput(19.16,21.79)(0.09,-0.46){2}{\line(0,-1){0.46}}
\multiput(19.34,20.87)(0.11,-0.46){2}{\line(0,-1){0.46}}
\multiput(19.57,19.96)(0.14,-0.45){2}{\line(0,-1){0.45}}
\multiput(19.84,19.06)(0.11,-0.30){3}{\line(0,-1){0.30}}
\multiput(20.15,18.17)(0.12,-0.29){3}{\line(0,-1){0.29}}
\multiput(20.51,17.30)(0.13,-0.29){3}{\line(0,-1){0.29}}
\multiput(20.91,16.44)(0.11,-0.21){4}{\line(0,-1){0.21}}
\multiput(21.36,15.61)(0.12,-0.20){4}{\line(0,-1){0.20}}
\multiput(21.84,14.79)(0.13,-0.20){4}{\line(0,-1){0.20}}
\multiput(22.36,13.99)(0.11,-0.15){5}{\line(0,-1){0.15}}
\multiput(22.92,13.22)(0.12,-0.15){5}{\line(0,-1){0.15}}
\multiput(23.52,12.48)(0.13,-0.14){5}{\line(0,-1){0.14}}
\multiput(24.16,11.76)(0.11,-0.12){6}{\line(0,-1){0.12}}
\multiput(24.83,11.07)(0.12,-0.11){6}{\line(1,0){0.12}}
\multiput(25.53,10.40)(0.15,-0.13){5}{\line(1,0){0.15}}
\multiput(26.26,9.77)(0.15,-0.12){5}{\line(1,0){0.15}}
\multiput(27.03,9.18)(0.16,-0.11){5}{\line(1,0){0.16}}
\multiput(27.82,8.61)(0.20,-0.13){4}{\line(1,0){0.20}}
\multiput(28.64,8.08)(0.21,-0.12){4}{\line(1,0){0.21}}
\multiput(29.48,7.59)(0.22,-0.11){4}{\line(1,0){0.22}}
\multiput(30.35,7.14)(0.30,-0.14){3}{\line(1,0){0.30}}
\multiput(31.24,6.72)(0.30,-0.13){3}{\line(1,0){0.30}}
\multiput(32.15,6.34)(0.31,-0.11){3}{\line(1,0){0.31}}
\multiput(33.07,6.00)(0.47,-0.15){2}{\line(1,0){0.47}}
\multiput(34.01,5.71)(0.48,-0.13){2}{\line(1,0){0.48}}
\multiput(34.97,5.45)(0.48,-0.11){2}{\line(1,0){0.48}}
\multiput(35.94,5.24)(0.98,-0.17){1}{\line(1,0){0.98}}
\multiput(36.91,5.06)(0.98,-0.13){1}{\line(1,0){0.98}}
\multiput(37.89,4.93)(0.99,-0.09){1}{\line(1,0){0.99}}
\multiput(38.88,4.85)(0.99,-0.04){1}{\line(1,0){0.99}}
\put(39.88,4.81){\line(1,0){0.99}}
\multiput(40.87,4.81)(0.99,0.04){1}{\line(1,0){0.99}}
\multiput(41.86,4.85)(0.99,0.09){1}{\line(1,0){0.99}}
\multiput(42.85,4.93)(0.98,0.13){1}{\line(1,0){0.98}}
\multiput(43.83,5.06)(0.98,0.17){1}{\line(1,0){0.98}}
\multiput(44.81,5.24)(0.48,0.11){2}{\line(1,0){0.48}}
\multiput(45.77,5.45)(0.48,0.13){2}{\line(1,0){0.48}}
\multiput(46.73,5.71)(0.47,0.15){2}{\line(1,0){0.47}}
\multiput(47.67,6.00)(0.31,0.11){3}{\line(1,0){0.31}}
\multiput(48.60,6.34)(0.30,0.13){3}{\line(1,0){0.30}}
\multiput(49.51,6.72)(0.30,0.14){3}{\line(1,0){0.30}}
\multiput(50.39,7.14)(0.22,0.11){4}{\line(1,0){0.22}}
\multiput(51.26,7.59)(0.21,0.12){4}{\line(1,0){0.21}}
\multiput(52.11,8.08)(0.20,0.13){4}{\line(1,0){0.20}}
\multiput(52.93,8.61)(0.16,0.11){5}{\line(1,0){0.16}}
\multiput(53.72,9.18)(0.15,0.12){5}{\line(1,0){0.15}}
\multiput(54.48,9.77)(0.15,0.13){5}{\line(1,0){0.15}}
\multiput(55.22,10.40)(0.12,0.11){6}{\line(1,0){0.12}}
\multiput(55.92,11.07)(0.11,0.12){6}{\line(0,1){0.12}}
\multiput(56.59,11.76)(0.13,0.14){5}{\line(0,1){0.14}}
\multiput(57.22,12.48)(0.12,0.15){5}{\line(0,1){0.15}}
\multiput(57.82,13.22)(0.11,0.15){5}{\line(0,1){0.15}}
\multiput(58.38,13.99)(0.13,0.20){4}{\line(0,1){0.20}}
\multiput(58.90,14.79)(0.12,0.20){4}{\line(0,1){0.20}}
\multiput(59.39,15.61)(0.11,0.21){4}{\line(0,1){0.21}}
\multiput(59.83,16.44)(0.13,0.29){3}{\line(0,1){0.29}}
\multiput(60.23,17.30)(0.12,0.29){3}{\line(0,1){0.29}}
\multiput(60.59,18.17)(0.11,0.30){3}{\line(0,1){0.30}}
\multiput(60.90,19.06)(0.14,0.45){2}{\line(0,1){0.45}}
\multiput(61.18,19.96)(0.11,0.46){2}{\line(0,1){0.46}}
\multiput(61.40,20.87)(0.09,0.46){2}{\line(0,1){0.46}}
\multiput(61.59,21.79)(0.14,0.93){1}{\line(0,1){0.93}}
\multiput(61.72,22.72)(0.09,0.93){1}{\line(0,1){0.93}}
\multiput(61.81,23.65)(0.05,0.93){1}{\line(0,1){0.93}}

\linethickness{0.15mm}
\put(19.17,27.67){\line(1,0){42.43}}

\linethickness{0.15mm}
\put(19.17,22.28){\line(1,0){42.43}}

\linethickness{0.35mm}
\multiput(40.92,22.35)(0.07,5.32){1}{\line(0,1){5.32}}
\put(41.00,27.67){\vector(0,1){0.12}}
\put(40.92,22.35){\vector(-0,-1){0.12}}

\linethickness{0.25mm}
\multiput(59.40,36.20)(0.16,-0.12){92}{\line(1,0){0.16}}
\put(59.40,36.20){\vector(-4,3){0.12}}

\linethickness{0.25mm}
\put(62.60,25.16){\line(1,0){11.28}}
\put(62.60,25.16){\vector(-1,0){0.12}}

\linethickness{0.25mm}
\multiput(59.80,14.60)(0.16,0.12){87}{\line(1,0){0.16}}
\put(59.80,14.60){\vector(-4,-3){0.12}}

\put(81.20,25.10){\makebox(0,0)[cc]{$u_0=0$}}

\linethickness{0.25mm}
\put(12.98,25.01){\line(1,0){5.16}}
\put(18.14,25.01){\vector(1,0){0.12}}

\put(44.90,24.90){\makebox(0,0)[cc]{$2\varepsilon$}}

\put(44.90,24.90){\makebox(0,0)[cc]{}}

\linethickness{0.15mm}
\multiput(13.20,13.00)(0.12,0.13){108}{\line(0,1){0.13}}
\put(26.10,27.30){\vector(1,1){0.12}}

\linethickness{0.15mm}
\multiput(13.00,12.90)(0.18,0.12){76}{\line(1,0){0.18}}
\put(26.80,22.00){\vector(3,2){0.12}}

\put(9.50,11.40){\makebox(0,0)[cc]{$K$}}

\put(4.50,25.10){\makebox(0,0)[cc]{$u_0=M$}}

\end{picture}

%% file: dYcas1.txt
\unitlength 1mm
\begin{picture}(31.04,27.67)(0,0)

\linethickness{0.15mm}
\put(23.03,13.25){\line(0,1){0.98}}
\multiput(22.94,15.22)(0.09,-0.98){1}{\line(0,-1){0.98}}
\multiput(22.77,16.19)(0.17,-0.97){1}{\line(0,-1){0.97}}
\multiput(22.52,17.13)(0.13,-0.47){2}{\line(0,-1){0.47}}
\multiput(22.18,18.05)(0.11,-0.31){3}{\line(0,-1){0.31}}
\multiput(21.76,18.93)(0.14,-0.29){3}{\line(0,-1){0.29}}
\multiput(21.27,19.77)(0.12,-0.21){4}{\line(0,-1){0.21}}
\multiput(20.70,20.55)(0.11,-0.16){5}{\line(0,-1){0.16}}
\multiput(20.07,21.28)(0.13,-0.15){5}{\line(0,-1){0.15}}
\multiput(19.38,21.94)(0.12,-0.11){6}{\line(1,0){0.12}}
\multiput(18.64,22.54)(0.15,-0.12){5}{\line(1,0){0.15}}
\multiput(17.84,23.05)(0.20,-0.13){4}{\line(1,0){0.20}}
\multiput(17.01,23.49)(0.21,-0.11){4}{\line(1,0){0.21}}
\multiput(16.13,23.85)(0.29,-0.12){3}{\line(1,0){0.29}}
\multiput(15.23,24.12)(0.45,-0.13){2}{\line(1,0){0.45}}
\multiput(14.32,24.30)(0.46,-0.09){2}{\line(1,0){0.46}}
\multiput(13.38,24.39)(0.93,-0.09){1}{\line(1,0){0.93}}
\put(12.45,24.39){\line(1,0){0.94}}
\multiput(11.52,24.30)(0.93,0.09){1}{\line(1,0){0.93}}
\multiput(10.60,24.12)(0.46,0.09){2}{\line(1,0){0.46}}
\multiput(9.70,23.85)(0.45,0.13){2}{\line(1,0){0.45}}
\multiput(8.83,23.49)(0.29,0.12){3}{\line(1,0){0.29}}
\multiput(7.99,23.05)(0.21,0.11){4}{\line(1,0){0.21}}
\multiput(7.20,22.54)(0.20,0.13){4}{\line(1,0){0.20}}
\multiput(6.45,21.94)(0.15,0.12){5}{\line(1,0){0.15}}
\multiput(5.76,21.28)(0.12,0.11){6}{\line(1,0){0.12}}
\multiput(5.13,20.55)(0.13,0.15){5}{\line(0,1){0.15}}
\multiput(4.57,19.77)(0.11,0.16){5}{\line(0,1){0.16}}
\multiput(4.07,18.93)(0.12,0.21){4}{\line(0,1){0.21}}
\multiput(3.66,18.05)(0.14,0.29){3}{\line(0,1){0.29}}
\multiput(3.32,17.13)(0.11,0.31){3}{\line(0,1){0.31}}
\multiput(3.06,16.19)(0.13,0.47){2}{\line(0,1){0.47}}
\multiput(2.89,15.22)(0.17,0.97){1}{\line(0,1){0.97}}
\multiput(2.80,14.24)(0.09,0.98){1}{\line(0,1){0.98}}
\put(2.80,13.25){\line(0,1){0.98}}
\multiput(2.80,13.25)(0.09,-0.98){1}{\line(0,-1){0.98}}
\multiput(2.89,12.28)(0.17,-0.97){1}{\line(0,-1){0.97}}
\multiput(3.06,11.31)(0.13,-0.47){2}{\line(0,-1){0.47}}
\multiput(3.32,10.36)(0.11,-0.31){3}{\line(0,-1){0.31}}
\multiput(3.66,9.44)(0.14,-0.29){3}{\line(0,-1){0.29}}
\multiput(4.07,8.56)(0.12,-0.21){4}{\line(0,-1){0.21}}
\multiput(4.57,7.73)(0.11,-0.16){5}{\line(0,-1){0.16}}
\multiput(5.13,6.94)(0.13,-0.15){5}{\line(0,-1){0.15}}
\multiput(5.76,6.21)(0.12,-0.11){6}{\line(1,0){0.12}}
\multiput(6.45,5.55)(0.15,-0.12){5}{\line(1,0){0.15}}
\multiput(7.20,4.96)(0.20,-0.13){4}{\line(1,0){0.20}}
\multiput(7.99,4.44)(0.21,-0.11){4}{\line(1,0){0.21}}
\multiput(8.83,4.00)(0.29,-0.12){3}{\line(1,0){0.29}}
\multiput(9.70,3.65)(0.45,-0.13){2}{\line(1,0){0.45}}
\multiput(10.60,3.38)(0.46,-0.09){2}{\line(1,0){0.46}}
\multiput(11.52,3.20)(0.93,-0.09){1}{\line(1,0){0.93}}
\put(12.45,3.10){\line(1,0){0.94}}
\multiput(13.38,3.10)(0.93,0.09){1}{\line(1,0){0.93}}
\multiput(14.32,3.20)(0.46,0.09){2}{\line(1,0){0.46}}
\multiput(15.23,3.38)(0.45,0.13){2}{\line(1,0){0.45}}
\multiput(16.13,3.65)(0.29,0.12){3}{\line(1,0){0.29}}
\multiput(17.01,4.00)(0.21,0.11){4}{\line(1,0){0.21}}
\multiput(17.84,4.44)(0.20,0.13){4}{\line(1,0){0.20}}
\multiput(18.64,4.96)(0.15,0.12){5}{\line(1,0){0.15}}
\multiput(19.38,5.55)(0.12,0.11){6}{\line(1,0){0.12}}
\multiput(20.07,6.21)(0.13,0.15){5}{\line(0,1){0.15}}
\multiput(20.70,6.94)(0.11,0.16){5}{\line(0,1){0.16}}
\multiput(21.27,7.73)(0.12,0.21){4}{\line(0,1){0.21}}
\multiput(21.76,8.56)(0.14,0.29){3}{\line(0,1){0.29}}
\multiput(22.18,9.44)(0.11,0.31){3}{\line(0,1){0.31}}
\multiput(22.52,10.36)(0.13,0.47){2}{\line(0,1){0.47}}
\multiput(22.77,11.31)(0.17,0.97){1}{\line(0,1){0.97}}
\multiput(22.94,12.28)(0.09,0.98){1}{\line(0,1){0.98}}

\linethickness{0.15mm}
\multiput(0.94,-1.83)(0.12,0.16){91}{\line(0,1){0.16}}
\put(11.90,12.80){\line(1,0){19.14}}

\linethickness{0.15mm}
\multiput(0.14,27.67)(0.12,-0.15){98}{\line(0,-1){0.15}}

\end{picture}

%% file: dYcas2.txt
\unitlength 1mm
\begin{picture}(32.64,28.22)(0,0)

\linethickness{0.15mm}
\put(31.81,13.56){\line(0,1){1.00}}
\multiput(31.72,15.56)(0.09,-1.00){1}{\line(0,-1){1.00}}
\multiput(31.54,16.54)(0.18,-0.98){1}{\line(0,-1){0.98}}
\multiput(31.28,17.51)(0.13,-0.48){2}{\line(0,-1){0.48}}
\multiput(30.94,18.44)(0.12,-0.31){3}{\line(0,-1){0.31}}
\multiput(30.51,19.34)(0.11,-0.22){4}{\line(0,-1){0.22}}
\multiput(30.00,20.19)(0.13,-0.21){4}{\line(0,-1){0.21}}
\multiput(29.42,20.98)(0.12,-0.16){5}{\line(0,-1){0.16}}
\multiput(28.78,21.72)(0.13,-0.15){5}{\line(0,-1){0.15}}
\multiput(28.07,22.40)(0.12,-0.11){6}{\line(1,0){0.12}}
\multiput(27.30,23.00)(0.15,-0.12){5}{\line(1,0){0.15}}
\multiput(26.49,23.53)(0.20,-0.13){4}{\line(1,0){0.20}}
\multiput(25.63,23.97)(0.21,-0.11){4}{\line(1,0){0.21}}
\multiput(24.74,24.34)(0.30,-0.12){3}{\line(1,0){0.30}}
\multiput(23.81,24.61)(0.46,-0.14){2}{\line(1,0){0.46}}
\multiput(22.87,24.79)(0.47,-0.09){2}{\line(1,0){0.47}}
\multiput(21.91,24.89)(0.96,-0.09){1}{\line(1,0){0.96}}
\put(20.96,24.89){\line(1,0){0.96}}
\multiput(20.00,24.79)(0.96,0.09){1}{\line(1,0){0.96}}
\multiput(19.06,24.61)(0.47,0.09){2}{\line(1,0){0.47}}
\multiput(18.13,24.34)(0.46,0.14){2}{\line(1,0){0.46}}
\multiput(17.24,23.97)(0.30,0.12){3}{\line(1,0){0.30}}
\multiput(16.38,23.53)(0.21,0.11){4}{\line(1,0){0.21}}
\multiput(15.57,23.00)(0.20,0.13){4}{\line(1,0){0.20}}
\multiput(14.80,22.40)(0.15,0.12){5}{\line(1,0){0.15}}
\multiput(14.09,21.72)(0.12,0.11){6}{\line(1,0){0.12}}
\multiput(13.45,20.98)(0.13,0.15){5}{\line(0,1){0.15}}
\multiput(12.87,20.19)(0.12,0.16){5}{\line(0,1){0.16}}
\multiput(12.36,19.34)(0.13,0.21){4}{\line(0,1){0.21}}
\multiput(11.93,18.44)(0.11,0.22){4}{\line(0,1){0.22}}
\multiput(11.59,17.51)(0.12,0.31){3}{\line(0,1){0.31}}
\multiput(11.33,16.54)(0.13,0.48){2}{\line(0,1){0.48}}
\multiput(11.15,15.56)(0.18,0.98){1}{\line(0,1){0.98}}
\multiput(11.06,14.56)(0.09,1.00){1}{\line(0,1){1.00}}
\put(11.06,13.56){\line(0,1){1.00}}
\multiput(11.06,13.56)(0.09,-1.00){1}{\line(0,-1){1.00}}
\multiput(11.15,12.57)(0.18,-0.98){1}{\line(0,-1){0.98}}
\multiput(11.33,11.58)(0.13,-0.48){2}{\line(0,-1){0.48}}
\multiput(11.59,10.62)(0.12,-0.31){3}{\line(0,-1){0.31}}
\multiput(11.93,9.69)(0.11,-0.22){4}{\line(0,-1){0.22}}
\multiput(12.36,8.79)(0.13,-0.21){4}{\line(0,-1){0.21}}
\multiput(12.87,7.94)(0.12,-0.16){5}{\line(0,-1){0.16}}
\multiput(13.45,7.14)(0.13,-0.15){5}{\line(0,-1){0.15}}
\multiput(14.09,6.40)(0.12,-0.11){6}{\line(1,0){0.12}}
\multiput(14.80,5.73)(0.15,-0.12){5}{\line(1,0){0.15}}
\multiput(15.57,5.12)(0.20,-0.13){4}{\line(1,0){0.20}}
\multiput(16.38,4.60)(0.21,-0.11){4}{\line(1,0){0.21}}
\multiput(17.24,4.15)(0.30,-0.12){3}{\line(1,0){0.30}}
\multiput(18.13,3.79)(0.46,-0.14){2}{\line(1,0){0.46}}
\multiput(19.06,3.52)(0.47,-0.09){2}{\line(1,0){0.47}}
\multiput(20.00,3.33)(0.96,-0.09){1}{\line(1,0){0.96}}
\put(20.96,3.24){\line(1,0){0.96}}
\multiput(21.91,3.24)(0.96,0.09){1}{\line(1,0){0.96}}
\multiput(22.87,3.33)(0.47,0.09){2}{\line(1,0){0.47}}
\multiput(23.81,3.52)(0.46,0.14){2}{\line(1,0){0.46}}
\multiput(24.74,3.79)(0.30,0.12){3}{\line(1,0){0.30}}
\multiput(25.63,4.15)(0.21,0.11){4}{\line(1,0){0.21}}
\multiput(26.49,4.60)(0.20,0.13){4}{\line(1,0){0.20}}
\multiput(27.30,5.12)(0.15,0.12){5}{\line(1,0){0.15}}
\multiput(28.07,5.73)(0.12,0.11){6}{\line(1,0){0.12}}
\multiput(28.78,6.40)(0.13,0.15){5}{\line(0,1){0.15}}
\multiput(29.42,7.14)(0.12,0.16){5}{\line(0,1){0.16}}
\multiput(30.00,7.94)(0.13,0.21){4}{\line(0,1){0.21}}
\multiput(30.51,8.79)(0.11,0.22){4}{\line(0,1){0.22}}
\multiput(30.94,9.69)(0.12,0.31){3}{\line(0,1){0.31}}
\multiput(31.28,10.62)(0.13,0.48){2}{\line(0,1){0.48}}
\multiput(31.54,11.58)(0.18,0.98){1}{\line(0,1){0.98}}
\multiput(31.72,12.57)(0.09,1.00){1}{\line(0,1){1.00}}

\linethickness{0.15mm}
\multiput(1.75,-1.78)(0.12,0.16){94}{\line(0,1){0.16}}
\put(13.00,13.10){\line(1,0){19.64}}

\linethickness{0.15mm}
\multiput(0.94,28.22)(0.12,-0.15){101}{\line(0,-1){0.15}}

\end{picture}

%% file: dTcas1e.txt
\ifx\JPicScale\undefined\def\JPicScale{1}\fi
\unitlength \JPicScale mm
\begin{picture}(60,52.88)(0,0)
\linethickness{0.3mm}
\multiput(2,6)(0.12,0.17){282}{\line(0,1){0.17}}
\linethickness{0.3mm}
\multiput(35.88,52.88)(0.12,-0.26){201}{\line(0,-1){0.26}}
\linethickness{0.3mm}
\multiput(33,23)(0.12,1.24){24}{\line(0,1){1.24}}
\linethickness{0.3mm}
\multiput(2,6)(0.22,0.12){142}{\line(1,0){0.22}}
\linethickness{0.3mm}
\multiput(33,23)(0.15,-0.12){183}{\line(1,0){0.15}}
\linethickness{0.3mm}
\put(42.59,14.75){\line(0,1){0.51}}
\multiput(42.49,14.25)(0.1,0.5){1}{\line(0,1){0.5}}
\multiput(42.29,13.78)(0.1,0.23){2}{\line(0,1){0.23}}
\multiput(42.01,13.36)(0.14,0.21){2}{\line(0,1){0.21}}
\multiput(41.65,13.01)(0.12,0.12){3}{\line(1,0){0.12}}
\multiput(41.23,12.72)(0.21,0.14){2}{\line(1,0){0.21}}
\multiput(40.75,12.53)(0.24,0.1){2}{\line(1,0){0.24}}
\multiput(40.25,12.43)(0.5,0.1){1}{\line(1,0){0.5}}
\put(39.75,12.43){\line(1,0){0.51}}
\multiput(39.25,12.53)(0.5,-0.1){1}{\line(1,0){0.5}}
\multiput(38.77,12.72)(0.24,-0.1){2}{\line(1,0){0.24}}
\multiput(38.35,13.01)(0.21,-0.14){2}{\line(1,0){0.21}}
\multiput(37.99,13.36)(0.12,-0.12){3}{\line(1,0){0.12}}
\multiput(37.71,13.78)(0.14,-0.21){2}{\line(0,-1){0.21}}
\multiput(37.51,14.25)(0.1,-0.23){2}{\line(0,-1){0.23}}
\multiput(37.41,14.75)(0.1,-0.5){1}{\line(0,-1){0.5}}
\put(37.41,14.75){\line(0,1){0.51}}
\multiput(37.41,15.25)(0.1,0.5){1}{\line(0,1){0.5}}
\multiput(37.51,15.75)(0.1,0.23){2}{\line(0,1){0.23}}
\multiput(37.71,16.22)(0.14,0.21){2}{\line(0,1){0.21}}
\multiput(37.99,16.64)(0.12,0.12){3}{\line(1,0){0.12}}
\multiput(38.35,16.99)(0.21,0.14){2}{\line(1,0){0.21}}
\multiput(38.77,17.28)(0.24,0.1){2}{\line(1,0){0.24}}
\multiput(39.25,17.47)(0.5,0.1){1}{\line(1,0){0.5}}
\put(39.75,17.57){\line(1,0){0.51}}
\multiput(40.25,17.57)(0.5,-0.1){1}{\line(1,0){0.5}}
\multiput(40.75,17.47)(0.24,-0.1){2}{\line(1,0){0.24}}
\multiput(41.23,17.28)(0.21,-0.14){2}{\line(1,0){0.21}}
\multiput(41.65,16.99)(0.12,-0.12){3}{\line(1,0){0.12}}
\multiput(42.01,16.64)(0.14,-0.21){2}{\line(0,-1){0.21}}
\multiput(42.29,16.22)(0.1,-0.23){2}{\line(0,-1){0.23}}
\multiput(42.49,15.75)(0.1,-0.5){1}{\line(0,-1){0.5}}

\linethickness{0.3mm}
\put(46.05,14.75){\line(0,1){0.49}}
\multiput(46.01,14.26)(0.04,0.49){1}{\line(0,1){0.49}}
\multiput(45.93,13.77)(0.08,0.49){1}{\line(0,1){0.49}}
\multiput(45.81,13.29)(0.12,0.48){1}{\line(0,1){0.48}}
\multiput(45.65,12.83)(0.16,0.47){1}{\line(0,1){0.47}}
\multiput(45.44,12.37)(0.1,0.23){2}{\line(0,1){0.23}}
\multiput(45.21,11.94)(0.12,0.22){2}{\line(0,1){0.22}}
\multiput(44.93,11.53)(0.14,0.21){2}{\line(0,1){0.21}}
\multiput(44.62,11.14)(0.1,0.13){3}{\line(0,1){0.13}}
\multiput(44.29,10.77)(0.11,0.12){3}{\line(0,1){0.12}}
\multiput(43.92,10.44)(0.12,0.11){3}{\line(1,0){0.12}}
\multiput(43.52,10.13)(0.13,0.1){3}{\line(1,0){0.13}}
\multiput(43.1,9.86)(0.21,0.14){2}{\line(1,0){0.21}}
\multiput(42.66,9.63)(0.22,0.12){2}{\line(1,0){0.22}}
\multiput(42.2,9.43)(0.23,0.1){2}{\line(1,0){0.23}}
\multiput(41.73,9.27)(0.47,0.16){1}{\line(1,0){0.47}}
\multiput(41.24,9.15)(0.49,0.12){1}{\line(1,0){0.49}}
\multiput(40.75,9.07)(0.49,0.08){1}{\line(1,0){0.49}}
\multiput(40.25,9.03)(0.5,0.04){1}{\line(1,0){0.5}}
\put(39.75,9.03){\line(1,0){0.5}}
\multiput(39.25,9.07)(0.5,-0.04){1}{\line(1,0){0.5}}
\multiput(38.76,9.15)(0.49,-0.08){1}{\line(1,0){0.49}}
\multiput(38.27,9.27)(0.49,-0.12){1}{\line(1,0){0.49}}
\multiput(37.8,9.43)(0.47,-0.16){1}{\line(1,0){0.47}}
\multiput(37.34,9.63)(0.23,-0.1){2}{\line(1,0){0.23}}
\multiput(36.9,9.86)(0.22,-0.12){2}{\line(1,0){0.22}}
\multiput(36.48,10.13)(0.21,-0.14){2}{\line(1,0){0.21}}
\multiput(36.08,10.44)(0.13,-0.1){3}{\line(1,0){0.13}}
\multiput(35.71,10.77)(0.12,-0.11){3}{\line(1,0){0.12}}
\multiput(35.38,11.14)(0.11,-0.12){3}{\line(0,-1){0.12}}
\multiput(35.07,11.53)(0.1,-0.13){3}{\line(0,-1){0.13}}
\multiput(34.79,11.94)(0.14,-0.21){2}{\line(0,-1){0.21}}
\multiput(34.56,12.37)(0.12,-0.22){2}{\line(0,-1){0.22}}
\multiput(34.35,12.83)(0.1,-0.23){2}{\line(0,-1){0.23}}
\multiput(34.19,13.29)(0.16,-0.47){1}{\line(0,-1){0.47}}
\multiput(34.07,13.77)(0.12,-0.48){1}{\line(0,-1){0.48}}
\multiput(33.99,14.26)(0.08,-0.49){1}{\line(0,-1){0.49}}
\multiput(33.95,14.75)(0.04,-0.49){1}{\line(0,-1){0.49}}
\put(33.95,14.75){\line(0,1){0.49}}
\multiput(33.95,15.25)(0.04,0.49){1}{\line(0,1){0.49}}
\multiput(33.99,15.74)(0.08,0.49){1}{\line(0,1){0.49}}
\multiput(34.07,16.23)(0.12,0.48){1}{\line(0,1){0.48}}
\multiput(34.19,16.71)(0.16,0.47){1}{\line(0,1){0.47}}
\multiput(34.35,17.17)(0.1,0.23){2}{\line(0,1){0.23}}
\multiput(34.56,17.63)(0.12,0.22){2}{\line(0,1){0.22}}
\multiput(34.79,18.06)(0.14,0.21){2}{\line(0,1){0.21}}
\multiput(35.07,18.47)(0.1,0.13){3}{\line(0,1){0.13}}
\multiput(35.38,18.86)(0.11,0.12){3}{\line(0,1){0.12}}
\multiput(35.71,19.23)(0.12,0.11){3}{\line(1,0){0.12}}
\multiput(36.08,19.56)(0.13,0.1){3}{\line(1,0){0.13}}
\multiput(36.48,19.87)(0.21,0.14){2}{\line(1,0){0.21}}
\multiput(36.9,20.14)(0.22,0.12){2}{\line(1,0){0.22}}
\multiput(37.34,20.37)(0.23,0.1){2}{\line(1,0){0.23}}
\multiput(37.8,20.57)(0.47,0.16){1}{\line(1,0){0.47}}
\multiput(38.27,20.73)(0.49,0.12){1}{\line(1,0){0.49}}
\multiput(38.76,20.85)(0.49,0.08){1}{\line(1,0){0.49}}
\multiput(39.25,20.93)(0.5,0.04){1}{\line(1,0){0.5}}
\put(39.75,20.97){\line(1,0){0.5}}
\multiput(40.25,20.97)(0.5,-0.04){1}{\line(1,0){0.5}}
\multiput(40.75,20.93)(0.49,-0.08){1}{\line(1,0){0.49}}
\multiput(41.24,20.85)(0.49,-0.12){1}{\line(1,0){0.49}}
\multiput(41.73,20.73)(0.47,-0.16){1}{\line(1,0){0.47}}
\multiput(42.2,20.57)(0.23,-0.1){2}{\line(1,0){0.23}}
\multiput(42.66,20.37)(0.22,-0.12){2}{\line(1,0){0.22}}
\multiput(43.1,20.14)(0.21,-0.14){2}{\line(1,0){0.21}}
\multiput(43.52,19.87)(0.13,-0.1){3}{\line(1,0){0.13}}
\multiput(43.92,19.56)(0.12,-0.11){3}{\line(1,0){0.12}}
\multiput(44.29,19.23)(0.11,-0.12){3}{\line(0,-1){0.12}}
\multiput(44.62,18.86)(0.1,-0.13){3}{\line(0,-1){0.13}}
\multiput(44.93,18.47)(0.14,-0.21){2}{\line(0,-1){0.21}}
\multiput(45.21,18.06)(0.12,-0.22){2}{\line(0,-1){0.22}}
\multiput(45.44,17.63)(0.1,-0.23){2}{\line(0,-1){0.23}}
\multiput(45.65,17.17)(0.16,-0.47){1}{\line(0,-1){0.47}}
\multiput(45.81,16.71)(0.12,-0.48){1}{\line(0,-1){0.48}}
\multiput(45.93,16.23)(0.08,-0.49){1}{\line(0,-1){0.49}}
\multiput(46.01,15.74)(0.04,-0.49){1}{\line(0,-1){0.49}}

\put(25,7){\makebox(0,0)[cc]{$B(r_0)$}}

\linethickness{0.3mm}
\multiput(26,10)(0.39,0.12){33}{\line(1,0){0.39}}
\linethickness{0.3mm}
\multiput(2,6)(1.38,-0.12){42}{\line(1,0){1.38}}
\linethickness{0.3mm}
\multiput(37,15)(0.25,-0.12){8}{\line(1,0){0.25}}
\linethickness{0.3mm}
\multiput(38,12)(0.12,0.25){8}{\line(0,1){0.25}}
\put(41,23){\makebox(0,0)[cc]{$B(Vr_0)$}}

\linethickness{0.3mm}
\multiput(16,43)(0.12,-0.14){125}{\line(0,-1){0.14}}
\linethickness{0.3mm}
\put(31,25){\line(0,1){3}}
\linethickness{0.3mm}
\put(29,25){\line(1,0){2}}
\put(13,48){\makebox(0,0)[cc]{Points of type $\T$ }}

\linethickness{0.3mm}
\multiput(52,10)(0.12,0.34){67}{\line(0,1){0.34}}
\linethickness{0.3mm}
\multiput(51,13)(0.12,-0.38){8}{\line(0,-1){0.38}}
\linethickness{0.3mm}
\put(52,10){\line(1,0){3}}
\put(60,35){\makebox(0,0)[cc]{Points of type $\Y$}}

\end{picture}

%% file: dTcas2.txt
\ifx\JPicScale\undefined\def\JPicScale{1}\fi
\unitlength \JPicScale mm
\begin{picture}(60,52.88)(0,0)
\linethickness{0.3mm}
\multiput(2,6)(0.12,0.17){282}{\line(0,1){0.17}}
\linethickness{0.3mm}
\multiput(35.88,52.88)(0.12,-0.26){201}{\line(0,-1){0.26}}
\linethickness{0.3mm}
\multiput(33,23)(0.12,1.24){24}{\line(0,1){1.24}}
\linethickness{0.3mm}
\multiput(2,6)(0.22,0.12){142}{\line(1,0){0.22}}
\linethickness{0.3mm}
\multiput(33,23)(0.15,-0.12){183}{\line(1,0){0.15}}
\linethickness{0.3mm}
\put(37.59,18.75){\line(0,1){0.51}}
\multiput(37.49,18.25)(0.1,0.5){1}{\line(0,1){0.5}}
\multiput(37.29,17.78)(0.1,0.23){2}{\line(0,1){0.23}}
\multiput(37.01,17.36)(0.14,0.21){2}{\line(0,1){0.21}}
\multiput(36.65,17.01)(0.12,0.12){3}{\line(1,0){0.12}}
\multiput(36.23,16.72)(0.21,0.14){2}{\line(1,0){0.21}}
\multiput(35.75,16.53)(0.24,0.1){2}{\line(1,0){0.24}}
\multiput(35.25,16.43)(0.5,0.1){1}{\line(1,0){0.5}}
\put(34.75,16.43){\line(1,0){0.51}}
\multiput(34.25,16.53)(0.5,-0.1){1}{\line(1,0){0.5}}
\multiput(33.77,16.72)(0.24,-0.1){2}{\line(1,0){0.24}}
\multiput(33.35,17.01)(0.21,-0.14){2}{\line(1,0){0.21}}
\multiput(32.99,17.36)(0.12,-0.12){3}{\line(1,0){0.12}}
\multiput(32.71,17.78)(0.14,-0.21){2}{\line(0,-1){0.21}}
\multiput(32.51,18.25)(0.1,-0.23){2}{\line(0,-1){0.23}}
\multiput(32.41,18.75)(0.1,-0.5){1}{\line(0,-1){0.5}}
\put(32.41,18.75){\line(0,1){0.51}}
\multiput(32.41,19.25)(0.1,0.5){1}{\line(0,1){0.5}}
\multiput(32.51,19.75)(0.1,0.23){2}{\line(0,1){0.23}}
\multiput(32.71,20.22)(0.14,0.21){2}{\line(0,1){0.21}}
\multiput(32.99,20.64)(0.12,0.12){3}{\line(1,0){0.12}}
\multiput(33.35,20.99)(0.21,0.14){2}{\line(1,0){0.21}}
\multiput(33.77,21.28)(0.24,0.1){2}{\line(1,0){0.24}}
\multiput(34.25,21.47)(0.5,0.1){1}{\line(1,0){0.5}}
\put(34.75,21.57){\line(1,0){0.51}}
\multiput(35.25,21.57)(0.5,-0.1){1}{\line(1,0){0.5}}
\multiput(35.75,21.47)(0.24,-0.1){2}{\line(1,0){0.24}}
\multiput(36.23,21.28)(0.21,-0.14){2}{\line(1,0){0.21}}
\multiput(36.65,20.99)(0.12,-0.12){3}{\line(1,0){0.12}}
\multiput(37.01,20.64)(0.14,-0.21){2}{\line(0,-1){0.21}}
\multiput(37.29,20.22)(0.1,-0.23){2}{\line(0,-1){0.23}}
\multiput(37.49,19.75)(0.1,-0.5){1}{\line(0,-1){0.5}}

\linethickness{0.3mm}
\put(41.05,18.75){\line(0,1){0.49}}
\multiput(41.01,18.26)(0.04,0.49){1}{\line(0,1){0.49}}
\multiput(40.93,17.77)(0.08,0.49){1}{\line(0,1){0.49}}
\multiput(40.81,17.29)(0.12,0.48){1}{\line(0,1){0.48}}
\multiput(40.65,16.83)(0.16,0.47){1}{\line(0,1){0.47}}
\multiput(40.44,16.37)(0.1,0.23){2}{\line(0,1){0.23}}
\multiput(40.21,15.94)(0.12,0.22){2}{\line(0,1){0.22}}
\multiput(39.93,15.53)(0.14,0.21){2}{\line(0,1){0.21}}
\multiput(39.62,15.14)(0.1,0.13){3}{\line(0,1){0.13}}
\multiput(39.29,14.77)(0.11,0.12){3}{\line(0,1){0.12}}
\multiput(38.92,14.44)(0.12,0.11){3}{\line(1,0){0.12}}
\multiput(38.52,14.13)(0.13,0.1){3}{\line(1,0){0.13}}
\multiput(38.1,13.86)(0.21,0.14){2}{\line(1,0){0.21}}
\multiput(37.66,13.63)(0.22,0.12){2}{\line(1,0){0.22}}
\multiput(37.2,13.43)(0.23,0.1){2}{\line(1,0){0.23}}
\multiput(36.73,13.27)(0.47,0.16){1}{\line(1,0){0.47}}
\multiput(36.24,13.15)(0.49,0.12){1}{\line(1,0){0.49}}
\multiput(35.75,13.07)(0.49,0.08){1}{\line(1,0){0.49}}
\multiput(35.25,13.03)(0.5,0.04){1}{\line(1,0){0.5}}
\put(34.75,13.03){\line(1,0){0.5}}
\multiput(34.25,13.07)(0.5,-0.04){1}{\line(1,0){0.5}}
\multiput(33.76,13.15)(0.49,-0.08){1}{\line(1,0){0.49}}
\multiput(33.27,13.27)(0.49,-0.12){1}{\line(1,0){0.49}}
\multiput(32.8,13.43)(0.47,-0.16){1}{\line(1,0){0.47}}
\multiput(32.34,13.63)(0.23,-0.1){2}{\line(1,0){0.23}}
\multiput(31.9,13.86)(0.22,-0.12){2}{\line(1,0){0.22}}
\multiput(31.48,14.13)(0.21,-0.14){2}{\line(1,0){0.21}}
\multiput(31.08,14.44)(0.13,-0.1){3}{\line(1,0){0.13}}
\multiput(30.71,14.77)(0.12,-0.11){3}{\line(1,0){0.12}}
\multiput(30.38,15.14)(0.11,-0.12){3}{\line(0,-1){0.12}}
\multiput(30.07,15.53)(0.1,-0.13){3}{\line(0,-1){0.13}}
\multiput(29.79,15.94)(0.14,-0.21){2}{\line(0,-1){0.21}}
\multiput(29.56,16.37)(0.12,-0.22){2}{\line(0,-1){0.22}}
\multiput(29.35,16.83)(0.1,-0.23){2}{\line(0,-1){0.23}}
\multiput(29.19,17.29)(0.16,-0.47){1}{\line(0,-1){0.47}}
\multiput(29.07,17.77)(0.12,-0.48){1}{\line(0,-1){0.48}}
\multiput(28.99,18.26)(0.08,-0.49){1}{\line(0,-1){0.49}}
\multiput(28.95,18.75)(0.04,-0.49){1}{\line(0,-1){0.49}}
\put(28.95,18.75){\line(0,1){0.49}}
\multiput(28.95,19.25)(0.04,0.49){1}{\line(0,1){0.49}}
\multiput(28.99,19.74)(0.08,0.49){1}{\line(0,1){0.49}}
\multiput(29.07,20.23)(0.12,0.48){1}{\line(0,1){0.48}}
\multiput(29.19,20.71)(0.16,0.47){1}{\line(0,1){0.47}}
\multiput(29.35,21.17)(0.1,0.23){2}{\line(0,1){0.23}}
\multiput(29.56,21.63)(0.12,0.22){2}{\line(0,1){0.22}}
\multiput(29.79,22.06)(0.14,0.21){2}{\line(0,1){0.21}}
\multiput(30.07,22.47)(0.1,0.13){3}{\line(0,1){0.13}}
\multiput(30.38,22.86)(0.11,0.12){3}{\line(0,1){0.12}}
\multiput(30.71,23.23)(0.12,0.11){3}{\line(1,0){0.12}}
\multiput(31.08,23.56)(0.13,0.1){3}{\line(1,0){0.13}}
\multiput(31.48,23.87)(0.21,0.14){2}{\line(1,0){0.21}}
\multiput(31.9,24.14)(0.22,0.12){2}{\line(1,0){0.22}}
\multiput(32.34,24.37)(0.23,0.1){2}{\line(1,0){0.23}}
\multiput(32.8,24.57)(0.47,0.16){1}{\line(1,0){0.47}}
\multiput(33.27,24.73)(0.49,0.12){1}{\line(1,0){0.49}}
\multiput(33.76,24.85)(0.49,0.08){1}{\line(1,0){0.49}}
\multiput(34.25,24.93)(0.5,0.04){1}{\line(1,0){0.5}}
\put(34.75,24.97){\line(1,0){0.5}}
\multiput(35.25,24.97)(0.5,-0.04){1}{\line(1,0){0.5}}
\multiput(35.75,24.93)(0.49,-0.08){1}{\line(1,0){0.49}}
\multiput(36.24,24.85)(0.49,-0.12){1}{\line(1,0){0.49}}
\multiput(36.73,24.73)(0.47,-0.16){1}{\line(1,0){0.47}}
\multiput(37.2,24.57)(0.23,-0.1){2}{\line(1,0){0.23}}
\multiput(37.66,24.37)(0.22,-0.12){2}{\line(1,0){0.22}}
\multiput(38.1,24.14)(0.21,-0.14){2}{\line(1,0){0.21}}
\multiput(38.52,23.87)(0.13,-0.1){3}{\line(1,0){0.13}}
\multiput(38.92,23.56)(0.12,-0.11){3}{\line(1,0){0.12}}
\multiput(39.29,23.23)(0.11,-0.12){3}{\line(0,-1){0.12}}
\multiput(39.62,22.86)(0.1,-0.13){3}{\line(0,-1){0.13}}
\multiput(39.93,22.47)(0.14,-0.21){2}{\line(0,-1){0.21}}
\multiput(40.21,22.06)(0.12,-0.22){2}{\line(0,-1){0.22}}
\multiput(40.44,21.63)(0.1,-0.23){2}{\line(0,-1){0.23}}
\multiput(40.65,21.17)(0.16,-0.47){1}{\line(0,-1){0.47}}
\multiput(40.81,20.71)(0.12,-0.48){1}{\line(0,-1){0.48}}
\multiput(40.93,20.23)(0.08,-0.49){1}{\line(0,-1){0.49}}
\multiput(41.01,19.74)(0.04,-0.49){1}{\line(0,-1){0.49}}

\put(25,7){\makebox(0,0)[cc]{$B(r_0)$}}

\linethickness{0.3mm}
\multiput(26,10)(0.12,0.12){67}{\line(1,0){0.12}}
\linethickness{0.3mm}
\put(34,16){\line(0,1){2}}
\linethickness{0.3mm}
\put(32,18){\line(1,0){2}}
\linethickness{0.3mm}
\multiput(2,6)(1.38,-0.12){42}{\line(1,0){1.38}}
\linethickness{0.3mm}
\put(47,19.75){\line(0,1){0.5}}
\multiput(46.98,19.26)(0.02,0.5){1}{\line(0,1){0.5}}
\multiput(46.94,18.76)(0.04,0.49){1}{\line(0,1){0.49}}
\multiput(46.88,18.27)(0.06,0.49){1}{\line(0,1){0.49}}
\multiput(46.8,17.78)(0.08,0.49){1}{\line(0,1){0.49}}
\multiput(46.7,17.3)(0.1,0.49){1}{\line(0,1){0.49}}
\multiput(46.59,16.81)(0.12,0.48){1}{\line(0,1){0.48}}
\multiput(46.45,16.34)(0.14,0.48){1}{\line(0,1){0.48}}
\multiput(46.3,15.87)(0.15,0.47){1}{\line(0,1){0.47}}
\multiput(46.12,15.41)(0.17,0.46){1}{\line(0,1){0.46}}
\multiput(45.93,14.95)(0.1,0.23){2}{\line(0,1){0.23}}
\multiput(45.72,14.51)(0.1,0.22){2}{\line(0,1){0.22}}
\multiput(45.5,14.07)(0.11,0.22){2}{\line(0,1){0.22}}
\multiput(45.26,13.65)(0.12,0.21){2}{\line(0,1){0.21}}
\multiput(45,13.23)(0.13,0.21){2}{\line(0,1){0.21}}
\multiput(44.72,12.83)(0.14,0.2){2}{\line(0,1){0.2}}
\multiput(44.43,12.44)(0.15,0.2){2}{\line(0,1){0.2}}
\multiput(44.12,12.06)(0.1,0.13){3}{\line(0,1){0.13}}
\multiput(43.8,11.69)(0.11,0.12){3}{\line(0,1){0.12}}
\multiput(43.46,11.34)(0.11,0.12){3}{\line(0,1){0.12}}
\multiput(43.11,11.01)(0.12,0.11){3}{\line(1,0){0.12}}
\multiput(42.75,10.68)(0.12,0.11){3}{\line(1,0){0.12}}
\multiput(42.37,10.38)(0.13,0.1){3}{\line(1,0){0.13}}
\multiput(41.99,10.09)(0.19,0.14){2}{\line(1,0){0.19}}
\multiput(41.59,9.82)(0.2,0.14){2}{\line(1,0){0.2}}
\multiput(41.18,9.57)(0.2,0.13){2}{\line(1,0){0.2}}
\multiput(40.76,9.33)(0.21,0.12){2}{\line(1,0){0.21}}
\multiput(40.34,9.11)(0.21,0.11){2}{\line(1,0){0.21}}
\multiput(39.9,8.91)(0.22,0.1){2}{\line(1,0){0.22}}
\multiput(39.46,8.73)(0.22,0.09){2}{\line(1,0){0.22}}
\multiput(39.01,8.57)(0.45,0.16){1}{\line(1,0){0.45}}
\multiput(38.55,8.43)(0.46,0.14){1}{\line(1,0){0.46}}
\multiput(38.09,8.31)(0.46,0.12){1}{\line(1,0){0.46}}
\multiput(37.63,8.21)(0.47,0.1){1}{\line(1,0){0.47}}
\multiput(37.16,8.13)(0.47,0.08){1}{\line(1,0){0.47}}
\multiput(36.69,8.06)(0.47,0.06){1}{\line(1,0){0.47}}
\multiput(36.21,8.02)(0.47,0.04){1}{\line(1,0){0.47}}
\multiput(35.74,8)(0.47,0.02){1}{\line(1,0){0.47}}
\put(35.26,8){\line(1,0){0.48}}
\multiput(34.79,8.02)(0.47,-0.02){1}{\line(1,0){0.47}}
\multiput(34.31,8.06)(0.47,-0.04){1}{\line(1,0){0.47}}
\multiput(33.84,8.13)(0.47,-0.06){1}{\line(1,0){0.47}}
\multiput(33.37,8.21)(0.47,-0.08){1}{\line(1,0){0.47}}
\multiput(32.91,8.31)(0.47,-0.1){1}{\line(1,0){0.47}}
\multiput(32.45,8.43)(0.46,-0.12){1}{\line(1,0){0.46}}
\multiput(31.99,8.57)(0.46,-0.14){1}{\line(1,0){0.46}}
\multiput(31.54,8.73)(0.45,-0.16){1}{\line(1,0){0.45}}
\multiput(31.1,8.91)(0.22,-0.09){2}{\line(1,0){0.22}}
\multiput(30.66,9.11)(0.22,-0.1){2}{\line(1,0){0.22}}
\multiput(30.24,9.33)(0.21,-0.11){2}{\line(1,0){0.21}}
\multiput(29.82,9.57)(0.21,-0.12){2}{\line(1,0){0.21}}
\multiput(29.41,9.82)(0.2,-0.13){2}{\line(1,0){0.2}}
\multiput(29.01,10.09)(0.2,-0.14){2}{\line(1,0){0.2}}
\multiput(28.63,10.38)(0.19,-0.14){2}{\line(1,0){0.19}}
\multiput(28.25,10.68)(0.13,-0.1){3}{\line(1,0){0.13}}
\multiput(27.89,11.01)(0.12,-0.11){3}{\line(1,0){0.12}}
\multiput(27.54,11.34)(0.12,-0.11){3}{\line(1,0){0.12}}
\multiput(27.2,11.69)(0.11,-0.12){3}{\line(0,-1){0.12}}
\multiput(26.88,12.06)(0.11,-0.12){3}{\line(0,-1){0.12}}
\multiput(26.57,12.44)(0.1,-0.13){3}{\line(0,-1){0.13}}
\multiput(26.28,12.83)(0.15,-0.2){2}{\line(0,-1){0.2}}
\multiput(26,13.23)(0.14,-0.2){2}{\line(0,-1){0.2}}
\multiput(25.74,13.65)(0.13,-0.21){2}{\line(0,-1){0.21}}
\multiput(25.5,14.07)(0.12,-0.21){2}{\line(0,-1){0.21}}
\multiput(25.28,14.51)(0.11,-0.22){2}{\line(0,-1){0.22}}
\multiput(25.07,14.95)(0.1,-0.22){2}{\line(0,-1){0.22}}
\multiput(24.88,15.41)(0.1,-0.23){2}{\line(0,-1){0.23}}
\multiput(24.7,15.87)(0.17,-0.46){1}{\line(0,-1){0.46}}
\multiput(24.55,16.34)(0.15,-0.47){1}{\line(0,-1){0.47}}
\multiput(24.41,16.81)(0.14,-0.48){1}{\line(0,-1){0.48}}
\multiput(24.3,17.3)(0.12,-0.48){1}{\line(0,-1){0.48}}
\multiput(24.2,17.78)(0.1,-0.49){1}{\line(0,-1){0.49}}
\multiput(24.12,18.27)(0.08,-0.49){1}{\line(0,-1){0.49}}
\multiput(24.06,18.76)(0.06,-0.49){1}{\line(0,-1){0.49}}
\multiput(24.02,19.26)(0.04,-0.49){1}{\line(0,-1){0.49}}
\multiput(24,19.75)(0.02,-0.5){1}{\line(0,-1){0.5}}
\put(24,19.75){\line(0,1){0.5}}
\multiput(24,20.25)(0.02,0.5){1}{\line(0,1){0.5}}
\multiput(24.02,20.74)(0.04,0.49){1}{\line(0,1){0.49}}
\multiput(24.06,21.24)(0.06,0.49){1}{\line(0,1){0.49}}
\multiput(24.12,21.73)(0.08,0.49){1}{\line(0,1){0.49}}
\multiput(24.2,22.22)(0.1,0.49){1}{\line(0,1){0.49}}
\multiput(24.3,22.7)(0.12,0.48){1}{\line(0,1){0.48}}
\multiput(24.41,23.19)(0.14,0.48){1}{\line(0,1){0.48}}
\multiput(24.55,23.66)(0.15,0.47){1}{\line(0,1){0.47}}
\multiput(24.7,24.13)(0.17,0.46){1}{\line(0,1){0.46}}
\multiput(24.88,24.59)(0.1,0.23){2}{\line(0,1){0.23}}
\multiput(25.07,25.05)(0.1,0.22){2}{\line(0,1){0.22}}
\multiput(25.28,25.49)(0.11,0.22){2}{\line(0,1){0.22}}
\multiput(25.5,25.93)(0.12,0.21){2}{\line(0,1){0.21}}
\multiput(25.74,26.35)(0.13,0.21){2}{\line(0,1){0.21}}
\multiput(26,26.77)(0.14,0.2){2}{\line(0,1){0.2}}
\multiput(26.28,27.17)(0.15,0.2){2}{\line(0,1){0.2}}
\multiput(26.57,27.56)(0.1,0.13){3}{\line(0,1){0.13}}
\multiput(26.88,27.94)(0.11,0.12){3}{\line(0,1){0.12}}
\multiput(27.2,28.31)(0.11,0.12){3}{\line(0,1){0.12}}
\multiput(27.54,28.66)(0.12,0.11){3}{\line(1,0){0.12}}
\multiput(27.89,28.99)(0.12,0.11){3}{\line(1,0){0.12}}
\multiput(28.25,29.32)(0.13,0.1){3}{\line(1,0){0.13}}
\multiput(28.63,29.62)(0.19,0.14){2}{\line(1,0){0.19}}
\multiput(29.01,29.91)(0.2,0.14){2}{\line(1,0){0.2}}
\multiput(29.41,30.18)(0.2,0.13){2}{\line(1,0){0.2}}
\multiput(29.82,30.43)(0.21,0.12){2}{\line(1,0){0.21}}
\multiput(30.24,30.67)(0.21,0.11){2}{\line(1,0){0.21}}
\multiput(30.66,30.89)(0.22,0.1){2}{\line(1,0){0.22}}
\multiput(31.1,31.09)(0.22,0.09){2}{\line(1,0){0.22}}
\multiput(31.54,31.27)(0.45,0.16){1}{\line(1,0){0.45}}
\multiput(31.99,31.43)(0.46,0.14){1}{\line(1,0){0.46}}
\multiput(32.45,31.57)(0.46,0.12){1}{\line(1,0){0.46}}
\multiput(32.91,31.69)(0.47,0.1){1}{\line(1,0){0.47}}
\multiput(33.37,31.79)(0.47,0.08){1}{\line(1,0){0.47}}
\multiput(33.84,31.87)(0.47,0.06){1}{\line(1,0){0.47}}
\multiput(34.31,31.94)(0.47,0.04){1}{\line(1,0){0.47}}
\multiput(34.79,31.98)(0.47,0.02){1}{\line(1,0){0.47}}
\put(35.26,32){\line(1,0){0.48}}
\multiput(35.74,32)(0.47,-0.02){1}{\line(1,0){0.47}}
\multiput(36.21,31.98)(0.47,-0.04){1}{\line(1,0){0.47}}
\multiput(36.69,31.94)(0.47,-0.06){1}{\line(1,0){0.47}}
\multiput(37.16,31.87)(0.47,-0.08){1}{\line(1,0){0.47}}
\multiput(37.63,31.79)(0.47,-0.1){1}{\line(1,0){0.47}}
\multiput(38.09,31.69)(0.46,-0.12){1}{\line(1,0){0.46}}
\multiput(38.55,31.57)(0.46,-0.14){1}{\line(1,0){0.46}}
\multiput(39.01,31.43)(0.45,-0.16){1}{\line(1,0){0.45}}
\multiput(39.46,31.27)(0.22,-0.09){2}{\line(1,0){0.22}}
\multiput(39.9,31.09)(0.22,-0.1){2}{\line(1,0){0.22}}
\multiput(40.34,30.89)(0.21,-0.11){2}{\line(1,0){0.21}}
\multiput(40.76,30.67)(0.21,-0.12){2}{\line(1,0){0.21}}
\multiput(41.18,30.43)(0.2,-0.13){2}{\line(1,0){0.2}}
\multiput(41.59,30.18)(0.2,-0.14){2}{\line(1,0){0.2}}
\multiput(41.99,29.91)(0.19,-0.14){2}{\line(1,0){0.19}}
\multiput(42.37,29.62)(0.13,-0.1){3}{\line(1,0){0.13}}
\multiput(42.75,29.32)(0.12,-0.11){3}{\line(1,0){0.12}}
\multiput(43.11,28.99)(0.12,-0.11){3}{\line(1,0){0.12}}
\multiput(43.46,28.66)(0.11,-0.12){3}{\line(0,-1){0.12}}
\multiput(43.8,28.31)(0.11,-0.12){3}{\line(0,-1){0.12}}
\multiput(44.12,27.94)(0.1,-0.13){3}{\line(0,-1){0.13}}
\multiput(44.43,27.56)(0.15,-0.2){2}{\line(0,-1){0.2}}
\multiput(44.72,27.17)(0.14,-0.2){2}{\line(0,-1){0.2}}
\multiput(45,26.77)(0.13,-0.21){2}{\line(0,-1){0.21}}
\multiput(45.26,26.35)(0.12,-0.21){2}{\line(0,-1){0.21}}
\multiput(45.5,25.93)(0.11,-0.22){2}{\line(0,-1){0.22}}
\multiput(45.72,25.49)(0.1,-0.22){2}{\line(0,-1){0.22}}
\multiput(45.93,25.05)(0.1,-0.23){2}{\line(0,-1){0.23}}
\multiput(46.12,24.59)(0.17,-0.46){1}{\line(0,-1){0.46}}
\multiput(46.3,24.13)(0.15,-0.47){1}{\line(0,-1){0.47}}
\multiput(46.45,23.66)(0.14,-0.48){1}{\line(0,-1){0.48}}
\multiput(46.59,23.19)(0.12,-0.48){1}{\line(0,-1){0.48}}
\multiput(46.7,22.7)(0.1,-0.49){1}{\line(0,-1){0.49}}
\multiput(46.8,22.22)(0.08,-0.49){1}{\line(0,-1){0.49}}
\multiput(46.88,21.73)(0.06,-0.49){1}{\line(0,-1){0.49}}
\multiput(46.94,21.24)(0.04,-0.49){1}{\line(0,-1){0.49}}
\multiput(46.98,20.74)(0.02,-0.5){1}{\line(0,-1){0.5}}

\put(41,27){\makebox(0,0)[cc]{$B(V^2r_0)$}}

\end{picture}

%% file: dx.txt
\unitlength 1mm
\begin{picture}(124.60,54.60)(0,0)

\linethickness{0.15mm}
\put(59.61,26.94){\line(0,1){0.98}}
\multiput(59.57,28.90)(0.04,-0.98){1}{\line(0,-1){0.98}}
\multiput(59.49,29.87)(0.08,-0.98){1}{\line(0,-1){0.98}}
\multiput(59.36,30.84)(0.12,-0.97){1}{\line(0,-1){0.97}}
\multiput(59.20,31.81)(0.16,-0.97){1}{\line(0,-1){0.97}}
\multiput(58.99,32.77)(0.10,-0.48){2}{\line(0,-1){0.48}}
\multiput(58.75,33.71)(0.12,-0.47){2}{\line(0,-1){0.47}}
\multiput(58.46,34.65)(0.14,-0.47){2}{\line(0,-1){0.47}}
\multiput(58.14,35.58)(0.11,-0.31){3}{\line(0,-1){0.31}}
\multiput(57.77,36.49)(0.12,-0.30){3}{\line(0,-1){0.30}}
\multiput(57.37,37.38)(0.13,-0.30){3}{\line(0,-1){0.30}}
\multiput(56.93,38.26)(0.11,-0.22){4}{\line(0,-1){0.22}}
\multiput(56.46,39.12)(0.12,-0.22){4}{\line(0,-1){0.22}}
\multiput(55.95,39.96)(0.13,-0.21){4}{\line(0,-1){0.21}}
\multiput(55.40,40.78)(0.11,-0.16){5}{\line(0,-1){0.16}}
\multiput(54.82,41.58)(0.12,-0.16){5}{\line(0,-1){0.16}}
\multiput(54.20,42.35)(0.12,-0.15){5}{\line(0,-1){0.15}}
\multiput(53.56,43.10)(0.13,-0.15){5}{\line(0,-1){0.15}}
\multiput(52.88,43.82)(0.11,-0.12){6}{\line(0,-1){0.12}}
\multiput(52.17,44.51)(0.12,-0.12){6}{\line(1,0){0.12}}
\multiput(51.44,45.17)(0.12,-0.11){6}{\line(1,0){0.12}}
\multiput(50.68,45.80)(0.15,-0.13){5}{\line(1,0){0.15}}
\multiput(49.89,46.41)(0.16,-0.12){5}{\line(1,0){0.16}}
\multiput(49.07,46.97)(0.16,-0.11){5}{\line(1,0){0.16}}
\multiput(48.24,47.51)(0.21,-0.13){4}{\line(1,0){0.21}}
\multiput(47.38,48.01)(0.21,-0.13){4}{\line(1,0){0.21}}
\multiput(46.50,48.48)(0.22,-0.12){4}{\line(1,0){0.22}}
\multiput(45.60,48.91)(0.22,-0.11){4}{\line(1,0){0.22}}
\multiput(44.69,49.30)(0.31,-0.13){3}{\line(1,0){0.31}}
\multiput(43.75,49.65)(0.31,-0.12){3}{\line(1,0){0.31}}
\multiput(42.81,49.97)(0.32,-0.11){3}{\line(1,0){0.32}}
\multiput(41.85,50.25)(0.48,-0.14){2}{\line(1,0){0.48}}
\multiput(40.88,50.49)(0.48,-0.12){2}{\line(1,0){0.48}}
\multiput(39.90,50.69)(0.49,-0.10){2}{\line(1,0){0.49}}
\multiput(38.92,50.85)(0.99,-0.16){1}{\line(1,0){0.99}}
\multiput(37.92,50.97)(0.99,-0.12){1}{\line(1,0){0.99}}
\multiput(36.93,51.05)(1.00,-0.08){1}{\line(1,0){1.00}}
\multiput(35.93,51.09)(1.00,-0.04){1}{\line(1,0){1.00}}
\put(34.93,51.09){\line(1,0){1.00}}
\multiput(33.93,51.05)(1.00,0.04){1}{\line(1,0){1.00}}
\multiput(32.94,50.97)(1.00,0.08){1}{\line(1,0){1.00}}
\multiput(31.94,50.85)(0.99,0.12){1}{\line(1,0){0.99}}
\multiput(30.96,50.69)(0.99,0.16){1}{\line(1,0){0.99}}
\multiput(29.98,50.49)(0.49,0.10){2}{\line(1,0){0.49}}
\multiput(29.01,50.25)(0.48,0.12){2}{\line(1,0){0.48}}
\multiput(28.05,49.97)(0.48,0.14){2}{\line(1,0){0.48}}
\multiput(27.11,49.65)(0.32,0.11){3}{\line(1,0){0.32}}
\multiput(26.17,49.30)(0.31,0.12){3}{\line(1,0){0.31}}
\multiput(25.26,48.91)(0.31,0.13){3}{\line(1,0){0.31}}
\multiput(24.36,48.48)(0.22,0.11){4}{\line(1,0){0.22}}
\multiput(23.48,48.01)(0.22,0.12){4}{\line(1,0){0.22}}
\multiput(22.62,47.51)(0.21,0.13){4}{\line(1,0){0.21}}
\multiput(21.79,46.97)(0.21,0.13){4}{\line(1,0){0.21}}
\multiput(20.97,46.41)(0.16,0.11){5}{\line(1,0){0.16}}
\multiput(20.18,45.80)(0.16,0.12){5}{\line(1,0){0.16}}
\multiput(19.42,45.17)(0.15,0.13){5}{\line(1,0){0.15}}
\multiput(18.69,44.51)(0.12,0.11){6}{\line(1,0){0.12}}
\multiput(17.98,43.82)(0.12,0.12){6}{\line(1,0){0.12}}
\multiput(17.30,43.10)(0.11,0.12){6}{\line(0,1){0.12}}
\multiput(16.66,42.35)(0.13,0.15){5}{\line(0,1){0.15}}
\multiput(16.04,41.58)(0.12,0.15){5}{\line(0,1){0.15}}
\multiput(15.46,40.78)(0.12,0.16){5}{\line(0,1){0.16}}
\multiput(14.91,39.96)(0.11,0.16){5}{\line(0,1){0.16}}
\multiput(14.40,39.12)(0.13,0.21){4}{\line(0,1){0.21}}
\multiput(13.93,38.26)(0.12,0.22){4}{\line(0,1){0.22}}
\multiput(13.49,37.38)(0.11,0.22){4}{\line(0,1){0.22}}
\multiput(13.09,36.49)(0.13,0.30){3}{\line(0,1){0.30}}
\multiput(12.72,35.58)(0.12,0.30){3}{\line(0,1){0.30}}
\multiput(12.40,34.65)(0.11,0.31){3}{\line(0,1){0.31}}
\multiput(12.11,33.71)(0.14,0.47){2}{\line(0,1){0.47}}
\multiput(11.87,32.77)(0.12,0.47){2}{\line(0,1){0.47}}
\multiput(11.66,31.81)(0.10,0.48){2}{\line(0,1){0.48}}
\multiput(11.50,30.84)(0.16,0.97){1}{\line(0,1){0.97}}
\multiput(11.37,29.87)(0.12,0.97){1}{\line(0,1){0.97}}
\multiput(11.29,28.90)(0.08,0.98){1}{\line(0,1){0.98}}
\multiput(11.25,27.92)(0.04,0.98){1}{\line(0,1){0.98}}
\put(11.25,26.94){\line(0,1){0.98}}
\multiput(11.25,26.94)(0.04,-0.98){1}{\line(0,-1){0.98}}
\multiput(11.29,25.96)(0.08,-0.98){1}{\line(0,-1){0.98}}
\multiput(11.37,24.99)(0.12,-0.97){1}{\line(0,-1){0.97}}
\multiput(11.50,24.02)(0.16,-0.97){1}{\line(0,-1){0.97}}
\multiput(11.66,23.05)(0.10,-0.48){2}{\line(0,-1){0.48}}
\multiput(11.87,22.09)(0.12,-0.47){2}{\line(0,-1){0.47}}
\multiput(12.11,21.15)(0.14,-0.47){2}{\line(0,-1){0.47}}
\multiput(12.40,20.21)(0.11,-0.31){3}{\line(0,-1){0.31}}
\multiput(12.72,19.28)(0.12,-0.30){3}{\line(0,-1){0.30}}
\multiput(13.09,18.37)(0.13,-0.30){3}{\line(0,-1){0.30}}
\multiput(13.49,17.48)(0.11,-0.22){4}{\line(0,-1){0.22}}
\multiput(13.93,16.60)(0.12,-0.22){4}{\line(0,-1){0.22}}
\multiput(14.40,15.74)(0.13,-0.21){4}{\line(0,-1){0.21}}
\multiput(14.91,14.90)(0.11,-0.16){5}{\line(0,-1){0.16}}
\multiput(15.46,14.08)(0.12,-0.16){5}{\line(0,-1){0.16}}
\multiput(16.04,13.28)(0.12,-0.15){5}{\line(0,-1){0.15}}
\multiput(16.66,12.51)(0.13,-0.15){5}{\line(0,-1){0.15}}
\multiput(17.30,11.76)(0.11,-0.12){6}{\line(0,-1){0.12}}
\multiput(17.98,11.04)(0.12,-0.12){6}{\line(1,0){0.12}}
\multiput(18.69,10.35)(0.12,-0.11){6}{\line(1,0){0.12}}
\multiput(19.42,9.69)(0.15,-0.13){5}{\line(1,0){0.15}}
\multiput(20.18,9.06)(0.16,-0.12){5}{\line(1,0){0.16}}
\multiput(20.97,8.45)(0.16,-0.11){5}{\line(1,0){0.16}}
\multiput(21.79,7.89)(0.21,-0.13){4}{\line(1,0){0.21}}
\multiput(22.62,7.35)(0.21,-0.13){4}{\line(1,0){0.21}}
\multiput(23.48,6.85)(0.22,-0.12){4}{\line(1,0){0.22}}
\multiput(24.36,6.38)(0.22,-0.11){4}{\line(1,0){0.22}}
\multiput(25.26,5.95)(0.31,-0.13){3}{\line(1,0){0.31}}
\multiput(26.17,5.56)(0.31,-0.12){3}{\line(1,0){0.31}}
\multiput(27.11,5.21)(0.32,-0.11){3}{\line(1,0){0.32}}
\multiput(28.05,4.89)(0.48,-0.14){2}{\line(1,0){0.48}}
\multiput(29.01,4.61)(0.48,-0.12){2}{\line(1,0){0.48}}
\multiput(29.98,4.37)(0.49,-0.10){2}{\line(1,0){0.49}}
\multiput(30.96,4.17)(0.99,-0.16){1}{\line(1,0){0.99}}
\multiput(31.94,4.01)(0.99,-0.12){1}{\line(1,0){0.99}}
\multiput(32.94,3.89)(1.00,-0.08){1}{\line(1,0){1.00}}
\multiput(33.93,3.81)(1.00,-0.04){1}{\line(1,0){1.00}}
\put(34.93,3.77){\line(1,0){1.00}}
\multiput(35.93,3.77)(1.00,0.04){1}{\line(1,0){1.00}}
\multiput(36.93,3.81)(1.00,0.08){1}{\line(1,0){1.00}}
\multiput(37.92,3.89)(0.99,0.12){1}{\line(1,0){0.99}}
\multiput(38.92,4.01)(0.99,0.16){1}{\line(1,0){0.99}}
\multiput(39.90,4.17)(0.49,0.10){2}{\line(1,0){0.49}}
\multiput(40.88,4.37)(0.48,0.12){2}{\line(1,0){0.48}}
\multiput(41.85,4.61)(0.48,0.14){2}{\line(1,0){0.48}}
\multiput(42.81,4.89)(0.32,0.11){3}{\line(1,0){0.32}}
\multiput(43.75,5.21)(0.31,0.12){3}{\line(1,0){0.31}}
\multiput(44.69,5.56)(0.31,0.13){3}{\line(1,0){0.31}}
\multiput(45.60,5.95)(0.22,0.11){4}{\line(1,0){0.22}}
\multiput(46.50,6.38)(0.22,0.12){4}{\line(1,0){0.22}}
\multiput(47.38,6.85)(0.21,0.13){4}{\line(1,0){0.21}}
\multiput(48.24,7.35)(0.21,0.13){4}{\line(1,0){0.21}}
\multiput(49.07,7.89)(0.16,0.11){5}{\line(1,0){0.16}}
\multiput(49.89,8.45)(0.16,0.12){5}{\line(1,0){0.16}}
\multiput(50.68,9.06)(0.15,0.13){5}{\line(1,0){0.15}}
\multiput(51.44,9.69)(0.12,0.11){6}{\line(1,0){0.12}}
\multiput(52.17,10.35)(0.12,0.12){6}{\line(1,0){0.12}}
\multiput(52.88,11.04)(0.11,0.12){6}{\line(0,1){0.12}}
\multiput(53.56,11.76)(0.13,0.15){5}{\line(0,1){0.15}}
\multiput(54.20,12.51)(0.12,0.15){5}{\line(0,1){0.15}}
\multiput(54.82,13.28)(0.12,0.16){5}{\line(0,1){0.16}}
\multiput(55.40,14.08)(0.11,0.16){5}{\line(0,1){0.16}}
\multiput(55.95,14.90)(0.13,0.21){4}{\line(0,1){0.21}}
\multiput(56.46,15.74)(0.12,0.22){4}{\line(0,1){0.22}}
\multiput(56.93,16.60)(0.11,0.22){4}{\line(0,1){0.22}}
\multiput(57.37,17.48)(0.13,0.30){3}{\line(0,1){0.30}}
\multiput(57.77,18.37)(0.12,0.30){3}{\line(0,1){0.30}}
\multiput(58.14,19.28)(0.11,0.31){3}{\line(0,1){0.31}}
\multiput(58.46,20.21)(0.14,0.47){2}{\line(0,1){0.47}}
\multiput(58.75,21.15)(0.12,0.47){2}{\line(0,1){0.47}}
\multiput(58.99,22.09)(0.10,0.48){2}{\line(0,1){0.48}}
\multiput(59.20,23.05)(0.16,0.97){1}{\line(0,1){0.97}}
\multiput(59.36,24.02)(0.12,0.97){1}{\line(0,1){0.97}}
\multiput(59.49,24.99)(0.08,0.98){1}{\line(0,1){0.98}}
\multiput(59.57,25.96)(0.04,0.98){1}{\line(0,1){0.98}}

\linethickness{0.15mm}
\put(85.59,25.53){\line(0,1){0.99}}
\multiput(85.55,27.52)(0.04,-0.99){1}{\line(0,-1){0.99}}
\multiput(85.47,28.51)(0.08,-0.99){1}{\line(0,-1){0.99}}
\multiput(85.35,29.49)(0.12,-0.99){1}{\line(0,-1){0.99}}
\multiput(85.19,30.47)(0.16,-0.98){1}{\line(0,-1){0.98}}
\multiput(84.99,31.44)(0.10,-0.49){2}{\line(0,-1){0.49}}
\multiput(84.75,32.41)(0.12,-0.48){2}{\line(0,-1){0.48}}
\multiput(84.47,33.36)(0.14,-0.47){2}{\line(0,-1){0.47}}
\multiput(84.15,34.29)(0.11,-0.31){3}{\line(0,-1){0.31}}
\multiput(83.80,35.21)(0.12,-0.31){3}{\line(0,-1){0.31}}
\multiput(83.41,36.12)(0.13,-0.30){3}{\line(0,-1){0.30}}
\multiput(82.98,37.01)(0.11,-0.22){4}{\line(0,-1){0.22}}
\multiput(82.51,37.87)(0.12,-0.22){4}{\line(0,-1){0.22}}
\multiput(82.02,38.72)(0.12,-0.21){4}{\line(0,-1){0.21}}
\multiput(81.48,39.54)(0.13,-0.21){4}{\line(0,-1){0.21}}
\multiput(80.92,40.34)(0.11,-0.16){5}{\line(0,-1){0.16}}
\multiput(80.32,41.11)(0.12,-0.15){5}{\line(0,-1){0.15}}
\multiput(79.69,41.86)(0.13,-0.15){5}{\line(0,-1){0.15}}
\multiput(79.04,42.58)(0.13,-0.14){5}{\line(0,-1){0.14}}
\multiput(78.35,43.26)(0.11,-0.11){6}{\line(0,-1){0.11}}
\multiput(77.64,43.92)(0.14,-0.13){5}{\line(1,0){0.14}}
\multiput(76.90,44.54)(0.15,-0.12){5}{\line(1,0){0.15}}
\multiput(76.13,45.13)(0.15,-0.12){5}{\line(1,0){0.15}}
\multiput(75.34,45.69)(0.16,-0.11){5}{\line(1,0){0.16}}
\multiput(74.53,46.21)(0.20,-0.13){4}{\line(1,0){0.20}}
\multiput(73.70,46.70)(0.21,-0.12){4}{\line(1,0){0.21}}
\multiput(72.86,47.14)(0.21,-0.11){4}{\line(1,0){0.21}}
\multiput(71.99,47.55)(0.29,-0.14){3}{\line(1,0){0.29}}
\multiput(71.11,47.92)(0.29,-0.12){3}{\line(1,0){0.29}}
\multiput(70.21,48.25)(0.30,-0.11){3}{\line(1,0){0.30}}
\multiput(69.30,48.54)(0.45,-0.15){2}{\line(1,0){0.45}}
\multiput(68.38,48.80)(0.46,-0.13){2}{\line(1,0){0.46}}
\multiput(67.45,49.00)(0.46,-0.10){2}{\line(1,0){0.46}}
\multiput(66.52,49.17)(0.94,-0.17){1}{\line(1,0){0.94}}
\multiput(65.57,49.30)(0.94,-0.13){1}{\line(1,0){0.94}}
\multiput(64.63,49.38)(0.95,-0.08){1}{\line(1,0){0.95}}
\multiput(63.68,49.42)(0.95,-0.04){1}{\line(1,0){0.95}}
\put(62.72,49.42){\line(1,0){0.95}}
\multiput(61.77,49.38)(0.95,0.04){1}{\line(1,0){0.95}}
\multiput(60.83,49.30)(0.95,0.08){1}{\line(1,0){0.95}}
\multiput(59.88,49.17)(0.94,0.13){1}{\line(1,0){0.94}}
\multiput(58.95,49.00)(0.94,0.17){1}{\line(1,0){0.94}}
\multiput(58.02,48.80)(0.46,0.10){2}{\line(1,0){0.46}}
\multiput(57.10,48.54)(0.46,0.13){2}{\line(1,0){0.46}}
\multiput(56.19,48.25)(0.45,0.15){2}{\line(1,0){0.45}}
\multiput(55.29,47.92)(0.30,0.11){3}{\line(1,0){0.30}}
\multiput(54.41,47.55)(0.29,0.12){3}{\line(1,0){0.29}}
\multiput(53.54,47.14)(0.29,0.14){3}{\line(1,0){0.29}}
\multiput(52.70,46.70)(0.21,0.11){4}{\line(1,0){0.21}}
\multiput(51.87,46.21)(0.21,0.12){4}{\line(1,0){0.21}}
\multiput(51.06,45.69)(0.20,0.13){4}{\line(1,0){0.20}}
\multiput(50.27,45.13)(0.16,0.11){5}{\line(1,0){0.16}}
\multiput(49.50,44.54)(0.15,0.12){5}{\line(1,0){0.15}}
\multiput(48.76,43.92)(0.15,0.12){5}{\line(1,0){0.15}}
\multiput(48.05,43.26)(0.14,0.13){5}{\line(1,0){0.14}}
\multiput(47.36,42.58)(0.11,0.11){6}{\line(0,1){0.11}}
\multiput(46.71,41.86)(0.13,0.14){5}{\line(0,1){0.14}}
\multiput(46.08,41.11)(0.13,0.15){5}{\line(0,1){0.15}}
\multiput(45.48,40.34)(0.12,0.15){5}{\line(0,1){0.15}}
\multiput(44.92,39.54)(0.11,0.16){5}{\line(0,1){0.16}}
\multiput(44.38,38.72)(0.13,0.21){4}{\line(0,1){0.21}}
\multiput(43.89,37.87)(0.12,0.21){4}{\line(0,1){0.21}}
\multiput(43.42,37.01)(0.12,0.22){4}{\line(0,1){0.22}}
\multiput(42.99,36.12)(0.11,0.22){4}{\line(0,1){0.22}}
\multiput(42.60,35.21)(0.13,0.30){3}{\line(0,1){0.30}}
\multiput(42.25,34.29)(0.12,0.31){3}{\line(0,1){0.31}}
\multiput(41.93,33.36)(0.11,0.31){3}{\line(0,1){0.31}}
\multiput(41.65,32.41)(0.14,0.47){2}{\line(0,1){0.47}}
\multiput(41.41,31.44)(0.12,0.48){2}{\line(0,1){0.48}}
\multiput(41.21,30.47)(0.10,0.49){2}{\line(0,1){0.49}}
\multiput(41.05,29.49)(0.16,0.98){1}{\line(0,1){0.98}}
\multiput(40.93,28.51)(0.12,0.99){1}{\line(0,1){0.99}}
\multiput(40.85,27.52)(0.08,0.99){1}{\line(0,1){0.99}}
\multiput(40.81,26.53)(0.04,0.99){1}{\line(0,1){0.99}}
\put(40.81,25.53){\line(0,1){0.99}}
\multiput(40.81,25.53)(0.04,-0.99){1}{\line(0,-1){0.99}}
\multiput(40.85,24.54)(0.08,-0.99){1}{\line(0,-1){0.99}}
\multiput(40.93,23.55)(0.12,-0.99){1}{\line(0,-1){0.99}}
\multiput(41.05,22.57)(0.16,-0.98){1}{\line(0,-1){0.98}}
\multiput(41.21,21.59)(0.10,-0.49){2}{\line(0,-1){0.49}}
\multiput(41.41,20.62)(0.12,-0.48){2}{\line(0,-1){0.48}}
\multiput(41.65,19.65)(0.14,-0.47){2}{\line(0,-1){0.47}}
\multiput(41.93,18.70)(0.11,-0.31){3}{\line(0,-1){0.31}}
\multiput(42.25,17.77)(0.12,-0.31){3}{\line(0,-1){0.31}}
\multiput(42.60,16.85)(0.13,-0.30){3}{\line(0,-1){0.30}}
\multiput(42.99,15.94)(0.11,-0.22){4}{\line(0,-1){0.22}}
\multiput(43.42,15.05)(0.12,-0.22){4}{\line(0,-1){0.22}}
\multiput(43.89,14.19)(0.12,-0.21){4}{\line(0,-1){0.21}}
\multiput(44.38,13.34)(0.13,-0.21){4}{\line(0,-1){0.21}}
\multiput(44.92,12.52)(0.11,-0.16){5}{\line(0,-1){0.16}}
\multiput(45.48,11.72)(0.12,-0.15){5}{\line(0,-1){0.15}}
\multiput(46.08,10.95)(0.13,-0.15){5}{\line(0,-1){0.15}}
\multiput(46.71,10.20)(0.13,-0.14){5}{\line(0,-1){0.14}}
\multiput(47.36,9.48)(0.11,-0.11){6}{\line(0,-1){0.11}}
\multiput(48.05,8.80)(0.14,-0.13){5}{\line(1,0){0.14}}
\multiput(48.76,8.14)(0.15,-0.12){5}{\line(1,0){0.15}}
\multiput(49.50,7.52)(0.15,-0.12){5}{\line(1,0){0.15}}
\multiput(50.27,6.93)(0.16,-0.11){5}{\line(1,0){0.16}}
\multiput(51.06,6.37)(0.20,-0.13){4}{\line(1,0){0.20}}
\multiput(51.87,5.85)(0.21,-0.12){4}{\line(1,0){0.21}}
\multiput(52.70,5.36)(0.21,-0.11){4}{\line(1,0){0.21}}
\multiput(53.54,4.92)(0.29,-0.14){3}{\line(1,0){0.29}}
\multiput(54.41,4.51)(0.29,-0.12){3}{\line(1,0){0.29}}
\multiput(55.29,4.14)(0.30,-0.11){3}{\line(1,0){0.30}}
\multiput(56.19,3.81)(0.45,-0.15){2}{\line(1,0){0.45}}
\multiput(57.10,3.52)(0.46,-0.13){2}{\line(1,0){0.46}}
\multiput(58.02,3.26)(0.46,-0.10){2}{\line(1,0){0.46}}
\multiput(58.95,3.06)(0.94,-0.17){1}{\line(1,0){0.94}}
\multiput(59.88,2.89)(0.94,-0.13){1}{\line(1,0){0.94}}
\multiput(60.83,2.76)(0.95,-0.08){1}{\line(1,0){0.95}}
\multiput(61.77,2.68)(0.95,-0.04){1}{\line(1,0){0.95}}
\put(62.72,2.64){\line(1,0){0.95}}
\multiput(63.68,2.64)(0.95,0.04){1}{\line(1,0){0.95}}
\multiput(64.63,2.68)(0.95,0.08){1}{\line(1,0){0.95}}
\multiput(65.57,2.76)(0.94,0.13){1}{\line(1,0){0.94}}
\multiput(66.52,2.89)(0.94,0.17){1}{\line(1,0){0.94}}
\multiput(67.45,3.06)(0.46,0.10){2}{\line(1,0){0.46}}
\multiput(68.38,3.26)(0.46,0.13){2}{\line(1,0){0.46}}
\multiput(69.30,3.52)(0.45,0.15){2}{\line(1,0){0.45}}
\multiput(70.21,3.81)(0.30,0.11){3}{\line(1,0){0.30}}
\multiput(71.11,4.14)(0.29,0.12){3}{\line(1,0){0.29}}
\multiput(71.99,4.51)(0.29,0.14){3}{\line(1,0){0.29}}
\multiput(72.86,4.92)(0.21,0.11){4}{\line(1,0){0.21}}
\multiput(73.70,5.36)(0.21,0.12){4}{\line(1,0){0.21}}
\multiput(74.53,5.85)(0.20,0.13){4}{\line(1,0){0.20}}
\multiput(75.34,6.37)(0.16,0.11){5}{\line(1,0){0.16}}
\multiput(76.13,6.93)(0.15,0.12){5}{\line(1,0){0.15}}
\multiput(76.90,7.52)(0.15,0.12){5}{\line(1,0){0.15}}
\multiput(77.64,8.14)(0.14,0.13){5}{\line(1,0){0.14}}
\multiput(78.35,8.80)(0.11,0.11){6}{\line(0,1){0.11}}
\multiput(79.04,9.48)(0.13,0.14){5}{\line(0,1){0.14}}
\multiput(79.69,10.20)(0.13,0.15){5}{\line(0,1){0.15}}
\multiput(80.32,10.95)(0.12,0.15){5}{\line(0,1){0.15}}
\multiput(80.92,11.72)(0.11,0.16){5}{\line(0,1){0.16}}
\multiput(81.48,12.52)(0.13,0.21){4}{\line(0,1){0.21}}
\multiput(82.02,13.34)(0.12,0.21){4}{\line(0,1){0.21}}
\multiput(82.51,14.19)(0.12,0.22){4}{\line(0,1){0.22}}
\multiput(82.98,15.05)(0.11,0.22){4}{\line(0,1){0.22}}
\multiput(83.41,15.94)(0.13,0.30){3}{\line(0,1){0.30}}
\multiput(83.80,16.85)(0.12,0.31){3}{\line(0,1){0.31}}
\multiput(84.15,17.77)(0.11,0.31){3}{\line(0,1){0.31}}
\multiput(84.47,18.70)(0.14,0.47){2}{\line(0,1){0.47}}
\multiput(84.75,19.65)(0.12,0.48){2}{\line(0,1){0.48}}
\multiput(84.99,20.62)(0.10,0.49){2}{\line(0,1){0.49}}
\multiput(85.19,21.59)(0.16,0.98){1}{\line(0,1){0.98}}
\multiput(85.35,22.57)(0.12,0.99){1}{\line(0,1){0.99}}
\multiput(85.47,23.55)(0.08,0.99){1}{\line(0,1){0.99}}
\multiput(85.55,24.54)(0.04,0.99){1}{\line(0,1){0.99}}

\linethickness{0.15mm}
\multiput(4.45,25.56)(116.95,-0.14){1}{\line(1,0){116.95}}

\linethickness{0.15mm}
\put(36.59,45.69){\line(0,1){1.02}}
\multiput(36.18,47.65)(0.14,-0.31){3}{\line(0,-1){0.31}}
\multiput(35.43,48.37)(0.13,-0.12){6}{\line(1,0){0.13}}
\multiput(34.44,48.75)(0.33,-0.13){3}{\line(1,0){0.33}}
\put(33.36,48.75){\line(1,0){1.07}}
\multiput(32.37,48.37)(0.33,0.13){3}{\line(1,0){0.33}}
\multiput(31.62,47.65)(0.13,0.12){6}{\line(1,0){0.13}}
\multiput(31.21,46.71)(0.14,0.31){3}{\line(0,1){0.31}}
\put(31.21,45.69){\line(0,1){1.02}}
\multiput(31.21,45.69)(0.14,-0.31){3}{\line(0,-1){0.31}}
\multiput(31.62,44.75)(0.13,-0.12){6}{\line(1,0){0.13}}
\multiput(32.37,44.03)(0.33,-0.13){3}{\line(1,0){0.33}}
\put(33.36,43.65){\line(1,0){1.07}}
\multiput(34.44,43.65)(0.33,0.13){3}{\line(1,0){0.33}}
\multiput(35.43,44.03)(0.13,0.12){6}{\line(1,0){0.13}}
\multiput(36.18,44.75)(0.14,0.31){3}{\line(0,1){0.31}}

\linethickness{0.15mm}
\put(65.35,44.62){\line(0,1){0.96}}
\multiput(64.96,46.46)(0.13,-0.29){3}{\line(0,-1){0.29}}
\multiput(64.24,47.14)(0.12,-0.11){6}{\line(1,0){0.12}}
\multiput(63.31,47.51)(0.31,-0.12){3}{\line(1,0){0.31}}
\put(62.29,47.51){\line(1,0){1.01}}
\multiput(61.36,47.14)(0.31,0.12){3}{\line(1,0){0.31}}
\multiput(60.64,46.46)(0.12,0.11){6}{\line(1,0){0.12}}
\multiput(60.25,45.58)(0.13,0.29){3}{\line(0,1){0.29}}
\put(60.25,44.62){\line(0,1){0.96}}
\multiput(60.25,44.62)(0.13,-0.29){3}{\line(0,-1){0.29}}
\multiput(60.64,43.74)(0.12,-0.11){6}{\line(1,0){0.12}}
\multiput(61.36,43.06)(0.31,-0.12){3}{\line(1,0){0.31}}
\put(62.29,42.69){\line(1,0){1.01}}
\multiput(63.31,42.69)(0.31,0.12){3}{\line(1,0){0.31}}
\multiput(64.24,43.06)(0.12,0.11){6}{\line(1,0){0.12}}
\multiput(64.96,43.74)(0.13,0.29){3}{\line(0,1){0.29}}

\linethickness{0.35mm}
\put(28.35,41.93){\line(1,0){40.10}}
\put(28.35,41.93){\line(0,1){7.54}}
\put(68.45,41.93){\line(0,1){7.54}}
\put(28.35,49.47){\line(1,0){40.10}}

\linethickness{0.15mm}
\qbezier(7.52,26.30)(10.91,24.03)(13.83,24.74)

\linethickness{0.15mm}
\qbezier(13.83,24.74)(16.74,25.45)(20.03,24.14)

\linethickness{0.15mm}
\qbezier(20.18,24.14)(22.70,27.76)(27.44,25.77)

\linethickness{0.15mm}
\qbezier(27.44,25.77)(32.19,23.78)(44.41,23.80)

\linethickness{0.15mm}
\qbezier(44.20,23.80)(53.36,23.50)(51.27,24.90)

\linethickness{0.15mm}
\qbezier(51.27,24.90)(49.17,26.30)(53.33,26.30)

\linethickness{0.15mm}
\qbezier(53.20,26.30)(61.43,27.01)(72.40,26.84)

\linethickness{0.15mm}
\qbezier(72.40,26.84)(83.35,26.67)(80.22,23.26)

\linethickness{0.15mm}
\qbezier(80.22,23.06)(87.61,25.62)(91.67,24.83)

\linethickness{0.15mm}
\qbezier(91.67,24.83)(95.72,24.05)(96.49,26.30)

\linethickness{0.15mm}
\qbezier(96.69,26.10)(102.94,27.15)(108.09,25.49)

\linethickness{0.15mm}
\qbezier(108.09,25.49)(113.25,23.83)(116.43,24.48)

\linethickness{0.15mm}
\multiput(54.22,35.03)(0.13,-0.12){11}{\line(1,0){0.13}}

\linethickness{0.15mm}
\multiput(53.87,33.74)(0.12,0.12){12}{\line(1,0){0.12}}

\put(51.90,36.84){\makebox(0,0)[cc]{$x$}}

\put(52.10,52.10){\makebox(0,0)[cc]{$D^x_k$}}

\put(124.60,25.80){\makebox(0,0)[cc]{$Z_x$}}

\put(115.80,21.60){\makebox(0,0)[cc]{$K$}}

\put(83.70,54.60){\makebox(0,0)[cc]{$D_k^j$}}

\linethickness{0.15mm}
\multiput(62.90,45.40)(0.36,0.12){52}{\line(1,0){0.36}}

\linethickness{0.15mm}
\multiput(33.40,46.20)(1.00,0.12){48}{\line(1,0){1.00}}

\linethickness{0.15mm}
\multiput(33.10,46.30)(0.22,0.12){10}{\line(1,0){0.22}}

\linethickness{0.15mm}
\multiput(33.10,46.20)(0.21,-0.11){8}{\line(1,0){0.21}}

\linethickness{0.15mm}
\multiput(62.90,45.50)(0.13,0.37){3}{\line(0,1){0.37}}

\linethickness{0.15mm}
\multiput(62.90,45.30)(0.53,-0.13){3}{\line(1,0){0.53}}

\linethickness{0.15mm}
\put(75.89,27.08){\line(0,1){0.95}}
\multiput(75.79,28.97)(0.10,-0.94){1}{\line(0,-1){0.94}}
\multiput(75.59,29.89)(0.10,-0.46){2}{\line(0,-1){0.46}}
\multiput(75.30,30.80)(0.15,-0.45){2}{\line(0,-1){0.45}}
\multiput(74.91,31.68)(0.13,-0.29){3}{\line(0,-1){0.29}}
\multiput(74.43,32.51)(0.12,-0.21){4}{\line(0,-1){0.21}}
\multiput(73.87,33.30)(0.11,-0.16){5}{\line(0,-1){0.16}}
\multiput(73.23,34.03)(0.13,-0.15){5}{\line(0,-1){0.15}}
\multiput(72.52,34.70)(0.12,-0.11){6}{\line(1,0){0.12}}
\multiput(71.74,35.30)(0.16,-0.12){5}{\line(1,0){0.16}}
\multiput(70.90,35.83)(0.21,-0.13){4}{\line(1,0){0.21}}
\multiput(70.00,36.27)(0.22,-0.11){4}{\line(1,0){0.22}}
\multiput(69.07,36.64)(0.31,-0.12){3}{\line(1,0){0.31}}
\multiput(68.10,36.91)(0.48,-0.14){2}{\line(1,0){0.48}}
\multiput(67.11,37.10)(0.50,-0.09){2}{\line(1,0){0.50}}
\multiput(66.11,37.19)(1.01,-0.09){1}{\line(1,0){1.01}}
\put(65.09,37.19){\line(1,0){1.01}}
\multiput(64.09,37.10)(1.01,0.09){1}{\line(1,0){1.01}}
\multiput(63.10,36.91)(0.50,0.09){2}{\line(1,0){0.50}}
\multiput(62.13,36.64)(0.48,0.14){2}{\line(1,0){0.48}}
\multiput(61.20,36.27)(0.31,0.12){3}{\line(1,0){0.31}}
\multiput(60.30,35.83)(0.22,0.11){4}{\line(1,0){0.22}}
\multiput(59.46,35.30)(0.21,0.13){4}{\line(1,0){0.21}}
\multiput(58.68,34.70)(0.16,0.12){5}{\line(1,0){0.16}}
\multiput(57.97,34.03)(0.12,0.11){6}{\line(1,0){0.12}}
\multiput(57.33,33.30)(0.13,0.15){5}{\line(0,1){0.15}}
\multiput(56.77,32.51)(0.11,0.16){5}{\line(0,1){0.16}}
\multiput(56.29,31.68)(0.12,0.21){4}{\line(0,1){0.21}}
\multiput(55.90,30.80)(0.13,0.29){3}{\line(0,1){0.29}}
\multiput(55.61,29.89)(0.15,0.45){2}{\line(0,1){0.45}}
\multiput(55.41,28.97)(0.10,0.46){2}{\line(0,1){0.46}}
\multiput(55.31,28.02)(0.10,0.94){1}{\line(0,1){0.94}}
\put(55.31,27.08){\line(0,1){0.95}}
\multiput(55.31,27.08)(0.10,-0.94){1}{\line(0,-1){0.94}}
\multiput(55.41,26.13)(0.10,-0.46){2}{\line(0,-1){0.46}}
\multiput(55.61,25.21)(0.15,-0.45){2}{\line(0,-1){0.45}}
\multiput(55.90,24.30)(0.13,-0.29){3}{\line(0,-1){0.29}}
\multiput(56.29,23.42)(0.12,-0.21){4}{\line(0,-1){0.21}}
\multiput(56.77,22.59)(0.11,-0.16){5}{\line(0,-1){0.16}}
\multiput(57.33,21.80)(0.13,-0.15){5}{\line(0,-1){0.15}}
\multiput(57.97,21.07)(0.12,-0.11){6}{\line(1,0){0.12}}
\multiput(58.68,20.40)(0.16,-0.12){5}{\line(1,0){0.16}}
\multiput(59.46,19.80)(0.21,-0.13){4}{\line(1,0){0.21}}
\multiput(60.30,19.27)(0.22,-0.11){4}{\line(1,0){0.22}}
\multiput(61.20,18.83)(0.31,-0.12){3}{\line(1,0){0.31}}
\multiput(62.13,18.46)(0.48,-0.14){2}{\line(1,0){0.48}}
\multiput(63.10,18.19)(0.50,-0.09){2}{\line(1,0){0.50}}
\multiput(64.09,18.00)(1.01,-0.09){1}{\line(1,0){1.01}}
\put(65.09,17.91){\line(1,0){1.01}}
\multiput(66.11,17.91)(1.01,0.09){1}{\line(1,0){1.01}}
\multiput(67.11,18.00)(0.50,0.09){2}{\line(1,0){0.50}}
\multiput(68.10,18.19)(0.48,0.14){2}{\line(1,0){0.48}}
\multiput(69.07,18.46)(0.31,0.12){3}{\line(1,0){0.31}}
\multiput(70.00,18.83)(0.22,0.11){4}{\line(1,0){0.22}}
\multiput(70.90,19.27)(0.21,0.13){4}{\line(1,0){0.21}}
\multiput(71.74,19.80)(0.16,0.12){5}{\line(1,0){0.16}}
\multiput(72.52,20.40)(0.12,0.11){6}{\line(1,0){0.12}}
\multiput(73.23,21.07)(0.13,0.15){5}{\line(0,1){0.15}}
\multiput(73.87,21.80)(0.11,0.16){5}{\line(0,1){0.16}}
\multiput(74.43,22.59)(0.12,0.21){4}{\line(0,1){0.21}}
\multiput(74.91,23.42)(0.13,0.29){3}{\line(0,1){0.29}}
\multiput(75.30,24.30)(0.15,0.45){2}{\line(0,1){0.45}}
\multiput(75.59,25.21)(0.10,0.46){2}{\line(0,1){0.46}}
\multiput(75.79,26.13)(0.10,0.94){1}{\line(0,1){0.94}}

\linethickness{0.15mm}
\put(44.29,26.62){\line(0,1){0.95}}
\multiput(44.19,28.52)(0.10,-0.95){1}{\line(0,-1){0.95}}
\multiput(43.98,29.46)(0.10,-0.47){2}{\line(0,-1){0.47}}
\multiput(43.68,30.36)(0.10,-0.30){3}{\line(0,-1){0.30}}
\multiput(43.29,31.23)(0.13,-0.29){3}{\line(0,-1){0.29}}
\multiput(42.80,32.06)(0.12,-0.21){4}{\line(0,-1){0.21}}
\multiput(42.23,32.83)(0.11,-0.15){5}{\line(0,-1){0.15}}
\multiput(41.58,33.53)(0.13,-0.14){5}{\line(0,-1){0.14}}
\multiput(40.85,34.17)(0.14,-0.13){5}{\line(1,0){0.14}}
\multiput(40.07,34.73)(0.16,-0.11){5}{\line(1,0){0.16}}
\multiput(39.22,35.21)(0.21,-0.12){4}{\line(1,0){0.21}}
\multiput(38.33,35.60)(0.30,-0.13){3}{\line(1,0){0.30}}
\multiput(37.41,35.89)(0.46,-0.15){2}{\line(1,0){0.46}}
\multiput(36.45,36.09)(0.48,-0.10){2}{\line(1,0){0.48}}
\multiput(35.49,36.19)(0.97,-0.10){1}{\line(1,0){0.97}}
\put(34.51,36.19){\line(1,0){0.97}}
\multiput(33.55,36.09)(0.97,0.10){1}{\line(1,0){0.97}}
\multiput(32.59,35.89)(0.48,0.10){2}{\line(1,0){0.48}}
\multiput(31.67,35.60)(0.46,0.15){2}{\line(1,0){0.46}}
\multiput(30.78,35.21)(0.30,0.13){3}{\line(1,0){0.30}}
\multiput(29.93,34.73)(0.21,0.12){4}{\line(1,0){0.21}}
\multiput(29.15,34.17)(0.16,0.11){5}{\line(1,0){0.16}}
\multiput(28.42,33.53)(0.14,0.13){5}{\line(1,0){0.14}}
\multiput(27.77,32.83)(0.13,0.14){5}{\line(0,1){0.14}}
\multiput(27.20,32.06)(0.11,0.15){5}{\line(0,1){0.15}}
\multiput(26.71,31.23)(0.12,0.21){4}{\line(0,1){0.21}}
\multiput(26.32,30.36)(0.13,0.29){3}{\line(0,1){0.29}}
\multiput(26.02,29.46)(0.10,0.30){3}{\line(0,1){0.30}}
\multiput(25.81,28.52)(0.10,0.47){2}{\line(0,1){0.47}}
\multiput(25.71,27.58)(0.10,0.95){1}{\line(0,1){0.95}}
\put(25.71,26.62){\line(0,1){0.95}}
\multiput(25.71,26.62)(0.10,-0.95){1}{\line(0,-1){0.95}}
\multiput(25.81,25.68)(0.10,-0.47){2}{\line(0,-1){0.47}}
\multiput(26.02,24.74)(0.10,-0.30){3}{\line(0,-1){0.30}}
\multiput(26.32,23.84)(0.13,-0.29){3}{\line(0,-1){0.29}}
\multiput(26.71,22.97)(0.12,-0.21){4}{\line(0,-1){0.21}}
\multiput(27.20,22.14)(0.11,-0.15){5}{\line(0,-1){0.15}}
\multiput(27.77,21.37)(0.13,-0.14){5}{\line(0,-1){0.14}}
\multiput(28.42,20.67)(0.14,-0.13){5}{\line(1,0){0.14}}
\multiput(29.15,20.03)(0.16,-0.11){5}{\line(1,0){0.16}}
\multiput(29.93,19.47)(0.21,-0.12){4}{\line(1,0){0.21}}
\multiput(30.78,18.99)(0.30,-0.13){3}{\line(1,0){0.30}}
\multiput(31.67,18.60)(0.46,-0.15){2}{\line(1,0){0.46}}
\multiput(32.59,18.31)(0.48,-0.10){2}{\line(1,0){0.48}}
\multiput(33.55,18.11)(0.97,-0.10){1}{\line(1,0){0.97}}
\put(34.51,18.01){\line(1,0){0.97}}
\multiput(35.49,18.01)(0.97,0.10){1}{\line(1,0){0.97}}
\multiput(36.45,18.11)(0.48,0.10){2}{\line(1,0){0.48}}
\multiput(37.41,18.31)(0.46,0.15){2}{\line(1,0){0.46}}
\multiput(38.33,18.60)(0.30,0.13){3}{\line(1,0){0.30}}
\multiput(39.22,18.99)(0.21,0.12){4}{\line(1,0){0.21}}
\multiput(40.07,19.47)(0.16,0.11){5}{\line(1,0){0.16}}
\multiput(40.85,20.03)(0.14,0.13){5}{\line(1,0){0.14}}
\multiput(41.58,20.67)(0.13,0.14){5}{\line(0,1){0.14}}
\multiput(42.23,21.37)(0.11,0.15){5}{\line(0,1){0.15}}
\multiput(42.80,22.14)(0.12,0.21){4}{\line(0,1){0.21}}
\multiput(43.29,22.97)(0.13,0.29){3}{\line(0,1){0.29}}
\multiput(43.68,23.84)(0.10,0.30){3}{\line(0,1){0.30}}
\multiput(43.98,24.74)(0.10,0.47){2}{\line(0,1){0.47}}
\multiput(44.19,25.68)(0.10,0.95){1}{\line(0,1){0.95}}

\put(81.30,46.10){\makebox(0,0)[cc]{$10W_j$}}

\put(64.90,31.30){\makebox(0,0)[cc]{$6W_j$}}

\end{picture}

%% file: dxy.tex
\unitlength 1mm
\begin{picture}(84.70,71.40)(0,0)

\linethickness{0.25mm}
\multiput(12.60,64.50)(0.17,-0.12){201}{\line(1,0){0.17}}

\linethickness{0.25mm}
\put(46.20,12.20){\line(0,1){28.10}}

\linethickness{0.25mm}
\multiput(46.20,40.40)(0.13,0.12){218}{\line(1,0){0.13}}

\linethickness{0.15mm}
\qbezier(10.80,64.10)(17.18,60.31)(18.95,60.90)

\linethickness{0.15mm}
\qbezier(18.95,60.90)(20.72,61.49)(21.40,60.70)

\linethickness{0.15mm}
\qbezier(21.30,61.10)(25.12,50.62)(33.41,48.89)

\linethickness{0.15mm}
\qbezier(33.41,48.89)(41.71,47.16)(42.00,41.40)

\linethickness{0.15mm}
\qbezier(41.70,41.60)(44.17,41.21)(44.86,40.24)

\linethickness{0.15mm}
\qbezier(44.86,40.24)(45.55,39.27)(47.20,39.70)

\linethickness{0.15mm}
\qbezier(47.40,39.30)(41.93,27.72)(44.76,28.85)

\linethickness{0.15mm}
\qbezier(44.76,28.85)(47.60,29.98)(47.00,27.10)

\linethickness{0.15mm}
\qbezier(47.40,39.50)(60.20,57.84)(65.21,57.58)

\linethickness{0.15mm}
\qbezier(65.21,57.58)(70.23,57.31)(75.00,63.20)

\linethickness{0.15mm}
\qbezier(47.20,27.10)(48.78,15.59)(45.66,17.58)

\linethickness{0.15mm}
\qbezier(45.66,17.58)(42.54,19.56)(47.20,13.00)

\put(77.50,68.70){\makebox(0,0)[cc]{$Z_x$}}

\linethickness{0.20mm}
\multiput(12.50,60.90)(0.16,-0.12){174}{\line(1,0){0.16}}
\multiput(41.20,40.00)(0.10,-26.20){1}{\line(0,-1){26.20}}
\multiput(35.40,9.20)(0.16,0.12){38}{\line(1,0){0.16}}
\multiput(6.20,56.20)(0.12,-0.19){243}{\line(0,-1){0.19}}
\multiput(6.20,56.20)(0.16,0.12){39}{\line(1,0){0.16}}
\put(12.20,60.90){\line(1,0){0.20}}

\linethickness{0.20mm}
\multiput(50.30,37.40)(0.14,0.12){193}{\line(1,0){0.14}}
\multiput(49.90,14.20)(0.13,7.73){3}{\line(0,1){7.73}}
\multiput(49.90,14.20)(0.22,-0.12){29}{\line(1,0){0.22}}
\multiput(56.20,10.70)(0.12,0.19){238}{\line(0,1){0.19}}
\multiput(78.30,60.40)(0.21,-0.12){31}{\line(1,0){0.21}}

\linethickness{0.20mm}
\multiput(46.70,46.10)(0.14,0.12){172}{\line(1,0){0.14}}
\multiput(17.20,66.30)(0.18,-0.12){168}{\line(1,0){0.18}}
\put(17.20,66.30){\line(0,1){5.10}}
\put(17.20,71.40){\line(1,0){52.90}}
\put(70.10,67.00){\line(0,1){4.40}}

\put(38.40,67.70){\makebox(0,0)[cc]{$D_1^x$}}

\put(11.20,54.50){\makebox(0,0)[cc]{$D_2^x$}}

\put(58.10,21.20){\makebox(0,0)[cc]{$D_3^x$}}

\linethickness{0.20mm}
\put(65.94,44.05){\line(0,1){0.89}}
\multiput(65.90,45.84)(0.04,-0.89){1}{\line(0,-1){0.89}}
\multiput(65.82,46.73)(0.09,-0.89){1}{\line(0,-1){0.89}}
\multiput(65.69,47.62)(0.13,-0.89){1}{\line(0,-1){0.89}}
\multiput(65.51,48.50)(0.17,-0.88){1}{\line(0,-1){0.88}}
\multiput(65.30,49.37)(0.11,-0.44){2}{\line(0,-1){0.44}}
\multiput(65.05,50.24)(0.13,-0.43){2}{\line(0,-1){0.43}}
\multiput(64.75,51.09)(0.15,-0.43){2}{\line(0,-1){0.43}}
\multiput(64.41,51.93)(0.11,-0.28){3}{\line(0,-1){0.28}}
\multiput(64.03,52.76)(0.13,-0.28){3}{\line(0,-1){0.28}}
\multiput(63.61,53.57)(0.14,-0.27){3}{\line(0,-1){0.27}}
\multiput(63.16,54.36)(0.11,-0.20){4}{\line(0,-1){0.20}}
\multiput(62.66,55.13)(0.12,-0.19){4}{\line(0,-1){0.19}}
\multiput(62.13,55.89)(0.13,-0.19){4}{\line(0,-1){0.19}}
\multiput(61.57,56.62)(0.11,-0.15){5}{\line(0,-1){0.15}}
\multiput(60.96,57.33)(0.12,-0.14){5}{\line(0,-1){0.14}}
\multiput(60.33,58.02)(0.13,-0.14){5}{\line(0,-1){0.14}}
\multiput(59.66,58.68)(0.13,-0.13){5}{\line(1,0){0.13}}
\multiput(58.96,59.31)(0.14,-0.13){5}{\line(1,0){0.14}}
\multiput(58.23,59.91)(0.15,-0.12){5}{\line(1,0){0.15}}
\multiput(57.48,60.49)(0.15,-0.11){5}{\line(1,0){0.15}}
\multiput(56.69,61.03)(0.16,-0.11){5}{\line(1,0){0.16}}
\multiput(55.88,61.55)(0.20,-0.13){4}{\line(1,0){0.20}}
\multiput(55.05,62.03)(0.21,-0.12){4}{\line(1,0){0.21}}
\multiput(54.19,62.47)(0.21,-0.11){4}{\line(1,0){0.21}}
\multiput(53.32,62.89)(0.29,-0.14){3}{\line(1,0){0.29}}
\multiput(52.42,63.26)(0.30,-0.13){3}{\line(1,0){0.30}}
\multiput(51.51,63.61)(0.30,-0.11){3}{\line(1,0){0.30}}
\multiput(50.58,63.91)(0.31,-0.10){3}{\line(1,0){0.31}}
\multiput(49.64,64.18)(0.47,-0.13){2}{\line(1,0){0.47}}
\multiput(48.68,64.41)(0.48,-0.12){2}{\line(1,0){0.48}}
\multiput(47.72,64.61)(0.48,-0.10){2}{\line(1,0){0.48}}
\multiput(46.75,64.76)(0.97,-0.16){1}{\line(1,0){0.97}}
\multiput(45.77,64.88)(0.98,-0.12){1}{\line(1,0){0.98}}
\multiput(44.78,64.96)(0.98,-0.08){1}{\line(1,0){0.98}}
\multiput(43.79,65.00)(0.99,-0.04){1}{\line(1,0){0.99}}
\put(42.81,65.00){\line(1,0){0.99}}
\multiput(41.82,64.96)(0.99,0.04){1}{\line(1,0){0.99}}
\multiput(40.83,64.88)(0.98,0.08){1}{\line(1,0){0.98}}
\multiput(39.85,64.76)(0.98,0.12){1}{\line(1,0){0.98}}
\multiput(38.88,64.61)(0.97,0.16){1}{\line(1,0){0.97}}
\multiput(37.92,64.41)(0.48,0.10){2}{\line(1,0){0.48}}
\multiput(36.96,64.18)(0.48,0.12){2}{\line(1,0){0.48}}
\multiput(36.02,63.91)(0.47,0.13){2}{\line(1,0){0.47}}
\multiput(35.09,63.61)(0.31,0.10){3}{\line(1,0){0.31}}
\multiput(34.18,63.26)(0.30,0.11){3}{\line(1,0){0.30}}
\multiput(33.28,62.89)(0.30,0.13){3}{\line(1,0){0.30}}
\multiput(32.41,62.47)(0.29,0.14){3}{\line(1,0){0.29}}
\multiput(31.55,62.03)(0.21,0.11){4}{\line(1,0){0.21}}
\multiput(30.72,61.55)(0.21,0.12){4}{\line(1,0){0.21}}
\multiput(29.91,61.03)(0.20,0.13){4}{\line(1,0){0.20}}
\multiput(29.12,60.49)(0.16,0.11){5}{\line(1,0){0.16}}
\multiput(28.37,59.91)(0.15,0.11){5}{\line(1,0){0.15}}
\multiput(27.64,59.31)(0.15,0.12){5}{\line(1,0){0.15}}
\multiput(26.94,58.68)(0.14,0.13){5}{\line(1,0){0.14}}
\multiput(26.27,58.02)(0.13,0.13){5}{\line(1,0){0.13}}
\multiput(25.64,57.33)(0.13,0.14){5}{\line(0,1){0.14}}
\multiput(25.03,56.62)(0.12,0.14){5}{\line(0,1){0.14}}
\multiput(24.47,55.89)(0.11,0.15){5}{\line(0,1){0.15}}
\multiput(23.94,55.13)(0.13,0.19){4}{\line(0,1){0.19}}
\multiput(23.44,54.36)(0.12,0.19){4}{\line(0,1){0.19}}
\multiput(22.99,53.57)(0.11,0.20){4}{\line(0,1){0.20}}
\multiput(22.57,52.76)(0.14,0.27){3}{\line(0,1){0.27}}
\multiput(22.19,51.93)(0.13,0.28){3}{\line(0,1){0.28}}
\multiput(21.85,51.09)(0.11,0.28){3}{\line(0,1){0.28}}
\multiput(21.55,50.24)(0.15,0.43){2}{\line(0,1){0.43}}
\multiput(21.30,49.37)(0.13,0.43){2}{\line(0,1){0.43}}
\multiput(21.09,48.50)(0.11,0.44){2}{\line(0,1){0.44}}
\multiput(20.91,47.62)(0.17,0.88){1}{\line(0,1){0.88}}
\multiput(20.78,46.73)(0.13,0.89){1}{\line(0,1){0.89}}
\multiput(20.70,45.84)(0.09,0.89){1}{\line(0,1){0.89}}
\multiput(20.66,44.95)(0.04,0.89){1}{\line(0,1){0.89}}
\put(20.66,44.05){\line(0,1){0.89}}
\multiput(20.66,44.05)(0.04,-0.89){1}{\line(0,-1){0.89}}
\multiput(20.70,43.16)(0.09,-0.89){1}{\line(0,-1){0.89}}
\multiput(20.78,42.27)(0.13,-0.89){1}{\line(0,-1){0.89}}
\multiput(20.91,41.38)(0.17,-0.88){1}{\line(0,-1){0.88}}
\multiput(21.09,40.50)(0.11,-0.44){2}{\line(0,-1){0.44}}
\multiput(21.30,39.63)(0.13,-0.43){2}{\line(0,-1){0.43}}
\multiput(21.55,38.76)(0.15,-0.43){2}{\line(0,-1){0.43}}
\multiput(21.85,37.91)(0.11,-0.28){3}{\line(0,-1){0.28}}
\multiput(22.19,37.07)(0.13,-0.28){3}{\line(0,-1){0.28}}
\multiput(22.57,36.24)(0.14,-0.27){3}{\line(0,-1){0.27}}
\multiput(22.99,35.43)(0.11,-0.20){4}{\line(0,-1){0.20}}
\multiput(23.44,34.64)(0.12,-0.19){4}{\line(0,-1){0.19}}
\multiput(23.94,33.87)(0.13,-0.19){4}{\line(0,-1){0.19}}
\multiput(24.47,33.11)(0.11,-0.15){5}{\line(0,-1){0.15}}
\multiput(25.03,32.38)(0.12,-0.14){5}{\line(0,-1){0.14}}
\multiput(25.64,31.67)(0.13,-0.14){5}{\line(0,-1){0.14}}
\multiput(26.27,30.98)(0.13,-0.13){5}{\line(1,0){0.13}}
\multiput(26.94,30.32)(0.14,-0.13){5}{\line(1,0){0.14}}
\multiput(27.64,29.69)(0.15,-0.12){5}{\line(1,0){0.15}}
\multiput(28.37,29.09)(0.15,-0.11){5}{\line(1,0){0.15}}
\multiput(29.12,28.51)(0.16,-0.11){5}{\line(1,0){0.16}}
\multiput(29.91,27.97)(0.20,-0.13){4}{\line(1,0){0.20}}
\multiput(30.72,27.45)(0.21,-0.12){4}{\line(1,0){0.21}}
\multiput(31.55,26.97)(0.21,-0.11){4}{\line(1,0){0.21}}
\multiput(32.41,26.53)(0.29,-0.14){3}{\line(1,0){0.29}}
\multiput(33.28,26.11)(0.30,-0.13){3}{\line(1,0){0.30}}
\multiput(34.18,25.74)(0.30,-0.11){3}{\line(1,0){0.30}}
\multiput(35.09,25.39)(0.31,-0.10){3}{\line(1,0){0.31}}
\multiput(36.02,25.09)(0.47,-0.13){2}{\line(1,0){0.47}}
\multiput(36.96,24.82)(0.48,-0.12){2}{\line(1,0){0.48}}
\multiput(37.92,24.59)(0.48,-0.10){2}{\line(1,0){0.48}}
\multiput(38.88,24.39)(0.97,-0.16){1}{\line(1,0){0.97}}
\multiput(39.85,24.24)(0.98,-0.12){1}{\line(1,0){0.98}}
\multiput(40.83,24.12)(0.98,-0.08){1}{\line(1,0){0.98}}
\multiput(41.82,24.04)(0.99,-0.04){1}{\line(1,0){0.99}}
\put(42.81,24.00){\line(1,0){0.99}}
\multiput(43.79,24.00)(0.99,0.04){1}{\line(1,0){0.99}}
\multiput(44.78,24.04)(0.98,0.08){1}{\line(1,0){0.98}}
\multiput(45.77,24.12)(0.98,0.12){1}{\line(1,0){0.98}}
\multiput(46.75,24.24)(0.97,0.16){1}{\line(1,0){0.97}}
\multiput(47.72,24.39)(0.48,0.10){2}{\line(1,0){0.48}}
\multiput(48.68,24.59)(0.48,0.12){2}{\line(1,0){0.48}}
\multiput(49.64,24.82)(0.47,0.13){2}{\line(1,0){0.47}}
\multiput(50.58,25.09)(0.31,0.10){3}{\line(1,0){0.31}}
\multiput(51.51,25.39)(0.30,0.11){3}{\line(1,0){0.30}}
\multiput(52.42,25.74)(0.30,0.13){3}{\line(1,0){0.30}}
\multiput(53.32,26.11)(0.29,0.14){3}{\line(1,0){0.29}}
\multiput(54.19,26.53)(0.21,0.11){4}{\line(1,0){0.21}}
\multiput(55.05,26.97)(0.21,0.12){4}{\line(1,0){0.21}}
\multiput(55.88,27.45)(0.20,0.13){4}{\line(1,0){0.20}}
\multiput(56.69,27.97)(0.16,0.11){5}{\line(1,0){0.16}}
\multiput(57.48,28.51)(0.15,0.11){5}{\line(1,0){0.15}}
\multiput(58.23,29.09)(0.15,0.12){5}{\line(1,0){0.15}}
\multiput(58.96,29.69)(0.14,0.13){5}{\line(1,0){0.14}}
\multiput(59.66,30.32)(0.13,0.13){5}{\line(1,0){0.13}}
\multiput(60.33,30.98)(0.13,0.14){5}{\line(0,1){0.14}}
\multiput(60.96,31.67)(0.12,0.14){5}{\line(0,1){0.14}}
\multiput(61.57,32.38)(0.11,0.15){5}{\line(0,1){0.15}}
\multiput(62.13,33.11)(0.13,0.19){4}{\line(0,1){0.19}}
\multiput(62.66,33.87)(0.12,0.19){4}{\line(0,1){0.19}}
\multiput(63.16,34.64)(0.11,0.20){4}{\line(0,1){0.20}}
\multiput(63.61,35.43)(0.14,0.27){3}{\line(0,1){0.27}}
\multiput(64.03,36.24)(0.13,0.28){3}{\line(0,1){0.28}}
\multiput(64.41,37.07)(0.11,0.28){3}{\line(0,1){0.28}}
\multiput(64.75,37.91)(0.15,0.43){2}{\line(0,1){0.43}}
\multiput(65.05,38.76)(0.13,0.43){2}{\line(0,1){0.43}}
\multiput(65.30,39.63)(0.11,0.44){2}{\line(0,1){0.44}}
\multiput(65.51,40.50)(0.17,0.88){1}{\line(0,1){0.88}}
\multiput(65.69,41.38)(0.13,0.89){1}{\line(0,1){0.89}}
\multiput(65.82,42.27)(0.09,0.89){1}{\line(0,1){0.89}}
\multiput(65.90,43.16)(0.04,0.89){1}{\line(0,1){0.89}}

\linethickness{0.20mm}
\put(45.78,57.35){\line(0,1){0.90}}
\multiput(45.61,59.13)(0.17,-0.88){1}{\line(0,-1){0.88}}
\multiput(45.27,59.98)(0.11,-0.28){3}{\line(0,-1){0.28}}
\multiput(44.77,60.75)(0.12,-0.19){4}{\line(0,-1){0.19}}
\multiput(44.13,61.44)(0.13,-0.14){5}{\line(0,-1){0.14}}
\multiput(43.37,62.02)(0.15,-0.12){5}{\line(1,0){0.15}}
\multiput(42.51,62.47)(0.22,-0.11){4}{\line(1,0){0.22}}
\multiput(41.58,62.77)(0.31,-0.10){3}{\line(1,0){0.31}}
\multiput(40.60,62.93)(0.98,-0.16){1}{\line(1,0){0.98}}
\put(39.60,62.93){\line(1,0){0.99}}
\multiput(38.62,62.77)(0.98,0.16){1}{\line(1,0){0.98}}
\multiput(37.69,62.47)(0.31,0.10){3}{\line(1,0){0.31}}
\multiput(36.83,62.02)(0.22,0.11){4}{\line(1,0){0.22}}
\multiput(36.07,61.44)(0.15,0.12){5}{\line(1,0){0.15}}
\multiput(35.43,60.75)(0.13,0.14){5}{\line(0,1){0.14}}
\multiput(34.93,59.98)(0.12,0.19){4}{\line(0,1){0.19}}
\multiput(34.59,59.13)(0.11,0.28){3}{\line(0,1){0.28}}
\multiput(34.42,58.25)(0.17,0.88){1}{\line(0,1){0.88}}
\put(34.42,57.35){\line(0,1){0.90}}
\multiput(34.42,57.35)(0.17,-0.88){1}{\line(0,-1){0.88}}
\multiput(34.59,56.47)(0.11,-0.28){3}{\line(0,-1){0.28}}
\multiput(34.93,55.62)(0.12,-0.19){4}{\line(0,-1){0.19}}
\multiput(35.43,54.85)(0.13,-0.14){5}{\line(0,-1){0.14}}
\multiput(36.07,54.16)(0.15,-0.12){5}{\line(1,0){0.15}}
\multiput(36.83,53.58)(0.22,-0.11){4}{\line(1,0){0.22}}
\multiput(37.69,53.13)(0.31,-0.10){3}{\line(1,0){0.31}}
\multiput(38.62,52.83)(0.98,-0.16){1}{\line(1,0){0.98}}
\put(39.60,52.67){\line(1,0){0.99}}
\multiput(40.60,52.67)(0.98,0.16){1}{\line(1,0){0.98}}
\multiput(41.58,52.83)(0.31,0.10){3}{\line(1,0){0.31}}
\multiput(42.51,53.13)(0.22,0.11){4}{\line(1,0){0.22}}
\multiput(43.37,53.58)(0.15,0.12){5}{\line(1,0){0.15}}
\multiput(44.13,54.16)(0.13,0.14){5}{\line(0,1){0.14}}
\multiput(44.77,54.85)(0.12,0.19){4}{\line(0,1){0.19}}
\multiput(45.27,55.62)(0.11,0.28){3}{\line(0,1){0.28}}
\multiput(45.61,56.47)(0.17,0.88){1}{\line(0,1){0.88}}

\linethickness{0.20mm}
\put(37.88,36.26){\line(0,1){0.98}}
\multiput(37.70,38.21)(0.18,-0.97){1}{\line(0,-1){0.97}}
\multiput(37.36,39.14)(0.12,-0.31){3}{\line(0,-1){0.31}}
\multiput(36.85,39.99)(0.13,-0.21){4}{\line(0,-1){0.21}}
\multiput(36.20,40.75)(0.13,-0.15){5}{\line(0,-1){0.15}}
\multiput(35.43,41.38)(0.15,-0.13){5}{\line(1,0){0.15}}
\multiput(34.55,41.87)(0.22,-0.12){4}{\line(1,0){0.22}}
\multiput(33.60,42.21)(0.32,-0.11){3}{\line(1,0){0.32}}
\multiput(32.61,42.38)(1.00,-0.17){1}{\line(1,0){1.00}}
\put(31.59,42.38){\line(1,0){1.01}}
\multiput(30.60,42.21)(1.00,0.17){1}{\line(1,0){1.00}}
\multiput(29.65,41.87)(0.32,0.11){3}{\line(1,0){0.32}}
\multiput(28.77,41.38)(0.22,0.12){4}{\line(1,0){0.22}}
\multiput(28.00,40.75)(0.15,0.13){5}{\line(1,0){0.15}}
\multiput(27.35,39.99)(0.13,0.15){5}{\line(0,1){0.15}}
\multiput(26.84,39.14)(0.13,0.21){4}{\line(0,1){0.21}}
\multiput(26.50,38.21)(0.12,0.31){3}{\line(0,1){0.31}}
\multiput(26.32,37.24)(0.18,0.97){1}{\line(0,1){0.97}}
\put(26.32,36.26){\line(0,1){0.98}}
\multiput(26.32,36.26)(0.18,-0.97){1}{\line(0,-1){0.97}}
\multiput(26.50,35.29)(0.12,-0.31){3}{\line(0,-1){0.31}}
\multiput(26.84,34.36)(0.13,-0.21){4}{\line(0,-1){0.21}}
\multiput(27.35,33.51)(0.13,-0.15){5}{\line(0,-1){0.15}}
\multiput(28.00,32.75)(0.15,-0.13){5}{\line(1,0){0.15}}
\multiput(28.77,32.12)(0.22,-0.12){4}{\line(1,0){0.22}}
\multiput(29.65,31.63)(0.32,-0.11){3}{\line(1,0){0.32}}
\multiput(30.60,31.29)(1.00,-0.17){1}{\line(1,0){1.00}}
\put(31.59,31.12){\line(1,0){1.01}}
\multiput(32.61,31.12)(1.00,0.17){1}{\line(1,0){1.00}}
\multiput(33.60,31.29)(0.32,0.11){3}{\line(1,0){0.32}}
\multiput(34.55,31.63)(0.22,0.12){4}{\line(1,0){0.22}}
\multiput(35.43,32.12)(0.15,0.13){5}{\line(1,0){0.15}}
\multiput(36.20,32.75)(0.13,0.15){5}{\line(0,1){0.15}}
\multiput(36.85,33.51)(0.13,0.21){4}{\line(0,1){0.21}}
\multiput(37.36,34.36)(0.12,0.31){3}{\line(0,1){0.31}}
\multiput(37.70,35.29)(0.18,0.97){1}{\line(0,1){0.97}}

\linethickness{0.20mm}
\put(61.03,34.85){\line(0,1){0.89}}
\multiput(60.84,36.62)(0.09,-0.44){2}{\line(0,-1){0.44}}
\multiput(60.47,37.44)(0.12,-0.27){3}{\line(0,-1){0.27}}
\multiput(59.93,38.19)(0.13,-0.19){4}{\line(0,-1){0.19}}
\multiput(59.24,38.82)(0.14,-0.13){5}{\line(1,0){0.14}}
\multiput(58.43,39.31)(0.20,-0.12){4}{\line(1,0){0.20}}
\multiput(57.54,39.65)(0.30,-0.11){3}{\line(1,0){0.30}}
\multiput(56.59,39.83)(0.95,-0.17){1}{\line(1,0){0.95}}
\put(55.61,39.83){\line(1,0){0.97}}
\multiput(54.66,39.65)(0.95,0.17){1}{\line(1,0){0.95}}
\multiput(53.77,39.31)(0.30,0.11){3}{\line(1,0){0.30}}
\multiput(52.96,38.82)(0.20,0.12){4}{\line(1,0){0.20}}
\multiput(52.27,38.19)(0.14,0.13){5}{\line(1,0){0.14}}
\multiput(51.73,37.44)(0.13,0.19){4}{\line(0,1){0.19}}
\multiput(51.36,36.62)(0.12,0.27){3}{\line(0,1){0.27}}
\multiput(51.17,35.75)(0.09,0.44){2}{\line(0,1){0.44}}
\put(51.17,34.85){\line(0,1){0.89}}
\multiput(51.17,34.85)(0.09,-0.44){2}{\line(0,-1){0.44}}
\multiput(51.36,33.98)(0.12,-0.27){3}{\line(0,-1){0.27}}
\multiput(51.73,33.16)(0.13,-0.19){4}{\line(0,-1){0.19}}
\multiput(52.27,32.41)(0.14,-0.13){5}{\line(1,0){0.14}}
\multiput(52.96,31.78)(0.20,-0.12){4}{\line(1,0){0.20}}
\multiput(53.77,31.29)(0.30,-0.11){3}{\line(1,0){0.30}}
\multiput(54.66,30.95)(0.95,-0.17){1}{\line(1,0){0.95}}
\put(55.61,30.77){\line(1,0){0.97}}
\multiput(56.59,30.77)(0.95,0.17){1}{\line(1,0){0.95}}
\multiput(57.54,30.95)(0.30,0.11){3}{\line(1,0){0.30}}
\multiput(58.43,31.29)(0.20,0.12){4}{\line(1,0){0.20}}
\multiput(59.24,31.78)(0.14,0.13){5}{\line(1,0){0.14}}
\multiput(59.93,32.41)(0.13,0.19){4}{\line(0,1){0.19}}
\multiput(60.47,33.16)(0.12,0.27){3}{\line(0,1){0.27}}
\multiput(60.84,33.98)(0.09,0.44){2}{\line(0,1){0.44}}

\put(39.60,57.50){\makebox(0,0)[cc]{$D_1^j$}}

\put(31.10,37.30){\makebox(0,0)[cc]{$D_2^j$}}

\put(55.40,35.00){\makebox(0,0)[cc]{$D_3^j$}}

\put(48.20,53.80){\makebox(0,0)[cc]{$x$}}

\linethickness{0.20mm}
\multiput(46.90,51.90)(0.13,-0.12){10}{\line(1,0){0.13}}

\linethickness{0.20mm}
\multiput(46.90,50.80)(0.13,0.12){9}{\line(1,0){0.13}}

\put(78.20,62.50){\makebox(0,0)[cc]{$K$}}

\end{picture}